\magnification=1200
%\nopagenumbers

\catcode`\ =\active \def { }
%         ^u            ^u^o 
% u = unbreakable space,  o = ordinary space

\hsize=11.25cm    
\vsize=18cm     
\parindent=12pt   \parskip=5pt     

\hoffset=.5cm   
\voffset=.8cm   

\pretolerance=500 \tolerance=1000  \brokenpenalty=5000

\catcode`\@=11

\font\eightrm=cmr8         \font\eighti=cmmi8
\font\eightsy=cmsy8        \font\eightbf=cmbx8
\font\eighttt=cmtt8        \font\eightit=cmti8
\font\eightsl=cmsl8        \font\sixrm=cmr6
\font\sixi=cmmi6           \font\sixsy=cmsy6
\font\sixbf=cmbx6

\font\tengoth=eufm10 
\font\eightgoth=eufm8  
\font\sevengoth=eufm7      
\font\sixgoth=eufm6        \font\fivegoth=eufm5

\skewchar\eighti='177 \skewchar\sixi='177
\skewchar\eightsy='60 \skewchar\sixsy='60

\newfam\gothfam           \newfam\bboardfam

\def\tenpoint{
  \textfont0=\tenrm \scriptfont0=\sevenrm \scriptscriptfont0=\fiverm
  \def\rm{\fam\z@\tenrm}
  \textfont1=\teni  \scriptfont1=\seveni  \scriptscriptfont1=\fivei
  \def\oldstyle{\fam\@ne\teni}\let\old=\oldstyle
  \textfont2=\tensy \scriptfont2=\sevensy \scriptscriptfont2=\fivesy
  \textfont\gothfam=\tengoth \scriptfont\gothfam=\sevengoth
  \scriptscriptfont\gothfam=\fivegoth
  \def\goth{\fam\gothfam\tengoth}
  
  \textfont\itfam=\tenit
  \def\it{\fam\itfam\tenit}
  \textfont\slfam=\tensl
  \def\sl{\fam\slfam\tensl}
  \textfont\bffam=\tenbf \scriptfont\bffam=\sevenbf
  \scriptscriptfont\bffam=\fivebf
  \def\bf{\fam\bffam\tenbf}
  \textfont\ttfam=\tentt
  \def\tt{\fam\ttfam\tentt}
  \abovedisplayskip=12pt plus 3pt minus 9pt
  \belowdisplayskip=\abovedisplayskip
  \abovedisplayshortskip=0pt plus 3pt
  \belowdisplayshortskip=4pt plus 3pt 
  \smallskipamount=3pt plus 1pt minus 1pt
  \medskipamount=6pt plus 2pt minus 2pt
  \bigskipamount=12pt plus 4pt minus 4pt
  \normalbaselineskip=12pt
  \setbox\strutbox=\hbox{\vrule height8.5pt depth3.5pt width0pt}
  \let\bigf@nt=\tenrm       \let\smallf@nt=\sevenrm
  \normalbaselines\rm}

\def\eightpoint{
  \textfont0=\eightrm \scriptfont0=\sixrm \scriptscriptfont0=\fiverm
  \def\rm{\fam\z@\eightrm}
  \textfont1=\eighti  \scriptfont1=\sixi  \scriptscriptfont1=\fivei
  \def\oldstyle{\fam\@ne\eighti}\let\old=\oldstyle
  \textfont2=\eightsy \scriptfont2=\sixsy \scriptscriptfont2=\fivesy
  \textfont\gothfam=\eightgoth \scriptfont\gothfam=\sixgoth
  \scriptscriptfont\gothfam=\fivegoth
  \def\goth{\fam\gothfam\eightgoth}
  
  \textfont\itfam=\eightit
  \def\it{\fam\itfam\eightit}
  \textfont\slfam=\eightsl
  \def\sl{\fam\slfam\eightsl}
  \textfont\bffam=\eightbf \scriptfont\bffam=\sixbf
  \scriptscriptfont\bffam=\fivebf
  \def\bf{\fam\bffam\eightbf}
  \textfont\ttfam=\eighttt
  \def\tt{\fam\ttfam\eighttt}
  \abovedisplayskip=9pt plus 3pt minus 9pt
  \belowdisplayskip=\abovedisplayskip
  \abovedisplayshortskip=0pt plus 3pt
  \belowdisplayshortskip=3pt plus 3pt 
  \smallskipamount=2pt plus 1pt minus 1pt
  \medskipamount=4pt plus 2pt minus 1pt
  \bigskipamount=9pt plus 3pt minus 3pt
  \normalbaselineskip=9pt
  \setbox\strutbox=\hbox{\vrule height7pt depth2pt width0pt}
  \let\bigf@nt=\eightrm     \let\smallf@nt=\sixrm
  \normalbaselines\rm}

\tenpoint

\def\pc#1{\bigf@nt#1\smallf@nt}         \def\pd#1 {{\pc#1} }

\catcode`\;=\active
\def;{\relax\ifhmode\ifdim\lastskip>\z@\unskip\fi
\kern\fontdimen2  -1.2 \fontdimen3 \string;}

\catcode`\:=\active
\def:{\relax\ifhmode\ifdim\lastskip>\z@\unskip\fi\penalty\@M\ \fi\string:}

\catcode`\!=\active
\def!{\relax\ifhmode\ifdim\lastskip>\z@
\unskip\fi\kern\fontdimen2  -1.1 \fontdimen3 \string!}

\catcode`\?=\active
\def?{\relax\ifhmode\ifdim\lastskip>\z@
\unskip\fi\kern\fontdimen2  -1.1 \fontdimen3 \string?}

\catcode`\«=\active 
\def«{\raise.4ex\hbox{%
 $\scriptscriptstyle\langle\!\langle$}}

\catcode`\»=\active 
\def»{\raise.4ex\hbox{%
 $\scriptscriptstyle\rangle\!\rangle$}}

\frenchspacing

\def\raggedbottom{\topskip 10pt plus 36pt\r@ggedbottomtrue}

\def\pointir{\unskip . --- \ignorespaces}

\def\Medbreak{\vskip-\lastskip\medbreak}

\long\def\th#1 #2\enonce#3\endth{
   \Medbreak\noindent
   {\pc#1} {#2\unskip}\pointir{\it #3}\smallskip}

\def\decale#1{\smallbreak\hskip 28pt\llap{#1}\kern 5pt}
\def\decaledecale#1{\smallbreak\hskip 34pt\llap{#1}\kern 5pt}
\def\puce{\smallbreak\hskip 6pt{$\scriptstyle\bullet$}\kern 5pt}

\def\eqalign#1{\null\,\vcenter{\openup\jot\m@th\ialign{
\strut\hfil$\displaystyle{##}$&$\displaystyle{{}##}$\hfil
&&\quad\strut\hfil$\displaystyle{##}$&$\displaystyle{{}##}$\hfil
\crcr#1\crcr}}\,}

\catcode`\@=12

\showboxbreadth=-1  \showboxdepth=-1

\newcount\numerodesection \numerodesection=1
\def\section#1{\bigbreak
 {\bf\number\numerodesection.\ \ #1}\nobreak\medskip
 \advance\numerodesection by1}

\mathcode`A="7041 \mathcode`B="7042 \mathcode`C="7043 \mathcode`D="7044
\mathcode`E="7045 \mathcode`F="7046 \mathcode`G="7047 \mathcode`H="7048
\mathcode`I="7049 \mathcode`J="704A \mathcode`K="704B \mathcode`L="704C
\mathcode`M="704D \mathcode`N="704E \mathcode`O="704F \mathcode`P="7050
\mathcode`Q="7051 \mathcode`R="7052 \mathcode`S="7053 \mathcode`T="7054
\mathcode`U="7055 \mathcode`V="7056 \mathcode`W="7057 \mathcode`X="7058
\mathcode`Y="7059 \mathcode`Z="705A

% handling accented characters in plain TeX :

\catcode`\À=\active \defÀ{\`A}    \catcode`\à=\active \defà{\`a} 
\catcode`\Â=\active \defÂ{\^A}    \catcode`\â=\active \defâ{\^a} 
\catcode`\Æ=\active \defÆ{\AE}    \catcode`\æ=\active \defæ{\ae}
\catcode`\Ç=\active \defÇ{\c C}   \catcode`\ç=\active \defç{\c c}
\catcode`\È=\active \defÈ{\`E}    \catcode`\è=\active \defè{\`e} 
\catcode`\É=\active \defÉ{\'E}    \catcode`\é=\active \defé{\'e} 
\catcode`\Ê=\active \defÊ{\^E}    \catcode`\ê=\active \defê{\^e} 
\catcode`\Ë=\active \defË{\"E}    \catcode`\ë=\active \defë{\"e} 
\catcode`\Î=\active \defÎ{\^I}    \catcode`\î=\active \defî{\^\i}
\catcode`\Ï=\active \defÏ{\"I}    \catcode`\ï=\active \defï{\"\i}
\catcode`\Ô=\active \defÔ{\^O}    \catcode`\ô=\active \defô{\^o} 
\catcode`\Ù=\active \defÙ{\`U}    \catcode`\ù=\active \defù{\`u} 
\catcode`\Û=\active \defÛ{\^U}    \catcode`\û=\active \defû{\^u} 
\catcode`\Ü=\active \defÜ{\"U}    \catcode`\ü=\active \defü{\"u} 

\def\diagram#1{\def\normalbaselines{\baselineskip=0pt\lineskip=5pt}
\matrix{#1}}

\def\vfl#1#2#3{\llap{$\textstyle #1$}
\left\downarrow\vbox to#3{}\right.\rlap{$\textstyle #2$}}

\def\hfl#1#2#3{\smash{\mathop{\hbox to#3{\rightarrowfill}}\limits
^{\textstyle#1}_{\textstyle#2}}}

\def\ogoth{{\goth o}}

\def\pgoth{{\goth p}}

\def\Q{{\bf Q}}
\def\Qp{\Q_p}
\def\R{{\bf R}}
\def\C{{\bf C}}
\def\N{{\bf N}}

\def\T{{\bf T}}

\def\Z{{\bf Z}}
\def\Zp{\Z_p}
\def\F{{\bf F}}
\def\Fp{\F_p}

\def\Hom{\mathop{\rm Hom}\nolimits}
\def\Ann{\mathop{\rm Ann}\nolimits}
\def\Id{\mathop{\rm Id}\nolimits}

\def\Card{\mathop{\rm Card}\nolimits}
\def\Gal{\mathop{\rm Gal}\nolimits}
\def\Ker{\mathop{\rm Ker}\nolimits}

\def\series#1{(\!(#1)\!)}

\def\to{\rightarrow}

\def\droite#1{\,\hfl{#1}{}{8mm}\,}

\def\mod{\mathop{\rm mod.}\nolimits}
\def\pmod#1{\;(\mod#1)}
\def\nomod#1{\;(#1)}

\def\sup{\mathop{\rm Sup}\displaylimits}
\def\inf{\mathop{\rm Inf}\displaylimits}
\def\Aut{\mathop{\rm Aut}\nolimits}
\def\Autord{\mathop{\rm Autord}\nolimits}
\def\res{\mathop{\rm res}\nolimits}
\def\dis{\mathop{\rm dis}\nolimits}

\centerline{\bf A first course in Local arithmetic}
\medskip
\centerline{\it Notes for a course at the H.-C. R. I.}
\medskip
\medskip
\centerline{\it Allahabad, 15 August 2008 -- 26 January 2009}
\bigskip
\centerline{Chandan Singh Dalawat}

\vfill

{\eightpoint Seit meiner ersten Besch{\"a}ftigung mit den Fragen der
  h{\"o}heren Zahlentheorie glaubte ich, da\ss\ die Methoden der
  Funktionentheorie auch auf dieses Gebiet anwendbar sein m{\"u}ssten, und
  da\ss\ sich auf dieser Grundlage eine in mancher Hinsicht einfachere Theorie
  der algebraischen Zahlen aufbauen lassen k{\"o}nnte.\hfill\break  
\rightline{--- Kurt Hensel, {\it Theorie der algebraischen Zahlen}, 1908,
  p.~iv.}} 

\vfill\eject

\centerline{\bf Lecture 1}
\medskip
\centerline{\it Valuations and absolute values}
\bigskip

Let $k$ be a field and denote its multiplicative group by $k^\times$.  The
multiplicative group of strictly positive reals is denoted
$\R^{\times\circ}$~; it is a totally ordered group.

\th DEFINITION 1
\enonce
A modulus\/ $|\ |$ on\/ $k$ is a homomorphism\/ $|\ 
|:k^\times\rightarrow\R^{\times\circ}$ for which there exists a constant $C>0$
such that 
$$
|x+y|\le C.\sup(|x|,|y|)
$$ 
for every\/ $x,y\in k$, with the convention that\/ $|0|=0$.  A modulus is
called trivial if the image of\/ $k^\times$ is\/ $\{1\}$, essential
otherwise~; an essential modulus is called an absolute value if the triangular
inequality
$$
|x+y|\le |x|+|y|
$$ 
holds for every\/ $x,y\in k$.  An absolute value is called ultrametric if we
may take\/ $C=1$, archimedean otherwise.  A field endowed with an absolute
value is called a valued field.
\endth

{\eightpoint The term {\it modulus\/} is not standard.  We shall abandon it
  soon, after showing that every essential modulus is equivalent to an
  absolute value.  To avoid circumlocution, the trivial modulus is not an
  absolute value, by definition.}

If $\zeta\in k^\times$ has finite order, then $|\zeta|=1$, for $1$ is the only
element of finite order in $\R^{\times\circ}$.  Consequently, every modulus on
a finite field is trivial.

Recall that the group $\Autord(\R^{\times\circ})$ of order-preserving
automorphisms of $\R^{\times\circ}$ is isomorphic to $\R^{\times\circ}$ by the
map $\gamma\mapsto(a\mapsto a^\gamma)$.  Via the exponential isomorphism
$\exp:\R\to\R^{\times\circ}$, this is equivalent to the fact that
$\R^\times\to\Aut\R$, $\gamma\mapsto(a\mapsto a\gamma)$ is an isomorphism.
Clearly, $a\mapsto a\gamma$ preserves the order if and only if $\gamma>0$.

Thus the group $\R^{\times\circ}$ acts on the set of moduli on the field $k$.
Two moduli $|\ |_1$, $|\ |_2$ are called equivalent if they are in the same
orbit~; concretely, they are equivalent if there exists a real $\gamma>0$ such
that $|\ |_1=|\ |_2^\gamma$.

\th PROPOSITION 2
\enonce
An essential modulus $|\ |$ is an absolute value if and only if we may take\/
$C=2$. 
\endth
If $|\ |$ is an absolute value, then $|x+y|\le|x|+|y|\le2\sup(|x|,|y|)$, so we
may take $C=2$.  To prove the converse, suppose that we may take $C=2$.  By
induction, we get  
$$
|x_1+x_2+\cdots+x_{2^n}|\le 2^n\sup_{j\in[1,2^n]}|x_j|
$$
For an integer $N>0$, by taking $n$ such that $N\in\,]2^{n-1},2^n]$, and
setting $x_{N+1}=\cdots=x_{2^n}=0$ in the above inequality, we get
$$
|x_1+x_2+\cdots+x_{N}|\le 2^n\sup_{j\in[1,2^n]}|x_j|
\le 2N\sup_{j\in[1,N]}|x_j|.
$$
In particular, taking $x_j=1$ for every $j\in[1,N]$, we get $|N|\le2N$ for
every integer $N>0$.  Now let $x,y\in k$ and let $s>0$  be an integer.  The
binomial theorem and the above estimates give  
$$\eqalign{
|x+y|^s&=\left|\sum_{r\in[0,s]}{s\choose r}x^ry^{s-r}\right|\cr
&\le2(s+1)\sup_{r\in[0,s]}\left|{s\choose r}x^ry^{s-r}\right|\cr
&\le2(s+1)\sup_{r\in[0,s]}
\left|{s\choose r}\right|\left|x\right|^r\left|y\right|^{s-r}\cr
&\le4(s+1)\sup_{r\in[0,s]} {s\choose r}\left|x\right|^r\left|y\right|^{s-r}\cr
&\le4(s+1)\sum_{r\in[0,s]} {s\choose r}\left|x\right|^r\left|y\right|^{s-r}\cr
&=4(s+1)\left(|x|+|y|\right)^s.\cr
}$$
Taking the $s$-th root and letting $s\to+\infty$, we get the triangular
inequality $|x+y|\le |x|+|y|$, as required.

\th COROLLARY 3
\enonce
Every  essential modulus is equivalent to an absolute value.
\endth

On the field $\R$ of real numbers, we have the archimedean absolute value
$|x|_\infty=\sup(x,-x)$. 

\th LEMMA 4
\enonce
For every absolute value\/ $|\ |$ and every\/ $x,y\in k$, we have
$$
\left|\,|x|-|y|\,\right|_\infty\le|x-y|.
$$
\endth
As $y=x+(y-x)$, we have $|y|\le|x|+|y-x|$.  Similarly, $|x|\le|y|+|x-y|$.  But
$|y-x|=|x-y|$, hence the result.

\th LEMMA 5
\enonce
Let\/ $|\ |$ be an ultrametric absolute value, and\/ $x,y\in k$ such that\/
$|x|<|y|$.  Then\/ $|x+y|=|y|$.
\endth
We have $|x+y|\le\sup(|x|,|y|)=|y|$.  On the other hand, $y=(x+y)+(-x)$, so
$|y|\le\sup(|x+y|,|x|)$, and the $\sup$ cannot be $|x|$ by hypothesis.

\th LEMMA 6
\enonce
An absolute value\/ $|\ |$ is ultrametric if and only if\/ $|\iota(n)|\le1$
for every\/ $n\in\Z$, where\/ $\iota(n)\in k$ is the image of\/ $n$.
\endth
If $|\ |$ is ultrametric, then $|\iota(n)|\le1$ by induction.  Conversely,
suppose that $|\iota(n)|\le1$ for every $n\in\Z$.  Let $x,y\in k$, and let
$s>0$ be an integer.  We have
$$
\eqalign{
|x+y|^s&=\left|\sum_{r\in[0,s]}{s\choose r}x^ry^{s-r}\right|\cr
&\le\sum_{r\in[0,s]}
 \left|{s\choose r}\right|\left|x\right|^r\left|y\right|^{s-r}\cr
&\le\sum_{r\in[0,s]} \left|x\right|^r\left|y\right|^{s-r}\cr
&\le(s+1)\sup(|x|,|y|)^s.
}$$
Taking the $s$-th root and letting $s\to+\infty$, we get
$|x+y|\le\sup(|x|,|y|)$, as required.

\th COROLLARY 7
\enonce
If the restriction of\/ $|\ |$ to a subfield is ultrametric, then so is\/ $|\
|$.  Every absolute value on a field of characteristic\/ $\neq0$ is\/
ultrametric. 
\endth
The second statement follows because that the restriction of\/ $|\ |$ to the
prime subfield is trivial.

\th DEFINITION 9
\enonce
A valuation\/ $v$ on $k$ is a homomorphism\ $v:k^\times\to\R$ such that
$$
v(x+y)\ge\inf(v(x),v(y))
$$
for every\/ $x,y\in k$, with the convention that\/ $v(0)=+\infty$.  A
valuation\/ $v$ is said to be trivial if\/ $v(k^\times)=\{0\}$, of height\/~$1$
otherwise.  A height-$1$ valuation\/ $v$ is called discrete if the subgroup\/
$v(k^\times)\subset\R$ is discrete. A discrete valuation\/ $v$ is said to be 
normalised if\/ $v(k^\times)=\Z$.  
\endth

{\eightpoint Our terminology follows Bourbaki~; the German word is {\it
    (Exponential)bewertung}.  Some authors use {\it valuation\/} to mean what
  we have been calling a modulus (def.~1).}

The group $\R^{\times\circ}$ of order-preserving automorphisms of $\R$ acts on
the set of valuations on the field $k$.  Two valuations $v_1$, $v_2$ are
called equivalent if they are in the same orbit~; concretely, they are
equivalent if there exists a real $\gamma>0$ such that $v_1=\gamma v_2$.

{\it Example\/}~8\pointir Let $p$ be a prime number and let
$v_p:\Q^\times\to\R$ be the unique homomorphism such that $v_p(p)=1$ and
$v_p(l)=0$ for every prime $l\neq p$.  Then $v_p$ is a normalised discrete
valuation, called the $p$-adic valuation of $\Q$.  If $l$ is a prime $\neq p$,
then $v_l$ is not equivalent to $v_p$.

\vfill\eject
\centerline{\bf Lecture 2}
\bigskip
\centerline{\it Absolute values on $\Q$ and on $k(T)$}
\bigskip

We say that $a\in\Q$ is a $p$-adic integer if $v_p(a)\ge0$.  They form the
subring $\Z_{(p)}=v_p^{-1}([0,+\infty])$, and $\Z=\cap_p\Z_{(p)}$.  Clearly,
$\Z_{(p)}$ is the smallest subring containing $l^{-1}$ for every prime $l\neq
p$.  The ideal $p\Z_{(p)}$ is prime, and the map $\Z/p\Z\to\Z_{(p)}/p\Z_{(p)}$
is an isomorphism of fields.

If $v$ is a height-$1$ valuation on $k$, then $|x|_v=\exp(-v(x))$ is an
ultrametric absolute value, and conversely, if $|\ |$ is an ultrametric
absolute value, then $v_{|\ |}(x)=-\log|x|$ is a height-$1$ valuation.
Equivalent valuations correspond to equivalent (ultrametric) absolute values.

Notice that if $v(x)<v(y)$, then $v(x+y)=v(x)$ (lemma~5).

For every prime $p$, the absolute value $|x|_p=p^{-v_p(x)}$ is called the
$p$-adic absolute value on the field $\Q$.  A modulus $|\ |$ is equivalent to
$|\ |_p$ if and only if $|p|<1$ whereas $|l|=1$ for every prime $l\neq p$.

\th THEOREM 1 (Ostrowski, 1918)
\enonce
Let\/ $|\phantom{x}|$ be a modulus on the field\/ $\Q$ of rational numbers.
Then\/ $|\phantom{x}|$ is either trivial, or equivalent to the archimedean
absolute value\/ $|\phantom{x}|_\infty$, or equivalent to the $p$-adic
absolute value\/ $|\phantom{x}|_p\;$ for some prime\/ $p$.
\endth
{\it Proof\/} (Artin, 1932)~: 
Clearly $|\phantom{x}|$ is trivial if $|p|=1$ for every prime $p$.  Assume
that $|p|\neq1$ for some prime $p$. We shall show that if $|p|>1$, then
$|l|>1$ for every prime $l$ and that $|\phantom{x}|$ is equivalent to
$|\phantom{x}|_\infty$.  On the other hand, if $|p|<1$, then such a prime is
unique, and that $|\phantom{x}|$ is equivalent to $|\phantom{x}|_p$.

For the time being, let $p$ and $l$ be any two integers $>1$, and
write $p$ as 
$$
p=a_0+a_1l+a_2l^2+\cdots+a_nl^n%
$$
in base~$l$, with digits $a_i\in[0,l-1]$ and $l^n\le p$, i.e.,~$n\le\alpha$,
with $\displaystyle \alpha={\log p\over\log l}$.  We may further assume that
$|\phantom{x}|$ is an absolute value (cor.~3).  The triangular inequality
gives
$$\eqalign{
|p|&\le |a_0|+|a_1||l|+|a_2||l|^2+\cdots+|a_n||l|^n\cr
&\le(|a_0|+|a_1|+|a_2|+\cdots+|a_n|)\,\sup(1,|l|^n)\cr
&\le (1+n)d\,\sup(1,|l|^n)
\qquad(\hbox{with\ } d=\sup(|0|, |1|, \ldots, |l-1|)\cr
&\le (1+{\alpha})d\,\sup(1,|l|^{\alpha})
\qquad\left(\hbox{since\ } n\le{\alpha}\right).\cr
}
$$
Replace $p$ by $p^s$ and extract the $s$-th root to get
$$
|p|\le \left(1+s{\alpha}\right)^{1/s}d^{1/ s}
\sup\left(1,|l|^{\alpha}\right),
$$
so that we obtain the estimate
$
|p|\le\sup\left(1,|l|^{\alpha}\right)
$
upon letting $s\to+\infty$.

Now suppose that $p$ and $l$ are prime numbers and that $|p|>1$.  As
$1<|p|\le\sup\left(1,|l|^{\alpha}\right)$, we see that $|l|^{\alpha}>1$ and
hence $|l|>1$ and $|p|\le|l|^{\alpha}$. Interchanging the role of $p$ and $l$,
we deduce $|p|=|l|^{\alpha}$.

Defining $\gamma$ ($>0$) by the equation $|p|=|p|_\infty^\gamma$ for the fixed
prime $p$, we see that $|l|=|l|_\infty^\gamma$ for every prime $l$,
and hence $|x|=|x|_\infty^\gamma$ for every $x\in\Q^\times$,
i.e.,~$|\phantom{x}|$ is equivalent to the archimedean absolute value
$|\phantom{x}|_\infty$.

Finally, assume that $|p|<1$ for some prime $p$.  
% We shall show that such a prime is unique, and that $|\phantom{x}|$ is
% equivalent to $|\phantom{x}|_p$.
We have already seen that then $|l|\le1$ for every prime $l$, and
hence $|n|\le1$ for every $n\in\Z$.  Let us show that $|l|=1$ for
every prime $l\neq p$.  For every integer $s>0$, writing
$1=a_sp+b_sl^s$ ($a_s,b_s\in\Z$), we have
$$
1=|1|\le|a_s|\,|p|+|b_s|\,|l|^s\le|p|+|l|^s,
$$
i.e.,~$|l|^s\ge1-|p|$.  This is possible for all integers $s>0$ only if
$|l|=1$.  This shows that $|\phantom{x}|$ is equivalent to
$|\phantom{x}|_p$.

\th COROLLARY 2
\enonce
If the image of an absolute value\/ $|\ |:k^\times\to\R^{\times\circ}$ is
discrete, then\/  $|\ |$ is ultrametric.
\endth

If $k$ has characteristic~$\neq0$, then every absolute value is ultrametric
(Lecture~1, cor.~7).  If $k$ has characteristic~$0$, the restriction of $|\ |$
to $\Q$ has discrete image, so must be trivial or equivalent to $|\ |_p$ for
some prime $p$.  Hence $|\ |$ is ultrametric (Lecture~1, cor.~7).

% That none of the absolute values on $\Q$ can be neglected is illustrated by
% the following result --- almost a tautology --- in which $\bar P$ denoted the
% set $P$ of primes together with an additional element $\infty$~:

Let $\bar P$ denote the set of primes together with an additional element
$\infty$.

\th THEOREM 3 (``Product formula'')
\enonce
For every $x\in\Q^\times$, one has\/ $|x|_v=1$ for almost all\/
$v\in\bar P$ and\/ $\prod_{v\in\bar P}|x|_v=1$.
\endth

Indeed, it is sufficient to verify this for $x=-1$ and for $x=p$ ($p\in P$).
Note that for the product formula to hold, we have to normalise the absolute
values on $\Q$ suitably.

\centerline{\hbox to5cm{\hrulefill}}

Let $k$ be a field and put $K=k(T)$.  Recall that the group
$K^\times\!/k^\times$ is the free commutative group on the set $P_K$
of monic irreducible polynomials $f$ in $k[T]$.  For each $f\in P_K$,
let $v_f:K^\times\rightarrow\Z$ be the unique homomorphism which is
trivial on $k^\times$, sends $f$ to $1$, and sends every other element
of $P_K$ to $0$.  It can be checked that $v_f$ is a discrete
valuation.

Also, the map $v_\infty:K^\times\rightarrow\Z$ which sends $a$ to
$-\deg(a)$ is a discrete valuation, trivial on $k^\times$.

It is clear that the discrete valuations $v_\infty$, $v_f$ ($f\in
P_K$) are mutually inequivalent.

\th THEOREM 4
\enonce
Up to equivalence, the only height-$1$ valuations on\/ $k(T)$, trivial
on\/ $k$, are\/ $v_\infty$ and the\/ $v_f$, one for each\/ $f\in P_K$.
\endth
Let $v$ be a height-$1$ valuation on $k(T)$, trivial on $k$.  We will
show that if $v(T)<0$, then $v$ is equivalent to $v_\infty$, whereas
if $v(T)\ge0$, then there is a unique $f\in P_K$ with $v(f)>0$ and $v$
is equivalent to $v_f$.

Suppose that $v(T)<0$.  It is sufficient to show that
$v(f)=v(T)\deg(f)$ for every $f\in P_K$.  This is clearly true for
$f=T$.  For any other $f$, write
$f=T^{n}(1+\alpha_1T^{-1}+\cdots+\alpha_nT^{-n})$ (with $n=\deg(f)$ and
$\alpha_i\in k$, at least one of them $\neq0$).  As $v(1)=0$ and
$v(\alpha_1T^{-1}+\cdots+\alpha_nT^{-n})>0$, we have
$v(1+\alpha_1T^{-1}+\cdots+\alpha_nT^{-n})=0$, and, finally,
$v(f)=v(T)\deg(f)$, i.e.~$v$ is equivalent to $v_\infty$.

Suppose now that $v(T)\ge0$~; then $v(a)\ge0$ for all $a\in k[T]$.  If
we had $v(f)=0$ for every $f\in P_K$, the valuation $v$ would be
trivial, not of height~$1$.  Pick $p\in P_K$ for which $v(p)>0$.
For every $q\neq p$ in $P_K$, write $1=ap+bq$ ($a,b\in k[T]$).  We
have
$$
0=v(1)\ge\inf\left(v(a)+v(p),v(b)+v(q)\right)
\ge\inf\left(v(p),v(q)\right)\ge0,
$$
which is possible only if $v(q)=0$.  It follows that $v$ is equivalent
to $v_p$.

It is instructive to compare this proof with Artin's proof classifying
absolute values on $\Q$.

\th COROLLARY 5
\enonce
Up to equivalence, the only absolute values on\/ $k(T)$,
trivial on\/ $k$, are\/ $|\phantom{x}|_\infty$ and the\/
$|\phantom{x}|_f$, one for each\/ $f\in P_K$.
\endth

Put $\bar P_K=P_K\cup\{\infty\}$ and $\deg(\infty)=1$.  For each
$p\in\bar P_K$, define the absolute value $|x|_p=e^{-v_p(x)\deg(p)}$.

\th THEOREM 6 (``Product formula'')
\enonce
For every\/ $x\in k(T)^\times$, one has\/ $|x|_p=1$ for almost all\/
$p\in\bar P_K$ and\/ $\prod_{p\in\bar P_K}|x|_p=1$.
\endth
This is clearly true for $x\in k^\times$, so it is sufficient to check
this for $x\in P_K$.  We have $v_\infty(x)=-\deg(x)$, $v_x(x)=1$, and
$v_q(x)=0$ for every $q\neq x$ in $P_K$, which gives the ``sum
formula'' $\sum_{p\in\bar P_K}\deg(p)v_p(x)=0$.  The result follows
from this upon exponentiating.

Note that here too, as in the case of $\Q$ earlier, it is necessary to
normalise the absolute values suitably for the product formula to
hold.

\vfill

% {\eightpoint
% Das System der Funktionalbedigungen
% $$
% \varphi(x)\varphi(y)=\varphi(xy)\ \ (I)\qquad 
% \varphi(x+y)\le\varphi(x)+\varphi(y)\ \ (II)
% $$
% wo die Argumente $x$, $y$ s{\"a}mtliche rationale Zahlen durchlaufen, die
% Funktion $\varphi(x)$ aber nur reelle Werte annimmt, hat nur folgende
% L{\"o}sungen
% $$
% a)\ \varphi(x)\equiv0\;;\quad b)\ \varphi(x)\equiv1\;;\quad 
% c)\ \varphi(0)=0\;,\ \varphi(x)=1\ \ (x\neq1)\leqno{1)}
% $$
% \vskip-.5cm
% $$
% \varphi(x)\equiv x^{\phantom{\rho}}\leqno{2)}
% $$
% \vskip-.5cm
% $$
% \varphi(x)\equiv|x|^\rho\leqno{3)}
% $$
% wo der positive Parameter $\rho$ nicht gr{\"o}sser als 1, sonst aber beliebig
% ist~: 
% $$
% \varphi(x)\equiv c^{a(p,x)}\leqno{4)}
% $$
% wo der positive Parameter $c$ kleiner als $1$ ist, $a(p,x)$ aber die Ordnung
% von $x$ in Bezug auf eine beliebige Primzahl $p$ bezeichnet.  Jeder Wahl von
% $c$ und $p$ entspricht eine Funktion $\varphi(x)$ von Typus 4.
% \hfill--- Ostrowski, {\it Acta}, {\bf 41} (1918), pp.~276--77.
% }

\vfill\eject

\centerline{\bf Lecture 3}
\bigskip
\centerline{\it Denominators of bernoullian numbers}
\bigskip

We shall be concerned with the denominators of bernoullian numbers.  These
rational numbers $B_k$ are defined in terms of the exponential series
$e^T=\exp(T)$ by the identity
$$
{T\over e^T-1}=B_0+B_1{T^1\over1!}+\sum_{k>1}B_k{T^k\over k!},
$$
so that $B_0=1$, $B_1=-1/2$.  As the function $B_0+B_1T-T/(e^T-1)$ is
invariant under $T\mapsto-T$, we have $B_k=0$ for $k$ odd $>1$.  Here are the
numerators $N_k$ and the denominators $D_k$ of the first few bernoulliian
numbers $B_k=N_k/D_k$ of even index $k>0$~:
$$
\vbox{\halign{&\hfil$#$\quad\cr
\multispan{11}\hrulefill\cr
k=&2&4&6&8&10&12&14&16&18&20\cr 
\noalign{\vskip-5pt}
\multispan{11}\hrulefill\cr
N_k=&1&-1&1&-1&5&-691&7&-3617&43867&-174611\cr
D_k=&6&30&42&30&66&2730&6&510&798&330\cr
\multispan{11}\hrulefill.\cr}}
$$

\th THEOREM 1 (von Staudt--Clausen, 1840)  
\enonce 
Let $k>0$ be an even integer, and let $l$ run through the primes.  Then the
number 
$$
W_k=B_k+\sum_{l-1\,|\,k}{1\over l}\leqno{(1)}
$$
is always an integer.  For example, $\displaystyle
B_{12}+{1\over 2}+{1\over3}+{1\over5}+{1\over7}+{1\over 13}=1$.  
\endth
{\eightpoint The British analyst Hardy says in his {\it Twelve lectures\/}
  (p.~11) that this theorem was rediscovered by Ramanujan ``\thinspace at a
  time of his life when he had hardly formed any definite concept of
  proof\thinspace''.}

{\it Proof\/} (Witt)~: The idea is to show that $W_k$ is a $p$-adic integer
for every prime $p$.  More precisely, we show that $B_k+p^{-1}$ (resp.~$B_k$)
is a $p$-adic integer if $p-1\,|\,k$ (resp.~if not).

For an integer $n>0$, let $S_k(n)=1^k+2^k+\cdots+(n-1)^k$.  Comparing the
coefficients on the two sides of
$$
1+e^T+e^{2T}+\cdots+e^{(n-1)T}={e^{nT}-1\over T}{T\over e^T-1},
$$
we get $\displaystyle S_k(n) =\sum_{m\in[0,k]}{k\choose m}{B_m\over
  k+1-m}n^{k+1-m}$.  To recover $B_k$ from $S_k(n)$, it is tempting to take
the limit $\displaystyle\lim_{n\to0}S_k(n)/n$, which doesn't make sense in the
archimedean world.  If, however, we make $n$ run through the powers $p^s$ of a
fixed prime $p$, then, $p$-adically, $p^s\to0$ as $s\to+\infty$, and
$$
\lim_{r\to+\infty}S_k(p^r)/p^r=B_k.\leqno{(2)}
$$
Let us compare $S_k(p^{s+1})/p^{s+1}$ with $S_k(p^s)/p^s$.  Every
$j\in[0,p^{s+1}[$ can be uniquely written as $j=up^s+v$, where $u\in[0,p[$ and
$v\in[0,p^s[$.  Now,
$$
\eqalign{
S_k(p^{s+1})
&=\sum_{j\in[0,p^{s+1}[}j^k=\sum_{u\in[0,p[}\sum_{v\in[0,p^s[}(up^s+v)^k\cr
% &\equiv p\left(\sum_{v\in[0,p^s[}v^k\right)
%  +kp^s\left(\sum_{u\in[0,p[}u\sum_{v\in[0,p^s[}v^{k-1}\right)
% \pmod{p^{2s}}\cr  
&\equiv p\left(\sum_v v^k\right)
 +kp^s\left(\sum_u u\sum_v v^{k-1}\right)\pmod{p^{2s}}\cr 
}$$
by the binomial theorem.  As $\sum_{v}v^k=S_k(p^s)$ and
$2\sum_uu=p(p-1)\equiv0\pmod p$, we get
$$
S_k(p^{s+1})\equiv pS_k(p^s)\pmod{p^{s+1}},
$$
where, for $p=2$, the fact that $k$ is even has been used.  Dividing
throughout by $p^{s+1}$, this can be expressed by saying that
$$
{S_k(p^{s+1})\over p^{s+1}}-{S_k(p^s)\over p^s}
$$
is a $p$-adic integer, and, since $\Z_{(p)}$ is a subring of $\Q$,
$$
{S_k(p^r)\over p^r}-{S_k(p^s)\over p^s}\in\Z_{(p)}
$$
for any two integers $r>0$, $s>0$.  Fixing $s=1$ and letting $r\to+\infty$,
we see that $B_k-S_k(p)/p\in\Z_{(p)}$, in view of $(2)$. (If a sequence
$(x_n)_n$ of $p$-adic integers tends to $a\in\Q$ as $n\to+\infty$, then
$a\in\Z_{(p)}$.  Why~?  Because $v_p(x_n)\to v_p(a)$.)  Now,
$$
S_k(p)=\sum_{j\in[1,p[}j^k\equiv\cases
{-1\pmod p&if $p-1\,|\,k$\cr
\phantom{-}0\pmod p&otherwise\ (*)\cr}
$$
and hence $B_k+p^{-1}\in\Z_{(p)}$ if $p-1\,|\,k$ and $B_k\in\Z_{(p)}$
otherwise.  In either case, the number $W_k$ (1), which can be written as
$$
W_k=\cases{
(B_k+p^{-1})+\sum_{l\neq p}l^{-1}&if $p-1\,|\,k$\cr
(B_k)+\sum_l l^{-1}&otherwise,\cr
}$$ 
(where $l$ runs through the primes for which $l-1\,|\,k$) turns out to be a
$p$-adic integer for every prime $p$.  Hence $W_k\in\Z$, as claimed.

{\eightpoint (*) To see that $\sum_{j\in[1,p[}j^k\equiv0\pmod p$ when $p-1$
  does not divide $k$, note that, $g$ being a generator of $(\Z/p\Z)^\times$,
  we have $g^k-1\not\equiv0$, whereas
$$
(g^k-1)\left(\sum_{j\in[1,p[}j^k\right)
\equiv (g^k-1)\left(\sum_{t\in[0,p-1[}g^{tk}\right)
\equiv g^{(p-1)k}-1\equiv0.
$$
}

% {\eightpoint
% Wenn\/ $B^{(n)}$ die $n$te {\it Bernoulli\/}schen Zahl bezeichnet, und
%   $\alpha$, $\beta$, $\gamma$, $\ldots$, $\lambda$ diejenigen Theiler von $n$
%   sind, f{\"u}r welche $2\alpha+1$, $2\beta+1$, $2\gamma+1$, $\ldots$,
%   $2\lambda+1$ Primzahlen werden, so ist
% $$
% B^{(n)}\pm\left({1\over2}
% +{1\over 2\alpha+1}
% +{1\over 2\beta+1}
% +{1\over 2\gamma+1}+\cdots
% +{1\over 2\lambda+1}\right)
% =\hbox{einer ganzen Zahl,}
% $$
% wo $\pm$ gilt, wenn $n$ ungerade (gerade) ist. \hfill--- von Staudt, {\it
%   Crelle}, {\bf 41} (1840), pp.~373--4. 
% }

\vfill\eject

\centerline{\bf Lecture 4}
\bigskip
\centerline{\it Independence of absolute values}
\bigskip

Any two absolute values on the same field $k$ are either equivalent or
``independent'' in a certain sense.  A precise version is proved as th.~4, for
which we need some lemmas, of interest in their own right.

\th LEMMA 1
\enonce
Let\/ $|\ |_1$ and\/ $|\ |_2$ be absolute values on\/ $k$.  If\/ $|x|_1<1$
implies\/ $|x|_2<1$ for every\/ $x\in k$, then\/ $|\ |_1$ and\/ $|\ |_2$ are
equivalent.
\endth
Taking $x=y^{-1}$, we see that $|y|_1>1$ implies\/ $|y|_2>1$ for every\/ $y\in
k^\times$.

% Suppose, if possible, that there is a $z\in k^\times$ such that $|z|_1=1$
% but $|z|_2\neq1$.  Replacing $z$ by $z^{-1}$ if necessary, assume that
% $|z|_2>1$.  
% 
% As $|\ |_1$ is not the trivial modulus, there is a $w\in k^\times$ such that
% $|w|_1<1$.  Then, $|wz^n|_1=|w|_1|z|_1^n=|w|_1<1$ for every $n\in\N$, whereas 
% $$
% |wz^n|_2=|w|_2|z|_2^n\to+\infty\qquad(n\to+\infty),
% $$
% which contradicts the hypothesis ($|wz^n|_1<1\Rightarrow|wz^n|_2<1$).  The
% conclusion so far is that
% $$
% |x|_1<1\Leftrightarrow|x|_2<1,\qquad 
% |x|_1=1\Leftrightarrow|x|_2=1,\qquad
% |x|_1>1\Leftrightarrow|x|_2>1.
% $$
Fix $z\in k^\times$ such that $|z|_1>1$, and take
$\rho=\log|z|_2/\log|z|_1$, which is ${}>0$.  Let $y\in k^\times$ and put
$|y|_1=|z|_1^\gamma$.  For integer $m$, $n$ such that $n>0$ and $m/n>\gamma$,
we have $|y|_1<|z|_1^{m/n}$, therefore $|y^n/z^m|_1<1$, and hence
$|y^n/z^m|_2<1$ and $|y|_2<|z|_2^{m/n}$.  Similarly, if $m/n<\gamma$, then
$|y|_2>|z|_2^{m/n}$.  Hence $|y|_2=|z|_2^\gamma$.  

In other words, $\log|y|_2=\gamma\log|z|_2=\gamma\rho\log|z|_1=\rho\log|y|_1$,
so $|y|_1^\rho=|y|_2$, and the two absolute values are equivalent.

% Now let $y,z\in k^\times$ and apply this remark to $x=y^mz^n$ ($m,n\in\Z$).
% Taking logarithms, we see that 
% $$
% m\log |y|_1+n\log|z|_1\phantom{.}
% $$
% is ${}<0$, or ${}=0$, or ${}>0$ if and only if the same is true of
% $$
% m\log |y|_2+n\log|z|_2.
% $$
% Fixing $z$ such that $|z|_1\neq1$, and taking
% $\lambda=\log|z|_1/\log|z|_2$, which is ${}>0$, we conclude that
% $\log|y|_1=\lambda\log|y|_2$ for all $y\in k$.  Exponentiating, we get
% $$
% |y|_1=|y|_2^\lambda
% $$
% for every $y\in k$, which shows that $|\ |_1$ and\/ $|\ |_2$ are equivalent. 
% %  ce qu'il fallait d{\'e}montrer. 

Every absolute value $|\ |$ defines a distance $d(x,y)=|x-y|$, and hence a
topology on $k$.  The triangular inequality $d(x,z)\le d(x,y)+d(y,z)$ follows
from the triangular inequality for $|\ |$.

\th LEMMA 2
\enonce
Two absolute values are equivalent if and only if they define the same
topology.
\endth
Let $|\ |_1$ and\/ $|\ |_2$ be two absolute values which induce the same
topology.  Apply lemma~1 to the following equivalences
$$\eqalign{
|x|_1<1&\Leftrightarrow|x^n|_1\to0\hbox{ (as }n\to+\infty)\cr
&\Leftrightarrow x^n\to0\cr
&\Leftrightarrow|x^n|_2\to0\cr
&\Leftrightarrow|x|_2<1.
}
$$
The following lemma implies, and will be superceded by, the ``weak
approximation'' theorem (th.~4). 
\th LEMMA 3
\enonce
Let\/ $|\ |_1$, $|\ |_2$, $\ldots$, $|\ |_n$ be pairwise inequivalent
absolute values on\/ $k$.  There exists\/ $x\in k$ such
that
$$
|x|_1>1,\qquad |x|_j<1\quad(j\in [2,n]).
$$
\endth
{\it The case $n=2$.}  By lemma~1, there is some $y\in
k$ such that $|y|_1<1$ and $|y|_2\ge1$.  Similarly, there is some $z\in
k$ such that $|z|_2<1$ and $|z|_1\ge1$.  Then $x=zy^{-1}$ works.

{\it The case $n>2$.}  Proceed by induction.  By the inductive hypothesis,
there is a $y\in k$ such that $|y|_1>1$ and $|y|_j<1$ ($j\in[2,n[$).  For the
same reason, since $|\ |_1$ and $|\ |_n$ are inequivalent, there is a $z\in k$
such that $|z|_1>1$ and $|z|_n<1$.  For sufficiently large $r$, take
$$
x=\cases{
y&if $\,|y|_n<1$,\cr
zy^r&if  $\,|y|_n=1$,\cr
zy^r/(1+y^r)&if  $|\,y|_n>1$.\cr
}
$$
To see that $x=zy^r/(1+y^r)$ works in the last case, note that, as
$r\to+\infty$, 
$$
{y^r\over 1+y^r}={1\over1+y^{-r}}\to\cases{
1&for $|\ |_1$ and $|\ |_n$,\cr
0&for $|\ |_j$ ($j\neq1,n$).\cr
}$$

\th THEOREM 4 (``Weak approximation'', Artin-Whaples, 1945)  
\enonce
Let\/ $|\ |_1$, $|\ |_2$, $\ldots$, $|\ |_n$ be  pairwise inequivalent
absolute values on\/ $k$.  Given\/ $\varepsilon>0$ and\/ $x_j\in k$  for every
$j\in[1,n]$, there exists a\/ $y\in k$ such that
$$
|x_j-y|_j<\varepsilon\qquad (j\in [1,n]).
$$
\endth
By lemma~3, there are $z_j\in k$ such that $|z_j|_j>1$ and $|z_j|_i<1$ for
$i\neq j$. For  sufficiently large $r$, take
$$
y=\sum_{j\in[1,n]}{z_j^r\over 1+z_j^r}x_j.
$$

\th COROLLARY 5
\enonce
Let\/ $|\ |_1$, $|\ |_2$, $\ldots$, $|\ |_n$ be pairwise
inequivalent absolute values on\/ $k$.  The relation
$$
|x|_1^{\rho_1}|x|_2^{\rho_2}\cdots|x|_n^{\rho_1}=1
$$
holds for all\/ $x\in k^\times$ if and only if\/
$\rho_1=\rho_2=\cdots=\rho_n=0$. 
\endth
If $\rho_i>0$ for some $i\in[1,n]$, the given relation cannot hold for an
$x\in k^\times$ for which $|x|_i$ is very large whereas $|x|_j$ is close to
$1$ for every $j\neq i$.

Compare this result with the ``product formula'' for $\Q$ (Lecture~2,
th.~3). 

% \vfill\eject
\bigskip
\centerline{\hbox to5cm{\hrulefill}}
\bigbreak
\centerline{\it Completions}
\bigskip

Every valued field --- a field $k$ endowed with an absolute value $|\ |$ ---
can be completed to get a topological field $\hat k$ together with a
continuous extension of $|\ |$ to $\hat k$.  The primordial example is the way
$\R$ is obtained from $\Q$ endowed with its archimedean absolute value $|\ 
|_\infty$.

Recall that a sequence $(x_n)_n$ of elements of $k$ is said to tend to $y\in
k$ if for every $\varepsilon>0$, there is an index $N$ such that for every
$n>N$, we have $|x_n-y|<\varepsilon$.  If $x_n\to y$, then $y$ is unique and
called the limit of $(x_n)_n$.

The sequence $(x_n)_n$ is called a {\it fundamental sequence\/} if for every
$\varepsilon>0$, there is an index $N$ such that for every $m,n>N$, we have
$|x_m-x_n|<\varepsilon$.  If $(x_n)_n$ has a limit, then it is a fundamental
sequence by the triangular inequality. 

\th DEFINITION 6
\enonce
A valued field\/ $k$, $|\ |$ is said to be complete
if every fundamental sequence has a limit in\/ $k$.
\endth

The fields $\R$, $\C$ are complete with respect to $|\ |_\infty$.  It can be
shown that these are the only commutative fields complete with respect to an
archimedean absolute value (Ostrowski, 1918).

{\eightpoint The trivial modulus on $k$ satisfies the triangular inequality,
  hence defines a distance, for which $k$ is complete.  The induced topology
  on $k$ is discrete.  We shall see that the completion for an absolute value
  is never discrete.}

\th DEFINITION 7
\enonce
An extension\/ $\hat k$ of\/ $k$, endowed with an absolute value\/ $|\!|\
|\!|$ extending\/ $|\ |$, is said to be a completion of\/ $k$ if\/ $k$ is
dense in\/ $\hat k$ and\/ $\hat k$ is complete.
\endth

\th THEOREM 8
\enonce
Every valued field\/ $k$, $|\ |$ has a completion\/ $\hat k$, $|\!|\
|\!|$ which is unique up to a canonical isomorphism.
\endth

The uniqueness of the completion means that if $K_i$, $|\!|\ |\!|_i$ ($i=1,2$)
are two completions of $k$, then there is a unique $k$-isomorphism $f:K_1\to
K_2$ such that $|\!|f(x)|\!|_2=|\!|x|\!|_1$ for every $x\in K_1$.

Let $\hat k$ be the completion of the metric space $k$ (endowed with the
distance $d(x,y)=|x-y|$), let $\hat d$ be the distance function on $\hat k$,
and put $|\!|x|\!|=\hat d(x,0)$.  We shall make $\hat k$ into an extension of
$k$ (the basic reason being that the functions $x+y$, $xy$ and $x^{-1}$ are
continous on $k^2$, resp.~$k^\times$) and show that $|\!|\ |\!|$ is an
absolute value.

$(x,y)\mapsto x+y$ {\it is continuous}.  Because
$|(x_1+\delta_1)+(x_2+\delta_2) -(x_1+x_2)|$ is smaller than
$|\delta_1|+|\delta_2|$, by the triangular inequality.

$(x,y)\mapsto xy$ {\it is continuous}.  Because $|(x_1+\delta_1)(x_2+\delta_2)
-x_1x_2|$ is smaller than $|\delta_1\delta_2|+|x_1\delta_2|+|x_2\delta_1|$, by
the same inequality.

$x\mapsto x^{-1}$ {\it is continuous}.  Because $|(x+\delta)^{-1}-x^{-1}|$
equals $|\delta|(|x+\delta||x|)^{-1}$, which tends to~$0$ as $\delta\to0$.

By continuity, addition and multiplication (resp.~inversion) extend to $\hat
k^2$ (resp.~$\hat k^\times$) and give $\hat k$ the structure of an extension
of~$k$.  Finally, $|\!|\ |\!|$ is an absolute value on $\hat k$, again by
continuity.  The uniqueness of the completion is a consequence of the
uniqueness of the completion of a metric space.

The weak approximation theorem can be expressed more suggestively in terms of
completions. 

\th THEOREM 9 (``Weak approximation'')
\enonce
Let\/ $|\ |_j$ ($j\in [1,n]$) be a family of pairwise inequivalent
absolute values on\/ $k$, and\/ $\hat k_j$ the completion of\/ $k$ with
respect to\/ $|\ |_j$.  Then the image of the diagonal embedding\/
$k\to\prod_j\hat k_j$ is everywhere dense. 
\endth
Let $(x_1,\ldots,x_n)$ be a point of the product, and let $\varepsilon>0$ be
given.  There is a point $(y_1,\ldots,y_n)\in k^n$ such that
$|x_j-y_j|_j<\varepsilon$ for every $j\in[1,n]$.  By th.~4,
there is a $z\in k$ such that $|y_j-z|_j<\varepsilon$ for every $j\in[1,n]$.
Now $|x_j-z|_j<2\varepsilon$, by the triangular inequality.

{\eightpoint There is a ``strong approximation theorem'' for  ``global
  fields''.}

Let $k$ be a field endowed with a height-$1$ valuation $v:k^\times\to\R$, and
$\hat k$ the completion of $k$ with respect to the corresponding absolute
value.

\th PROPOSITION 10
\enonce
There is a unique extension $\hat v:\hat k^\times\to\R$ of\/ $v$ to\/ $\hat
k$, and\/ $\hat v(\hat k^\times)=v(k^\times)$.
\endth
The existence and uniqueness of $\hat v$ follow from the continuity of $v$, as
does the fact that $\hat v$ is a valuation.  Or simply take $\hat v=-\log|\!|\
|\!|$.  

To see that $\hat v$ and $v$ have the same image, observe that, by the density
of $k$ in $\hat k$, every $x\in\hat k$ can be written as $x=y+a$, where $y\in
k$ and $\hat v(a)$ is as large as we please.  Take $\hat v(a)>\hat v(x)$.
Then (cf.~Lecture~1, lemma~5),
$$
v(y)=\hat v(x-a)=\hat v(x).
$$

We shall say that $\hat k$, endowed with $\hat v$, is the {\it completion\/}
of $k$, $v$.

\bigskip
\centerline{\hbox to5cm{\hrulefill}}

\vfill\eject

\centerline{\bf Lecture 5}
\bigskip
\centerline{\it Complete valued fields}
\bigskip

Let $K$ be a field endowed with a height-$1$ valuation $v:K^\times\to\R$.
Then $\ogoth=v^{-1}([0,+\infty])$ is a subring of $K$, called the ring of
$v$-integers of $K$, and $\pgoth=v^{-1}(]0,+\infty])$ is the unique maximal
ideal of $\ogoth$, because $\ogoth^\times=v^{-1}(0)$.  The field
$k=\ogoth/\pgoth$ is called the residue field of $v$, or of $\ogoth$, or of
$K$, by an abuse of language.  Notice that $K$ is the field of fractions of
$\ogoth$. 

Let $\hat K$ be the completion of $K$, $\hat v$ its valuation, $\hat\ogoth$
the ring of $\hat v$-integers of $K$, $\hat\pgoth$ the unique maximal ideal of
$\hat\ogoth$, and $\hat k$ the residue field of $\hat v$.  We have
$$
\ogoth=\hat\ogoth\cap K,\qquad
\pgoth=\hat\pgoth\cap K,\qquad
$$
and consequently an embedding $k\to\hat k$ of fields.

\th LEMMA 1
\enonce
We have\/ $\hat v(\hat K^\times)=v(K^\times)$, and \/ $k\to\hat k$ is
an isomorphism. 
\endth
Let $x\in\hat\ogoth$ be given.  There is a $y\in K$ such that $\hat
v(x-y)>\hat v(x)$.  Then $v(y)=\hat v(x)$, so $y\in\ogoth$ and
$x-y\in\hat\pgoth$. 

\th COROLLARY 2
\enonce
$\hat\ogoth$ (resp.~$\hat\pgoth$) is the closure of\/ $\ogoth$
(resp.~$\pgoth$) in\/ $\hat K$.
\endth
{\eightpoint The topological closure is meant, not the integral closure.
  French avoids the confusion by having two different words, {\it
    adh{\'e}rence\/} and {\it fermeture}, apart from {\it cl{\^o}ture}.}

\th LEMMA 3
\enonce
The valuation\/ $v$ is discrete precisely when the ideal\/ $\pgoth$ is
principal. 
\endth

We shall show that if $\pi$ generates the ideal $\pgoth$, then $v(\pi)$
generates the group $v(K^\times)$.  Conversely, if the discrete subgroup
$v(K^\times)\subset\R$ is generated by $a>0$, we shall show that every $\pi\in
v^{-1}(a)$ generates the ideal $\pgoth$.

Suppose that $\pi$ generates $\pgoth$, and let $x\in K^\times$.  If $v(x)>0$,
then $x\in\pgoth$, and $x=y\pi$ for some $y\in\ogoth$.  Therefore $v(x)\ge
v(\pi)$.  Similarly, if $v(x)<0$, then, applying the argument to $x^{-1}$, we
get $v(x)\le -v(\pi)$.  So the subgroup $v(K^\times)\subset\R$ is discrete,
and generated by $v(\pi)$.

Conversely, suppose that $a>0$ generates $v(K^\times)$ and let $\pi\in
v^{-1}(a)$~; we have $\pi\in\pgoth$.  For every $x\neq0$ in $\pgoth$, there is
an integer $n>0$ such $v(x)=na=v(\pi^n)$.  Consequently,
$x/\pi^n\in\ogoth^\times$, and $\pi$ generates the ideal~$\pgoth$.

When $v$ is discrete, every generator $\pi$ of the maximal ideal $\pgoth$ is
called a {\it uniformiser\/} for $v$.

\th COROLLARY 4
\enonce
Suppose that\/ $v$ is discrete and let\/ $\pi$ be a uniformiser.  Every\/
$x\in K^\times$ can be uniquely written as\/ $x=u\pi^n$, where\/
$u\in\ogoth^\times$ and $n\in\Z$.  If\/ $v$ is normalised, then $n=v(x)$, and
the map\/ $\Z\times\ogoth^\times\to K^\times$, $(n,u)\mapsto u\pi^n$, is an
isomorphism. 
\endth

Let $K$ be a valued field.  A series $\sum_{i>0} x_i$ of elements of $K$ is
said to {\it converge\/} to the sum $s\in K$ if the sequence
$s_n=\sum_{i\in[1,n]}x_i$ converges to the limit $s$ as $n\to+\infty$.  If
such is the case, and if the absolute value comes from a valuation $v$, then
$v(s)\ge\inf\limits_{i>0} v(x_i)$.

\th LEMMA 5
\enonce
Suppose that\/ $K$ is complete for a valuation.  Then a series\/ $\sum_{i>0}
x_i$ converges if and only if\/ $x_i\to0$ as\/ $i\to+\infty$.
\endth
Suppose that the series converges to $s$.  Then (even if $K$ is archimedean 
or incomplete)
$$
\lim_{i\to+\infty}x_i
=\lim_{i\to+\infty}(s_{i+1}-s_i)
=\lim_{i\to+\infty}s_{i+1}-\lim_{i\to+\infty}s_i
=s-s=0.
$$
Suppose that, given $\varepsilon>0$, there is an $N$ such that
$|x_n|<\varepsilon$ for all $n>N$.  Now, for $j>i>N$, we have
$$
|s_j-s_i|=|x_{i+1}+\cdots+x_j|\le\sup_{n\in]i,j]}|x_n|<\varepsilon,
$$
so $(s_i)_i$ is a fundamental sequence and, $K$ being complete, has a limit.

Suppose that $K$ is complete for a normalised discrete valuation $v$, let
$\pi$ be a uniformiser, let $A^*\subset\ogoth^\times$ be a system of
representatives of $k^\times$, and put $A=A^*\cup\{0\}$.

\th THEOREM 6
\enonce
Every unit\/ $u\in\ogoth^\times$ in a complete discretely valued field can be 
uniquely written as 
$$
u=\sum_{i\in[0,+\infty[}a_i\pi^i,\qquad (a_0\in A^*, a_i\in A),
$$
and conversely, every such series converges to an element of\/
$\ogoth^\times$.  
\endth
Every such series is convergent because $a_i\pi^i\to0$ as $i\to+\infty$
(lemma~5)~; the sum is a unit because $\ogoth^\times$ is closed in $K^\times$.

Let $a_0\in A^*$ be the representative of $\bar u\in k^\times$.  Then
$u=a_0+b_1\pi$ for some $b_1\in\ogoth$.  Let $a_1$ be the representative of
$\bar b_1\in k$, so that $b_1=a_1+b_2\pi$ and $u=a_0+a_1\pi+b_2\pi^2$ for some
$b_2\in\ogoth$.  Proceeding in this manner, we may write, for every $n\in\N$,
$$
u=a_0+a_1\pi+\cdots+a_n\pi^n+b_{n+1}\pi^{n+1}
$$
where $a_0\in A^*$, $a_1,\ldots,a_n\in A$ and $b_{n+1}\in\ogoth$ are
uniquely determined by $u$.  As $\lim_{n\to+\infty}b_n\pi^n=0$, the series
converges to $u$, as claimed.

\th COROLLARY 7
\enonce
Every\/ $x\in K^\times$ in a complete discretely valued field can be uniquely
written as 
$$
x=\sum_{i\in[n,+\infty[}a_i\pi^i,\qquad (a_n\in A^*, a_i\in A).
$$
and conversely, every such series has a sum of valuation\/ $n$, if\/
$v(\pi)=1$.   
\endth
Apply the preceding theorem to the unit $u=x/\pi^{v(x)}$ (cor.~4).

The displayed series of sum $x$ is called the $\pi$-{\it adic expansion\/} of
$x$ (relative to the system $A$ of representative of the residue field $k$.)
One may declare the $\pi$-adic expansion of $0$ to be $\sum_i0\pi^i$.

{\it Example}\pointir Take $K=\Q\series{T}$, $\pi=T$, $A=\Q$, $u=T/(e^T-1)$.
Then the $T$-adic expansion of $u$ is $u=\sum_{i\in[0,+\infty[}(B_i/i!)T^i$
(see Lecture~3).

{\it Example}\pointir Take $K=\Q_p$ ($p$ prime), $\pi=p$,
$A=\{0,1,\ldots,p-1\}$.  This is how Hensel had defined $p$-adic numbers, as
being ``formal'' $p$-adic expansions. We have $1/(1-p)=1+p+p^2+\cdots$.

\th LEMMA 8
\enonce
Suppose that\/ $K$ is a complete discretely valued field and that the residue
field\/ $k$ is finite.  Then\/ $\ogoth$ is compact. 
\endth
For a metric space such as $\ogoth$, compactness is equivalent to every
sequence $(x_j)_j$ having a convergent subsequence.  We choose a uniformiser
$\pi$, a set of representatives $A\subset\ogoth$ of $k$ such that $0\in A$,
and apply the diagonal argument to the $\pi$-adic expansions 
$$
x_j=\sum_{n\in[0,+\infty[}a_{j,n}\pi^n\qquad(a_{j,n}\in A)
$$
furnished by cor.~7.  Since $A$ is finite, there is some $b_0\in A$ which
occurs as $a_{j,0}$ for infinitely many $j$.  For the $x_j$ with
$a_{j,0}=b_0$, there is some $b_1\in A$ which occurs as $a_{j,1}$ for
infinitely many $j$. For the $x_j$ with $a_{j,0}=b_0$ and $a_{j,1}=b_1$, there
is some $b_2\in A$ which occurs as $a_{j,2}$ for infinitely many $j$.  And so
on.  There is thus a subsequence tending to $\sum_n b_n\pi^n$.

To see the same thing more directly, observe that, when $v$ is discrete and
$K$ complete, $\ogoth=\lim\limits_{\leftarrow n}\ogoth/\pgoth^n$.  If moreover
the residue field $k$ is finite, then so is each $\ogoth/\pgoth^n$, and
therefore $\ogoth$, being a profinite ring, is compact and totally
disconnected.

\th COROLLARY 9
\enonce
For a field\/ $K$ with a valuation\/ $v$ to be locally compact, it is
necessary and sufficient that\/ 1) $K$ be complete, 2) $v$ be discrete, and 3)
$k$ be finite.  
\endth

% The three conditions are clearly necessary~; sufficiency follows from
% the compactness of $\ogoth$ (lemma~8). 
% 
% There are valued fields satisfying any subset of these three
%   conditions, but not the remaining ones.

\vfill\eject

\centerline{\bf Lecture 6}
\bigskip
\centerline{\it Hensel's lemma for complete unarchimedean fields}
\bigskip

Let\/ $K$ be a field complete for a valuation $v$, with the corresponding
absolute value $|\ |$.  Let $\ogoth$ the ring of
$v$-integers, and $k$ the residue field.

\th THEOREM 1 (``Hensel's lemma'')
\enonce
Suppose that\/ $K$ is complete for the unarchimedean absolute value\/ $|\ |$.
Let\/ $f\in\ogoth[T]$ and\/ $a_0\in\ogoth$ be such that\/ 
$|f(a_0)|<|f'(a_0)|^2$.  Then the sequence
$$
a_{i+1}=a_i-{f(a_i)\over f'(a_i)}\qquad(i\ge0)
$$
converges to a root\/ $a\in\ogoth$ of\/ $f$.  Moreover, $|a-a_0|\le
|f(a_0)/f'(a_0)^2|<1$.   
\endth

{\it Proof\/} (Lang, 1952)~: 
Notice that $|f(a_0)|<1$ and $0<|f'(a_0)|\le1$, because $f\in\ogoth[T]$ and
$a_0\in\ogoth$.  Put $C=|f(a_0)/f'(a_0)^2|<1$.  We show inductively that
$$
1)\ |a_i|\le1,\quad 
2)\ |a_i-a_0|\le C,\quad 
3)\ f'(a_i)\neq0\hbox{ and }\left|f(a_i)/f'(a_i)^2\right|\le C^{2^i}. 
$$
If such is the case for all $i\ge0$, then $(a_i)_i$ is a fundamental sequence
because 
$$
|a_{i+1}-a_i|=\left|{f(a_i)\over f'(a_i)}\right|
=\left|{f(a_i)\over f'(a_i)^2}\right||f'(a_i)|\le C^{2^i}
$$
and consequently $|a_j-a_i|\le C^{2^i}$ for $j>i$, by the ultrametric
inequality.  The limit $a\in\ogoth$ satisfies $f(a)=0$, as one sees by letting 
$i\to+\infty$ in $f(a_i)=(a_i-a_{i+1})f'(a_i)$.

Let us carry out the induction~; the case $i=0$ is true by hypothesis.  Assume
that the three assertions are true for some $i\ge0$.  Then,

1). $\left|f(a_i)/f'(a_i)^2\right|\le C^{2^i}$ gives, as we have seen,
  $|a_{i+1}-a_i|\le C^{2^i}<1$, and hence $|a_{i+1}|\le1$. 

2). $|a_{i+1}-a_0|\le\sup(|a_{i+1}-a_i|, |a_i-a_0|)\le C$.

3). Define $f_j\in\ogoth[T]$ by the polynomial identity (``\thinspace Taylor 
  expansion\thinspace'')  
$$
f(T+S)=f(T)+f_1(T)S+f_2(T)S^2+\cdots
$$
so that $f'=f_1$, and substitute $T=a_i$, $S=-f(a_i)/f'(a_i)$ (which is a
$v$-integer, being equal to $a_{i+1}-a_i$) to get
$$
f(a_{i+1})=f(a_i)
+f'(a_i)\left({-f(a_i)\over f'(a_i)}\right)+
(v\hbox{-integer}).\left({-f(a_i)\over f'(a_i)}\right)^2,
$$
so that $|f(a_{i+1})|\le |f(a_i)/f'(a_i)|^2$.  Applying
the same method to $f'$ instead of $f$, we get 
$$
f'(a_{i+1})=f'(a_i)+(v\hbox{-integer}).\left({-f(a_i)\over f'(a_i)}\right)
$$
and consequently 
$$
{f'(a_{i+1})\over f'(a_i)}=1
+(v\hbox{-integer}).\left({-f(a_i)\over f'(a_i)^2}\right).
$$
Hence $|f'(a_{i+1})/f'(a_i)|=1$.  In particular, $f'(a_{i+1})\neq0$, and we
have
$$
\left|{f(a_{i+1})\over f'(a_{i+1})^2}\right|
\le\left|{f(a_i)^2\over f'(a_i)^2f'(a_{i+1})^2}\right|
=\left|{f(a_i)\over f'(a_i)^2}\right|^2\le C^{2^{i+1}},
$$
completing the induction.  Finally, the fact that $|a-a_0|\le C$ follows by
taking the limit $i\to+\infty$ in $2)$.  This completes the proof.

\th COROLLARY 2
\enonce
If\/ $\bar f$ has a simple root\/ $x\in k$, then\/ $f$ has a simple root\/
$a\in\ogoth$ such that\/ $\bar a=x$.
\endth

\th COROLLARY 3
\enonce
Let\/ $p\neq2$ be a prime, and suppose that\/ $u\in\Zp^\times$ is such that\/ 
$\bar u\in\Fp^{\times2}$.  Then $u=a^2$ for some $a\in\Zp^\times$.
\endth
Apply cor.~2, or take $a_0\in\Z_p^\times$ such that $\bar a_0^2=\bar u$, and
take $f=T^2-u$.  Then $|f'(a_0)|=|2a_0|=1$.

\th COROLLARY 4
\enonce
Let\/ $p\neq2$ be a prime, and let\/ $b\in\Zp^\times$ be such that\/ $\bar 
b\not\in\Fp^{\times2}$.  The\/ $\F_2$-space\/ $\Qp^\times\!/\Qp^{\times2}$ has 
dimension\/~$2$, and\/ $\hat b,\hat p$ is a basis.  The three quadratic
extensions of\/ $\Qp$ are obtained by adjoining $\sqrt b, \sqrt p,
\sqrt{bp}$. 
\endth

\th COROLLARY 5
\enonce
Suppose that\/ $u\in\Z_2^\times$ is such that\/ $u\equiv1\pmod{2^3}$.  Then
$u=a^2$ for some $a\in\Z_2^\times$. 
\endth
Take $a_0=1$ and $f=T^2-u$.  We have $v_2(f(a_0))\ge3$ and $v_2(f'(a_0))=1$, 
and therefore $v_2(f(a_0))>2v_2(f'(a_0))$.  Now apply th.~1.

\th COROLLARY 6
\enonce
The\/ $\F_2$-space\/ $\Q_2^\times\!/\Q_2^{\times2}$ has 
dimension\/~$3$, and\/ $\hat 5, \hat 3, \hat 2$ is a basis.  The seven
quadratic extensions of\/ $\Q_2$ are obtained by adjoining $\sqrt5,\sqrt3,
\sqrt{15}, \sqrt2, \sqrt{10}, \sqrt6, \sqrt{30}$.
\endth
For example, $\Q_2(\sqrt{-1})$ is the same as $\Q_2(\sqrt{15})$, because
$-1\equiv15\pmod{2^3}$. 

\th COROLLARY 7
\enonce
Let\/ $p$ be a prime.  The group\/ $\Zp^\times$ has an element of order\/
$p-1$. 
\endth

Fix a generator $g\in\Fp^\times$, take $a_0\in\Zp^\times$ to be any lift of
$g$, and take $f=T^{p-1}-1$.  Now $|f'(a_0)|_p=1$, whereas $|f(a_0)|_p\le
p^{-1}$, so th.~1 provides a root $a\in\Zp$ of $f$ such that $|a-a_0|_p\le
p^{-1}$.  This $a$ is in $\Zp^\times$ and has order $p-1$.  (In fact, the
torsion subgroup of $\Qp^\times$ is generated by $a$ for
$p\neq2$, by $-1$ for $p=2$.)

Let $p$ be a prime number, and let $U_n=\Ker(\Zp^\times\to(\Z/p^n\Z)^\times)$
($n>0$).  Here is a generalisation of cor.~5 from $2$ to $p$.  It can be
checked that raising-to-the-exponent-$p$ map $(\ )^p$ takes $U_n$ into
$U_{n+1}$.

\th PROPOSITION 8
\enonce
For every\/ $n>1/(p-1)$, the map\/ $(\ )^p:U_n\to U_{n+1}$ is surjective.   
\endth
Let $y\in U_{n+1}$ and write $y=1+bp^{n+1}$
(with $b\in\Zp$).  We seek a root of $x^p=y$ such that $x=1+ap^n$ for
some $a\in\Zp$.  This leads to the equation
$$
1+bp^{n+1}=1+pa.p^n+\cdots+pa^{p-1}.p^{n(p-1)}+a^p.p^{np},
$$
in which the {\it coefficients\/} of $p^i$ ($i=n,\ldots,n(p-1)$) are all
divisible by $p$.  Also, the exponents of $p$ in all but the first two terms
on the right are $>n+1$, unless $pn=n+1$, which cannot happen because $n$ is
$>1/(p-1)$.  The equation can therefore be rewritten as
$$
b=a+p g(a)
$$
for some polynomial $g\in\Zp[T]$.  Applying cor.~2 to $f=b-T-pg(T)$, we get a
simple root $a\in\Zp$ of $f$ such that $\bar a=\bar b$.  Then $x^p=y$, and
$(\ )^p$ is surjective.  (We shall see later that it is also injective.)

The only case excluded by prop.~8 is when $p=2$ and $n=1$.  Indeed, the map
$(\ )^2:U_1\to U_2$ is neither injective nor surjective~; the kernel is
generated by $-1$ and the cokernel by $1+2^2$.

\vfill\eject

\centerline{\bf Lecture 7}
\bigskip
\centerline{\it Hensel's lemma for complete discretely valued fields}
\bigskip

Let us recall the notion of the {\it resultant\/} $\res(g,h)\in A$ of two
polynomials $g,h\in A[T]$ over a commutative ring $A$.

Let $m$, $n$ be positive integers and let $g,h\in A[T]$ be two polynomials 
$$
g=b_mT^m+b_{m-1}T^{m-1}+\cdots+b_0,\quad 
h=c_nT^n+c_{n-1}T^{n-1}+\cdots+c_0,
$$
of degrees less than $m$, $n$ respectively.  The $(m,n)$-resultant
$\res_{m,n}(g,h)$ of $g,h$ is the determinant of the square matrix of size
$m+n$ which, for $m=2$ and $n=3$, would be written as
$$
\pmatrix{
b_2&0&0&c_3&0\cr
b_1&b_2&0&c_2&c_3\cr
b_0&b_1&b_2&c_1&c_2\cr
0&b_0&b_1&c_0&c_1\cr
0&0&b_0&0&c_0\cr
}.$$
If $m=\deg(g)$ and $n=\deg(h)$, then we put $\res(g,h)=\res_{m,n}(g,h)$ and
call it the resultant of $g$, $h$.  For a homomorphism $\rho:A\to B$ of
commutative rings, $\res_{m,n}({}^\rho g,{}^\rho h)=\rho(\res_{m,n}(g,h))$,
where ${}^\rho g,{}^\rho h\in B[T]$ are the images of $g,h$.

For $g$ of degree~$m$, the {\it discriminant\/} $\dis(g)$ satisfies
$$
\res_{m,m-1}(g,g')=\res_{m-1,m}(g',g)=(-1)^{m(m-1)/2}b_m\dis(g).
$$

When $m+n>0$, which we assume from now on, there exist $H,G\in A[T]$,
with $\deg(H)<\deg(h)$ and $\deg(G)<\deg(g)$, such that
$$
\res(g,h)=gH+hG.\leqno{(1)}
$$
In particular, $\res(g,h)$ is in the ideal $(g,h)\subset A[T]$ generated by
$g$ and $h$.  The pair $H,G$ is unique when $\res(g,h)$ is simplifiable in $A$
--- multiplication by $\res(g,h)$ is injective in $A$.  Such is the
case when $A$ is integral and $\res(g,h)\neq0$.

Now let $K$ be a discretely valued field $K$, $v$ its valuation, $\ogoth$ the
ring of $v$-integers, and $\pi$ a uniformiser. Take $A=\ogoth$, which is an
integral ring, and suppose that $\res(g,h)\not\equiv0\nomod{\pi^{\alpha+1}}$
for some $\alpha\in\N$.  Multiplying the relation $(1)$ by units and a
suitable power of $\pi$, we may assume that $H,G$ satisfy
$$
\pi^\alpha=gH+hG.\leqno{(2)}
$$

Let $F\in\ogoth[T]$ be some polynomial such that $\deg(F)<\deg(gh)$.
Multiplying $(2)$ by $F$, we get
$$
\pi^\alpha F=gHF+hGF.\leqno{(3)}
$$
If $\deg(HF)<\deg(h)$ and $\deg(GF)<\deg(g)$, take $H'=HF$, $G'=GF$.
Otherwise, let $H'$ be the part of $HF$ so that $\deg(H')<\deg(h)$ and
similarly let $G'$ be the part of $GF$ with $\deg(G')<\deg(g)$~; by
hypothesis, one of $HF-H'$, $GF-G'$ is $\neq0$ and indeed of degree at least
$\deg(h)$ or $\deg(g)$, as the case may be.  Rewrite $(3)$ as
$$
\pi^\alpha F=gH'+hG'+F',\quad F'=g(HF-H')+h(GF-G'),
$$
defining $F'$.  Now, $\deg(F')$ is at least $\deg(gh)$, unless $F'=0$.  As
the left-hand side of the first equation has degree $<\deg(gh)$, as do the
terms $gH'$, $hG'$ on the right, we must have $F'=0$, and we get the following
lemma.

\th LEMMA 1
\enonce
Let\/ $g,h,F\in\ogoth[T]$ be polynomials such that\/ $\deg(gh)>0$ and
$\deg(F)<\deg(gh)$.  Suppose that\/
$\res(g,h)\not\equiv0\nomod{\pi^{\alpha+1}}$ for some\/ $\alpha\in\N$.  Then
there exists a unique pair of polynomials\/ $H',G'\in\ogoth[T]$ such that
$$
\deg(H')<\deg(h),\quad \deg(G')<\deg(g),\quad \pi^\alpha F=gH'+hG'.
$$
\endth
As we have seen, $H',G'$ are truncations of $HF,GF$, where $H,G$ is the unique
solution of $(2)$, which depends only on $g,h$ (and the choice of $\pi$), not
on $F$.

Assume now that $K$ is complete for the discrete valuation $v$.

\th THEOREM 2 (``Hensel's lemma'')
\enonce
Let\/ $f,g_0,h_0\in\ogoth[T]$ be polynomials such that\/ $f$ and\/ $g_0h_0$
have the same leading term.  Suppose that
$$
\res(g_0,h_0)\not\equiv0\nomod{\pi^{\alpha+1}},\quad 
f\equiv g_0h_0\nomod{\pi^{2\alpha+1}}\leqno{(4)}
$$
for some $\alpha\in\N$.  Then there exists a unique pair $g,h\in\ogoth[T]$,
having the same leading terms as\/ $g_0,h_0$, and such that\/
$$
g\equiv g_0\nomod{\pi^{\alpha+1}},\quad 
h\equiv h_0\nomod{\pi^{\alpha+1}},\quad 
f=gh.
\leqno{(5)}
$$
\endth 
As a first step, let us construct polynomials $g_1,h_1\in\ogoth[T]$,
having the same leading terms as\/ $f_0,g_0$, and such that
$$
g_1\equiv g_0\nomod{\pi^{\alpha+1}},\quad
h_1\equiv h_0\nomod{\pi^{\alpha+1}},\quad
f\equiv g_1h_1\nomod{\pi^{2\alpha+1+1}},\leqno{(6)}
$$
thereby improving the last congruence in $(4)$. 
These requirements amount to finding $G_1,H_1\in\ogoth[T]$, of degrees
strictly less than those of $g_0$ and $h_0$ respectively, such that, upon
taking
$$
g_1=g_0+\pi^{\alpha+1}G_1,\quad
h_1=h_0+\pi^{\alpha+1}H_1,\quad
$$
the last congruence $(6)$ is satisfied.  Writing
$f-g_0h_0=\pi^{2\alpha+1}F_1$, where $F_1\in\ogoth[T]$ has degree
$<\deg(f)=\deg(g_0h_0)$, the required congruence can be translated as
$$
\pi^\alpha F_1\equiv g_0H_1+h_0G_1\nomod{\pi^{\alpha+1}}.\leqno{(7)}
$$
But lemma~$1$ is applicable to $g_0,h_0,F_1$~; the pair
$G_1,H_1$ furnished by it has the desired properties.

This process can be repeated.  More precisely, suppose that we have been able
to find, for some integer $i>0$, polynomials $g_i,h_i\in\ogoth[T]$,
having the same leading terms as\/ $g_{i-1},h_{i-1}$, and such that
$$
g_i\equiv g_{i-1}\nomod{\pi^{\alpha+i}},\quad
h_i\equiv h_{i-1}\nomod{\pi^{\alpha+i}},\quad
f\equiv g_ih_i\nomod{\pi^{2\alpha+1+i}},\leqno{(8)}
$$
as we have just done for $i=1$.  Writing
$f-g_ih_i=\pi^{2\alpha+1+i}F_{i+1}$, where $F_{i+1}\in\ogoth[T]$, and noting
that $\res(g_i,h_i)\equiv\res(g_0,h_0)\not\equiv0\nomod{\pi^{\alpha+1}}$,
there is a unique pair of polynomials $G_{i+1},H_{i+1}\in\ogoth[T]$, of
degrees strictly less than those of $g_i$ and $h_i$ respectively, such that
$$
\pi^\alpha F_{i+1}= g_iH_{i+1}+h_iG_{i+1},
$$
as follows from lemma~1.  It is now easily seen that
$g_{i+1}=g_i+\pi^{\alpha+i+1}G_{i+1}$, $h_{i+1}=h_i+\pi^{\alpha+i+1}H_{i+1}$ 
satisfy $(8)$ with $i$ replaced
by $i+1$.  Finally take % $g=\lim g_i$, $h=\lim h_i$, 
$$
\eqalign{
g&=\lim_{i\to+\infty} g_i=g_0+\pi^{\alpha+1}G_1+\pi^{\alpha+2}G_2+\cdots,\cr
h&=\lim_{i\to+\infty} h_i=h_0+\pi^{\alpha+1}H_1+\pi^{\alpha+2}H_2+\cdots,\cr
}
$$
so that $g,h$ have all the properties claimed for them, such as $f=gh$~;
this last equality follows because
$f=\lim_i(g_ih_i)=(\lim_ig_i)(\lim_ih_i)=gh$. 

\th COROLLARY 3
\enonce
Let\/ $f,g_0,h_0\in\ogoth[T]$ be unitary polynomials such that\/ 
$$
\bar f=\bar g_0\bar h_0,\quad \gcd(\bar g_0,\bar h_0)=1.
$$ 
Then there is a unique pair of unitary polynomials\/ $g,h\in\ogoth[T]$ such
that\/ $f=gh$ and\/ $\bar g=\bar g_0$, $\bar h=\bar h_0$.  
\endth 

That $f$ and $g_0h_0$ have the same leading term follows from the fact that
they are unitary and $\bar f=\bar g_0\bar h_0$.  Also, as $\bar g_0$, $\bar
h_0$ are relatively prime, $\res(\bar g_0,\bar h_0)\neq0$ in
$k=\ogoth/\pi\ogoth$.  So th.~2 can be applied with $\alpha=0$.

\vfill\eject

\centerline{\bf Lecture 8}
\bigskip
\centerline{\it Valuations on purely transcendental extensions}
\bigskip

Let\/ $K$ be a field with a valuation\/ $v:K^\times\to\R$.  Every extension
$w$ of $v$ to $K(T)$ satisfies $w(a_iT^i)=iw(T)+v(a_i)$ for every monomial
$a_iT^i$, and
$$
w(f)\ge \inf_{i\in[0,m]}(iw(T)+v(a_i)),\leqno{(1)}
$$
for every polynomial $f=a_0+a_1T+\cdots+a_mT^m$ in $K[T]$.  It turns out that
there are extensions $w$ of $v$ for which equality holds in $(1)$, and that
such an extension is unique if we specify $w(T)$, which can be arbitrary in
$\R$.  This observation is enshrined in the following lemma.

\th LEMMA 1
\enonce
Let\/ $C\in\R$ be a constant.  There is a unique valuation\/
$w:K(T)^\times\to\R$ extending\/ $v$ such that\/ $w(T)=C$ and\/ equality holds
in $(1)$ for every polynomial\/ $f=a_0+a_1T+\cdots+a_mT^m$ in\/ $K[T]$.
\endth 
Defining $w(f)$ by requiring that $(1)$ be an equality, it is clear that
$$ 
w(f+g)\ge\inf(w(f),w(g)).\leqno{(2)}
$$
for all $f,g\in K[T]$.  Let us show that we also have $w(fg)=w(f)+w(g)$. 

It is easy to see that $w(fg)\ge w(f)+w(g)$.  To show the converse, let
$r\in[0,m]$ be the smallest index such that $w(f)=rC+v(a_r)$.  Similarly for
$g=b_0+b_1T+\cdots+b_nT^n$, let $s\in[0,n]$ be the smallest index such that
$w(g)=sC+v(b_s)$.  The coefficient of the term of degree $r+s$ in $fg$ is
$d=\sum_{i+j=r+s}a_ib_j$.  Let us compute $v(d)$.

If $i<r$, then $w(a_iT^i)>w(f)$.  Hence $v(a_i)>-iC+w(f)$.  Also, as
$w(b_jT^j)\ge w(g)$ for any $j$ (such that $i+j=r+s)$, we have $v(b_j)\ge
-jC+w(g)$.  Consequently,
$$
v(a_ib_j)>(-r-s)C+w(f)+w(g).
$$
This is also valid for any $j<s$ and any $i$ (such that $i+j=r+s$), for
similar reasons.

If $i=r$ and $j=s$, then $v(a_r)=-rC+w(f)$ and $v(b_s)=-sC+w(g)$, so $v(a_r
b_s)=(-r-s)C+w(f)+w(g)$.  Hence $v(d)=v(a_rb_s)$ and
$$
w(fg)\le w(dT^{r+s})=w(f)+w(g),
$$
as was to be shown.  So $w(fg)= w(f)+w(g)$ holds for all $f,g\in K[T]$.  To
complete the proof, we have to show that $w$ can be extended from $K[T]$ to a
valuation on $K(T)$.  This is purely formal.

Let $A$ be an integral commutative ring, $M=\{f\in A\,|\,f\neq0\}$ its
multiplicative monoid, $k$ its field of fractions, and $w:M\to\R$ a
homomorphism of monoids such that $(2)$ holds for every $f,g\in A$.  Then the
unique extension $w$ of $w$ from $M$ to its ``group of fractions'' $k^\times$
is a valuation on $k$.  In other words, $(2)$ continues to hold for $f,g\in
k$.

From now on, $\T=(T_1,T_2,\ldots,T_n)$ is a family of indeterminates, and
denote the ring $A[T_1,T_2,\ldots,T_n]$ by $A[\T]$, for every commutative ring
$A$.  Also, for every field $k$, let $k(\T)$ stand for $k(T_1,T_2,\ldots,T_n)$.
For $\alpha\in\N^n$ and $C\in\R^n$, put $\alpha\cdot
C=i_1C_1+i_2C_2+\cdots+i_nC_n$.  Also put $\T^\alpha=T_1^{i_1}T_2^{i_2}\ldots
T_n^{i_n}$.

Recall that we are woking with a field $K$ with a height-$1$ valuation $v$.

\th COROLLARY 2
\enonce
For every\/ $C\in\R^n$, there is a unique extension of $v$ to a valuation\/
$w$ on\/ $K(\T)$ such that 
$$
w(f)=\inf_{\alpha}(\alpha\cdot C+v(b_\alpha))
$$
holds for every polynomial\/ $f=\sum_{\alpha\in\N^n}b_\alpha \T^\alpha$ in\/
$K[\T]$.  
\endth
Because $K(\T)=K(T_1,T_2,\ldots,T_{n-1})(T_n)$, the corollary
follows from lemma~1 by induction on $n$.

Recall that $\ogoth$ denotes the ring of $v$-integers in $K$.

\th LEMMA 3 
\enonce
Suppose that\/ $f\in\ogoth[\T]$ can be written\/ $f=gh$ for
some\/ $g,h\in K[\T]$.  Then there exists a\/ $b\in K^\times$
such that\/ $b^{-1}g,bh\in\ogoth[\T]$.
\endth
Extend the valuation $v$ to a valuation $w$ on $K(\T)$ by taking $C=0$ in
$\R^n$ (cor.~2).  Then $\ogoth[\T]$ is contained in the ring $\goth O$ of
$w$-integers of $K(\T)$, and in fact $\ogoth[\T]={\goth O}\cap K[\T]$.
Moreover, $w(K(\T)^\times)=v(K^\times)$.

Let $f\neq0$ be a polynomial in $\ogoth[\T]$, and write $f=gh$ as given.  There
is a $b\in K^\times$ such that $v(b)=w(g)$.  Now,
$w(b^{-1}g)=v(b^{-1})+w(g)=0$, so $b^{-1}g\in\ogoth[T]$.  Further,
$$
0\le w(f)=w(b^{-1}g)+w(bh)=w(bh),
$$
so $bh\in\ogoth[T]$ as well, proving the claim.  Finally, if $f=0$, then
either $g=0$ or $h=0$, and the existence of $b$ is clear.

\th COROLLARY 4 
\enonce
If\/ $f$ is irreducible in\/ $\ogoth[\T]$, then it is irreducible in\/ $K[\T]$.
\endth

\th LEMMA 5 (Gau\ss)
\enonce
Suppose that\/ $f\in\Z[\T]$ can be written\/ $f=gh$ for some\/ $g,h\in\Q[\T]$.
Then there exists a\/ $b\in\Q^\times$ such that\/ $b^{-1}g,bh\in\Z[\T]$.
\endth 

It is enough to show that there is a $b\in\Q^\times$ such that
$b^{-1}g,bh\in\Z_p[\T]$ for every prime $p$.

For every prime $p$, there is a $b_p\in\Qp^\times$ such that
$b_p^{-1}g,b_ph\in\Z_p[\T]$ (lemma~3).  We may assume that $b_p=p^{m_p}$
($m_p\in\Z$) and that $m_p=0$ for almost all $p$.  Take $b=\prod_pb_p$.

\th COROLLARY 6
\enonce
If\/ $f$ is irreducible in\/ $\Z[\T]$, then it is irreducible in\/ $\Q[\T]$.   
\endth

Now suppose that the valuation $v$ is discrete and normalised, denote by $k$
its residue field, and let there be only one indeterminate, $T$.

\th LEMMA 7 (Eisenstein)
\enonce
Let\/ $f=a_0+a_1T+\cdots+a_mT^m$ be a polynomial in\/ $\ogoth[T]$ such that
$$
v(a_0)=1,\quad v(a_i)>0\ \ (0<i<m),\quad v(a_m)=0.
$$
Then\/ $f$ is irreducible in\/ $K[T]$.
\endth
It is enough to show that $f$ is irreducible in $\ogoth[T]$ (cor.~4).  If not,
write $f=gh$, where
$$
g=b_rT^r+b_{r-1}T^{r-1}+\cdots+b_0,\quad 
h=c_sT^s+c_{s-1}T^{s-1}+\cdots+c_0,\quad 
$$
are in $\ogoth[T]$ and $r+s=m$.  Passing to the residue field $k$ of
$\ogoth$, we have $\bar f=\bar a_mT^m$, so $\bar g=\bar b_rT^r$ and $\bar
h=\bar c_sT^s$.  In particular, $v(b_0)\ge1$ and $v(c_0)\ge1$.  But then
$v(a_0)=v(b_0c_0)\ge2$, contradicting the hypothesis that $v(a_0)=1$.

Polynomials $f$ satisfying the hypothesis of lemma~7 are called {\it
  Eisenstein polynomials}.

Before giving the first application of this lemma, recall that, for every
prime number $p$, the binomial coefficients $p\choose i$ are divisible by~$p$
for $i\in[1,p[$.  This is clear for $i=1$, as ${p\choose 1}=p$.  For
$i\in[2,p[$, the identity $i{p\choose i}=(p-i+1){p\choose i-1}$ shows that if
$p$ divides ${p\choose i-1}$, then it divides ${p\choose i}$.

\th EXAMPLE 8
\enonce
For every prime\/ $p$, the cyclotomic polynomial\/
$\Phi_p(T)=T^{p-1}+T^{p-2}+\cdots+1$ is irreducible in $\Qp[T]$.
\endth
We have 
$$
\Phi_p(T)={T^p-1\over T-1},\quad
\Phi_p(S+1)=S^{p-1}+pS^{p-2}+\cdots+p,
$$
where the suppressed terms have coefficients divisible by $p$, as we have
just remarked.  Thus $\Phi_p(S+1)$ is an Eisenstein polynomial, so it is
irreducible in $\Qp[T]$, and hence so is $\Phi_p(T)$.

\th COROLLARY 9
\enonce
For every\/ $n>1/(p-1)$, the map\/ $(\ )^p:U_n\to U_{n+1}$ is bijective.   
\endth
We have already seen that it is surjective (Lecture~6, prop.~8).  Every
element $\zeta\neq1$ in the kernel of $(\ )^p$ is a root of the cyclotomic
polynomial $\Phi_p$, for $\Phi_p(\zeta)=(\zeta^p-1)/(\zeta-1)=0$.  But
$\Phi_p$ is irreducible (example~8), and cannot have a root if
$\deg(\Phi)=p-1>1$, which is the case if $p\neq2$.  If $p=2$, then $\zeta=-1$,
but $-1\notin U_2$, so $(\ )^2$ is injective on $U_n$ if $n>1$, which is the
case by hypothesis.

\th EXAMPLE 10
\enonce
For every prime\/ $p$ and every\/ $n>0$, the cyclotomic polynomial\/
$\Phi_{p^n}(T)=\Phi_p(T^{p^{n-1}})$ is irreducible in $\Qp[T]$.
\endth
Put $\theta(S)=\Phi_{p^n}(S+1)$.  Then $\theta(0)=\Phi_p(1)=p$.  Moreover,
from the identity $\Phi_p(T)=(T^p-1)/(T-1)$, we get
$\Phi_{p^n}(T)(T^{p^{n-1}}-1)=T^{p^n}-1$, and hence
$$
\theta(S)((S+1)^{p^{n-1}}-1)=(S+1)^{p^n}-1.
$$
Reading this modulo $p$, we get $\bar\theta(S)S^{p^{n-1}}=S^{p^n}$, because
${p^r\choose t}$ is divisible by $p$ for $r>0$, $0<t<p^r\ (*)$.  So
$\bar\theta(S)=S^{{p^n}-p^{n-1}}$, and hence $\theta$ is an Eisenstein
polynomial and therefore irreducible.  So $\Phi_{p^n}$ is irreducible.

{\eightpoint (*)  To see that ${p^r\choose i}$ is divisible by $p$, proceed by
  induction on $r$.  We have seen this for $r=1$~; it implies the identity
  $1+T^p=(1+T)^p$ in $\F_p(T)$.  By induction we get
$$
1+T^{p^r}=(1+T)^{p^r}=1+T^{p^r}+\sum_{i\in[1,p^r[}{p^r\choose i}T^i,
$$
which implies that $p|{p^r\choose i}$.}

\vfill\eject

\centerline{\bf Lecture 9}
\bigskip
\centerline{\it The Newton polygon}
\bigskip

Let $K$ be a field with a valuation $v$, and let $f=a_0+a_1T+\cdots+a_mT^m$
($m>0$) be a polynomial in\/ $K[T]$.  Suppose that $a_0a_m\neq0$, so $m$ is
the degree of $f$ and $T$ does not divide $f$.

For every $j\in[0,m]$ such that $a_j\neq0$, plot the point $P_j=(j,v(a_j))$ in
the plane $\R^2$, with coordinates $x,y$.  

\th DEFINITION 1
\enonce
The  Newton polygon\/ $\Pi_f$
of $f$ is the lower convex envelope of this set of points $P_j$.
\endth

{\eightpoint Imagine a piece of string affixed to the point $P_0$ and hanging
  vertically down along the $y$-axis.  Pull it tightly counter-clockwise so as
  to make it pass through the point $P_m$.  The piece-wise linear portion of
  the string between $P_0$ and $P_m$ is the polygon~$\Pi_f$.}

Every point $P_j$ lies on or above $\Pi_f$.  Suppose that $\Pi_f$ has $r+1$
vertices (the points $P_j$ on $\Pi_f$ such that the function on the real
interval $[0,m]$ whose graph is $\Pi_f$ is not differentiable at $j$), with
$x$-coordinates
$$
0=m_0<m_1<m_2<\cdots<m_{r-1}<m_r=m.
$$
Thus $\Pi_f$ has $r>0$ sides, joining the vertices
$P_{m_{i-1}}=(m_{i-1},v(a_{m_{i-1}}))$ and $P_{m_i}=(m_i,v(a_{m_i}))$ for
$i\in[1,r]$~; let $\gamma_i$ be the slopes of these sides, so that
$$
\gamma_1<\gamma_2<\cdots<\gamma_r~;\quad
\gamma_i={v(a_{m_i})-v(a_{m_{i-1}})\over m_i-m_{i-1}}.
$$
Let $l_i=m_i-m_{i-1}$ be the length of the projection of side $i$ onto the
$x$-axis, so that $l_1+l_2+\cdots+l_r=m$.

\th DEFINITION 2
\enonce
The type of\/ $f$ is the sequence\/ 
$(l_1,\gamma_1; l_2,\gamma_2;\ldots;l_r,\gamma_r)$.
\endth
The type of $f$ knows about the number $r$ of sides of the polygon $\Pi_f$,
the lengths $l_i$ of their horizontal projections, and their slopes
$\gamma_i$, for $i\in[1,r]$.  Knowing $\Pi_f$ is the same as knowing the type
of $f$ and the point $P_0$.

{\it Example}\pointir The type of the polynomial
$f=1+{1\over1!}T+{1\over2!}T^2+\cdots+{1\over7!}T^7$ over $K=\Q_2$ is
$(4,-{3\over4};2,-{1\over2};1,0)$.  The polygon $\Pi_f$ has three sides.

{\it Example}\pointir More generally, let $p$ be a prime, $m>0$ an
integer, and 
$$
m=b_rp^{n_r}+\cdots+b_2p^{n_2}+b_1p^{n_1},\ \ 
(n_r<\cdots<n_2<n_1,\ 0<b_i<p)
$$
the $p$-adic expansion of $m$, so that $r$ is the number of ``digits''
$\neq0$.  Put $m_0=0$, $m_i=b_ip^{n_i}+\cdots+b_1p^{n_1}$ for $i\in[1,r]$ and
take
$$
f=1+{1\over1!}T+{1\over2!}T^2+\cdots+{1\over m!}T^m
$$
over $K=\Q_p$.  Then the polygon $\Pi_f$ has $r$ sides~; we have
$l_i=b_ip^{n_i}$ for every $i\in[1,r]$, and the $r+1$ vertices and $r$ slopes are
given by
$$
P_{m_0}=(0,0),\ P_{m_i}=(m_i,-v_p(m_i!)),\quad
\gamma_i={-(p^{n_i}-1)\over p^{n_i}(p-1)}.%\quad (i\in[1,r]).
$$

{\it Remark}~3\pointir Let $C\in\R$ be a constant, and consider the unique
valuation $w$ on $K(T)$ which extends $v$, for which $w(T)=C$, and such that
the valuation of a polynomial is the smallest of the valuations of its terms
(Lecture~8).   In particular, $w(f)=\inf_{j\in[0,m]}w(a_jT^j)$.  For which $j$
does the infimum occur~?  For just one~$j$, or for many~?  The Newton polygon
$\Pi_f$ can answer these questions.  An understanding of this point is crucial
for what follows.

We need only consider $j$ for which $a_j\neq0$.  Among the lines of slope $-C$
in the plane, the one which passes through $P_j$ has $y$-intercept
$w(a_jT^j)$.  In other words, the point $P_j$ lies on the line
$$
y+Cx=w(a_jT^j).
$$
So the infimum in question occurs at those $j$ for which the $y$-intercept
of the line of slope $-C$ passing through $P_j$ is the smallest.  If $-C$ is
different from the slopes $\gamma_i$ of the sides of $\Pi_f$, then there is
just one $j$ for which the infimum occurs, and $P_j$ is a vertex of $\Pi_f$.
If, on the other hand, $-C=\gamma_i$ for some $i\in[1,r]$, then the infimum
occurs at $j=m_{i-1}$ and $j=m_i$ at least, and also at every other
$j\in[m_{i-1},m_i]$ for which $P_j$ lies on the side~$i$ of $\Pi_f$, but for
no other~$j$.  Writing $f_i=\sum_{j\in[m_{i-1},m_i]}a_jT^j$, we have
$$
w(f-f_i)>w(f).\leqno{(1)}
$$
\vskip-.5\baselineskip 
{\eightpoint Imagine a line of slope $-C$ whose
  $y$-intercept is so small, so close to $-\infty$, that the points $P_j$ lie
  above it, and move it slowly upwards, keeping the slope unchanged.  If the
  slope $-C$ is different from the $\gamma_i$, the first time the moving line
  meets any of the points $P_j$ happens at a unique vertex of $\Pi_f$, namely
  the vertex $P_{m_i}$ if $\gamma_i<-C<\gamma_{i+1}$, with the momentary
  convention that $\gamma_0=-\infty$, $\gamma_{r+1}=+\infty$.  If, however,
  $-C=\gamma_i$ for some $i\in[1,r]$, the moving line first meets any of the
  points $P_j$ along the side $i$ of $\Pi_f$.}

\th  DEFINITION 4
\enonce
We say that\/ $f$ is {\it pure\/} if\/ $r=1$, if\/ $\Pi_f$ is a line segment.
\endth
This is the same as saying that $\Pi_f$ is the segment joining the points
$P_0=(0,v(a_0))$ and $P_m=(m,v(a_m))$.  If $f$ is pure of type $(l,\gamma)$,
then $l=m$ and $\gamma=v(a_ma_0^{-1})/m$.  We also say that $f$ is pure of
{\it slope\/} $\gamma$.

\th THEOREM 5 (``Newton'')
\enonce
Suppose that the field\/ $K$ is complete for the valuation\/~$v$, and let\/
$f\in K[T]$ be a polynomial of degree\/~$>0$ such that\/ $f(0)\neq0$.  Let\/
$(l_1,\gamma_1;l_2,\gamma_2;\ldots;l_r,\gamma_r)$ be the type of\/ $f$.  Then
there exist\/ $r$ polynomials\/ $g_i\in K[T]$, each pure of type\/
$(l_i,\gamma_i)$, such that\/ $f=g_1g_2\ldots g_r$.  
\endth 
The rest of the lecture is devoted to a proof of this theorem.

\th LEMMA 6
\enonce
Suppose that\/ $f,g\in K[T]$ are pure of type\/ $(m,\gamma)$, $(n,\gamma)$
respectively.  Then\/ $fg$ is pure of type\/ $(m+n,\gamma)$.
\endth
Take $C=-\gamma$, and let $w$ be the corresponding valuation of $K(T)$.  Then 
$$
w(f)=w(a_0)=w(a_mT^m),\quad
w(g)=w(b_0)=w(b_nT^n),
$$
where $b_0$ (resp.~$b_n$) is the constant (resp.~the leading)
coefficient of~$g$.  The corresponding coefficients of~$fg$ are $c_0=a_0b_0$
and $c_{m+n}=a_mb_n$.  On adding the above equalities, we get
$$
w(fg)=w(a_0b_0)=w(a_mb_nT^{m+n}),
$$
so $fg$ is pure of type $(m+n,\gamma)$ by remark~3, proving the claim.

\th LEMMA 7
\enonce
Suppose that\/ $f$ has type\/ $\Gamma=(l_1,\gamma_1;
l_2,\gamma_2;\ldots;l_r,\gamma_r)$, and that\/ $g$ is pure of type\/
$(n,\gamma)$, where $\gamma>\gamma_r$. Then\/ $fg$ has type\/
$(\Gamma\,;n,\gamma)$.
\endth
If $-C=\gamma_i$ for some $i\in[1,r]$, then, by remark~3 and the
hypothesis that $\gamma_r<\gamma$, we have $w(g-b_0)>w(g)$.  Hence and by
$(1)$,
$$
w(fg-b_0f_i)>w(fg).
$$
where $f_i=\sum_{j\in[m_{i-1},m_i]}a_jT^j$.  Similarly,
$w(fg-a_mT^mg)>w(fg)$, if $-C=\gamma$.  These inequalities show (remark~3)
that $\Pi_{fg}$ has at least $r+1$ sides, of respective slopes
$\gamma_1,\gamma_2,\ldots,\gamma_r,\gamma$, having projections onto the
$x$-axis of lengths at least $m_1,m_2,\ldots,m_r,n$ respectively.  But the sum
of these lengths is already equal to the degree $m+n$ of $fg$, so there are no
more sides to $\Pi_{fg}$.  This proves that $fg$ has the type it is claimed to
have.

For the next two lemmas, choose a constant $C\in\R$ and let $w$ be the unique
valuation on $K(T)$ extending $v$, such that $w(T)=C$, and such that every
polynomial
$$
g=b_0+b_1T+\cdots+b_nT^n
$$
in $K[T]$ has valuation equal to the infimum of the valuations of its terms
(Lecture~8, lemma~1). 
\th LEMMA 8 
\enonce
Suppose that\/ $w(g)=w(b_nT^n)$, let\/ $F\in K[T]$ be another polynomial, and 
write\/ $F=qg+r$, where\/ $q,r\in K[T]$ and\/ $\deg(r)<\deg(g)$.  Then
$$
w(qg)\ge w(F),\quad w(r)\ge w(F).
$$
\endth
Let $l=\deg(q)$, so $\deg(F)=l+n$, and write
$F=c_0+c_1T+\cdots+c_{l+n}T^{l+n}$, $q=d_0+d_1T+\cdots+d_lT^l$. 

First consider the coefficients of $T^{l+n}$ on the two sides of the equality
$F=qg+r$.  We have $c_{l+n}=d_lb_n$, and hence $w(c_{l+n}T^{l+n})=w(d_lT^lg)$,
by the hypothesis that $w(g)=w(b_nT^n)$.  It follow that $w(d_lT^lg)\ge w(F)$.

Secondly, consider the coefficients of $T^{l+n-1}$ on the two sides.  We have
$c_{l+n-1}=d_lb_{n-1}+d_{l-1}b_n$, or equivalently
$d_{l-1}b_n=c_{l+n-1}-d_lb_{n-1}$, so
$$
w(d_{l-1}T^{l-1}g)\ge
\inf\left(w(c_{l+n-1}T^{l+n-1}),\; w(d_lT^l.b_{n-1}T^{n-1})\right).
$$
As $w(b_{n-1}T^{n-1})\ge w(g)$, we have $w(d_lT^l.b_{n-1}T^{n-1})\ge
w(d_lT^lg)\ge w(F)$.  It is thus clear that the infimum in question is $\ge
w(F)$.  It follow that $w(d_{l-1}T^{l-1}g)\ge w(F)$.

Proceeding in this manner, we see that $w(d_iT^ig)\ge w(F)$ for every $i$, and
hence that $w(qg)\ge w(F)$.  Finally, $w(r)\ge\inf(w(F),w(qg))=w(F)$.

\th LEMMA 9
\enonce
Suppose that for $f=a_0+a_1T+\cdots+a_mT^m$, there is some $n$ with $0<n<m$
such that $w(a_nT^n)=w(f)$ and $w(a_jT^j)>w(f)$ for $j>n$.  Then there exist
 $g,h\in K[T]$ of degrees $n,m-n$ such that $f=gh$.
\endth
By hypothesis, $w(f-f_n)>w(f)$, where $f_n=a_0+a_1T+\cdots+a_nT^n$, so there
is a $\Delta>0$ such that $w(f-f_n)=w(f)+\Delta$.

Let us consider pairs of polynomials $G\in K[T]$ of degree~$n$ and $H\in K[T]$
of degree $\le m-n$ such that
$$
w(f-G)\ge w(f)+\Delta,\quad w(H-1)\ge\Delta.\leqno{(2)}
$$
Adding $w(G)$ to the second inequality, we get $w(GH-G)\ge w(G)+\Delta$, and
hence $w(f-GH)\ge w(f)+\Delta$.  Defining $\delta$ by $w(f-GH)=w(f)+\delta$,
we have $\delta\ge\Delta$. 

One choice of $G,H$ would be $f_n,1$.  We shall show that if $f\neq GH$,
equivalently if $\delta<+\infty$, then we can find $G',H'$ which satisfy the
same conditions as $G,H$ but for which $\delta'\ge\delta+\Delta$.

It follows from the hypotheses on $f$, $n$ and $G$ that $w(G)=w(b_nT^n)$,
where $b_nT^n$ is the leading term of of $G$.  Indeed, as $w(f)<w(f-G)$, we
have $w(G)=w(f)$.  If we also had $w(G)<w(b_nT^n)$, then we would have
$w((a_n-b_n)T^n)=w(f)$~; but $w(f-G)\le w((a_n-b_n)T^n)=w(f)$ is a
contradiction.  Writing
$$
f-GH=qG+r,\quad\deg(q)\le m-n,\ \ \deg(r)<n,
$$
lemma~8 implies that $w(q)\ge\delta$ and $w(r)\ge w(f)+\delta$.  Now take
$G'=G+r$, $H'=H+q$.  Clearly $w(f-G')\ge w(f)+\Delta$ and $w(H'-1)\ge\Delta$.
It is also clear that $G',H'$ have degrees $n$ and $\le m-n$ respectively, so
they satisfy all the conditions that $G,H$ were required to satisfy.  But now,
$$\eqalign{
w(f-G'H')&=w((H-1)r+qr)\cr
&\ge\inf(w(H-1)+w(r),w(q)+w(r))\cr
&\ge w(f)+\delta+\Delta,\cr
}
$$
so $G'H'$ is a ``better approximation'' to $f$ than $GH$.  If $f\neq G'H'$,
the process can be repeated, to get an even better approximation.  In the
limit we get $g,h\in K[T]$ of degree $n$ and $m-n$ such that $f=gh$, as
claimed.  (This is the only place where the completeness of $K$ has been
used.) 

\th COROLLARY 10
\enonce
If\/ $f$ is an irreducible polynomial over a complete unarchimedean field,
then\/ $f$ is pure.  
\endth
The converse is of course not true~; see lemma~6.

\th COROLLARY 11
\enonce
If\/ $f=a_0+a_1T+\cdots+a_nT^n$ is pure and has degree~$n$, then
$v(a_i)\ge\inf(v(a_0),v(a_n))$.  In particular,  and if\/
$a_0,a_n\in\ogoth$, then $f\in\ogoth[T]$.  
\endth
As $\Pi_f$ is the line segment from $(0,v(a_0))$ to $(n,v(a_n))$, the
coefficients of $f$ have to have valuation at least $\inf(v(a_0),v(a_n))$, and
$f\in\ogoth[T]$ if $a_0,a_n\in\ogoth$. 

\medskip

{\it Proof of th.} 5\pointir Let $f=h_1h_2\ldots h_s$ be the factorisation of
$f$ into irreducible polynomials (with repetition), and let $\eta_i$ be the
slope of $h_i$, which is pure by cor.~10.  Take the smallest slope $\delta_1$
among the $\eta_i$ and let $g_1$ be the product of the $h_i$ with slope
$\delta_1$, so $g_1$ is pure of slope $\delta_1$ (lemma~6) and type
$(m_1,\delta_1)$, say.  Next, let $\delta_2>\delta_1$ be the smallest slope
among the remaining $\eta_i$ and let $g_2$ be the product of the $h_i$ with
slope $\delta_2$, so $g_2$ is pure of slope $\delta_2$, and say of type
$(m_2,\delta_2)$.  This process comes to an end after a certain number $t$ of
steps, so that we get a factorisation $f=g_1g_2\ldots g_t$, where each $g_j$
is of type $(m_j,\delta_j)$, and $\delta_1<\delta_2<\cdots<\delta_t$.
Consequently, $g_1g_2\ldots g_t$ is of type
$(m_1,\delta_1;m_2,\delta_2;\cdots;m_t,\delta_t)$, by repeated application of
lemma~7.  But the type of $f$ is
$(l_1,\gamma_1;l_2,\gamma_2;\ldots;l_r,\gamma_r)$, so the two types must be
the same, and the $g_j$ have been found.\hfil QED

\medskip
{\it Example}\pointir The polynomial
$1+{1\over1!}T+\cdots+{1\over7!}T^7$ has a root in $\Q_2$.
\medskip
Note finally that, at least when $v$ is discrete, we know how to factor a
polynomial in $\ogoth[T]$ if an approximate factorisation is given
(Lecture~7).  If moreover the residue field $\ogoth/\pgoth$ is finite, an
approximate factorisation can always be found by trial and error, as there are
only finitely many polynomials of a given degree over $\ogoth/\pgoth^n$.

\vfill\eject

\centerline{\bf Lecture 10}
\bigskip
\centerline{\it Factorisation~; Weierstra\ss\ preparation}
\bigskip

Let $K$ be a field complete for a valuation $v$ and denote by $\ogoth$ the
ring of $v$-integers.  Extend $v$ to a valuation $w$ on $K(T)$ by requiring
that $w(f)=\inf_i v(a_i)$ for every polynomial $f=\sum a_iT^i$ in $K[T]$
(Lecture~8, lemma~1).  Our first aim is to show how a sufficiently good
approximate factorisation of $f\in\ogoth[T]$ can be converted into an actual
factorisation~; see Lecture~7, theorem~2, when $v$ is discrete.

\th THEOREM 1
\enonce
Let\/ $f,G,H\in\ogoth[T]$ be polynomials of degrees\/ $n,s,t$ respectively such
that\/ $f$ and\/ $GH$ have the same leading terms (so that $n=s+t$).  Suppose
that\/ $w(f-GH)>2v(\res(G,H))$.  Then there exist\/ $g,h\in\ogoth[T]$,
having the same leading terms as\/ $G,H$, and such that\/ 
$$
w(G-g), w(H-h)\ge w\!\left({f-GH\over\res(G,H)}\right),\qquad f=gh.
$$
\endth
We look for polynomials $\gamma,\delta\in\ogoth[T]$ of degree $<s$
(resp.~$<t$) such that $G'H'$ is a better approximation to $f$ than $GH$,
where $G'=G+\gamma$ and $H'=H+\delta$.

We have $f-G'H'=(f-GH)-G\delta-H\gamma-\gamma\delta$.  Let
$\gamma,\delta$ be such that $(f-GH)-G\delta-H\gamma=0$~; let's see
why they are uniquely determined by this condition.

Think of the requirement $G\delta+H\gamma=f-GH$ as a system of equations in
which the coefficients of $\delta$ and $\gamma$ are the unknowns.  The general
idea is illustrated by the particular case $s=2$, $t=3$~; denoting the
coefficients of $G,H,f-GH,\gamma,\delta$ respectively by
$g_i,h_i,b_i,c_i,d_i$, the system of equations $G\delta+H\gamma=f-GH$ can be
written matricially as
$$
\pmatrix{
g_2&0&0&h_3&0\cr
g_1&g_2&0&h_2&h_3\cr
g_0&g_1&g_2&h_1&h_2\cr
0&g_0&g_1&h_0&h_1\cr
0&0&g_0&0&h_0\cr
}
\pmatrix{
d_2\cr
d_1\cr
d_0\cr
c_1\cr
c_0\cr
}
=
\pmatrix{
b_4\cr
b_3\cr
b_2\cr
b_1\cr
b_0\cr
}
.$$
The determinant of the square matrix is by definition $\res(G,H)$, which is
$\neq0$ by hypothesis.  Thus for example $d_2=d/\res(G,H)$, where
$$
d=\left|
\matrix{
b_4&0&0&h_3&0\cr
b_3&g_2&0&h_2&h_3\cr
b_2&g_1&g_2&h_1&h_2\cr
b_1&g_0&g_1&h_0&h_1\cr
b_0&0&g_0&0&h_0\cr
}\right|.
$$
Expanding along the first column, we see that $v(d)\ge w(f-GH)$.  The other
unknowns $d_i,c_i$ are given by similar formulas, and we conclude that
$w(\delta)$ and $w(\gamma)$ are both $\ge w(f-GH)-v(\res(G,H))$~; in
particualr, $\delta,\gamma\in\ogoth[T]$.

To see that $G'H'$ is a better approximation to $f$ than $GH$, note that 
$$
\eqalign{
w(f-G'H')&=w(\gamma\delta)\cr
&\ge 2w(f-GH)-2v(\res(G,H))\cr
&>w(f-GH).\cr
}$$
Moreover, the above estimates imply that $w(\gamma),w(\delta)>v(\res(G,H))$~; 
it follows (check it) that $v(\res(G',H'))=v(\res(G,H))$, and the whole
argument can be repeated, to get two polynomials $f,g\in\ogoth[T]$ in the
limit. We thus get the factorisation $f=gh$, as claimed.

\th COROLLARY 2
\enonce
The same conclusion holds if\/ $2v(\res(G,H))$ is replaced by\/
$v(\dis(f))$ in the hypothesis\/ $w(f-GH)>2v(\res(G,H))$.
\endth
First check that $v(\dis(GH))=v(\dis(f))$.  The discriminant of the product
$GH$ is given by $\dis(GH)=\dis(G)\dis(H)\res(G,H)^2$.  In view of the fact
that $G,H\in\ogoth[T]$, the discriminants $\dis(G)$ and $\dis(H)$ are in
$\ogoth$~; it follows that $v(\dis(f))=v(\dis(GH))\ge2v(\res(G,H))$, and we
are done.

\bigskip
\centerline{\hbox to5cm{\hrulefill}}
\bigbreak

Let $W\subset K[[T]]$ be the subring of those $h=\sum_ic_iT^i$ for which
$\lim_{i\to+\infty}c_i=0$, and, for $h\in W$, define $w(h)=\inf_iv(c_i)$~; in
particular, $w(0)=+\infty$ and $w(T)=0$.  It is easily seen that
$w(h_1h_2)=w(h_1)+w(h_2)$ and that $w(h_1+h_2)\ge\inf(w(h_1),w(h_2))$.  Thus
$w$ is a valuation on the field of fractions of $W$.

To say that a power series in $f=\sum_ia_iT^i$ in $W$ has $w(f)>0$ means that
the $a_i$ are in the maximal ideal of the ring of $v$-integers of $K$.
Suppose that $f\neq0$, and define $N$ by the requirement
$$
v(a_N)=w(f)\;;\qquad v(a_i)>v(a_N)\quad (i>N).
$$
\vskip-2\baselineskip
\th THEOREM 3 (``\thinspace Weierstra\ss\thinspace'')
\enonce
There exists a power series\/ $h=\sum_{i}c_iT^i$ in\/ $W$, with\/ $h(0)=1$,
$w(h-1)>0$, and a polynomial\/ $g=b_0+b_1T+\cdots+b_NT^N$ of degree\/~$N$ in\/
$K[T]$, with $v(b_N)=w(g)$, such that\/ $f=gh$.
\endth

\th LEMMA 4
\enonce
The ring\/ $W$ is complete for the distance defined by\/ $w$, and $K[T]$ is
a dense subring. 
\endth
Let $h_n=\sum_ic_i^{(n)}T^i$ ($n>0$) be a fundamental sequence of elements of 
$W$.  For each $i$, the sequence $c_i^{(n)}$ ($n>0$) is a fundamental sequence
in $K$, which is complete~; let $d_i$ be the limit of this sequence, and put
$h=\sum_id_iT^i$.  We will show that $h=\lim_{n\to+\infty}h_n$ is in $W$.

Let a real constant $R$ --- however large --- be given.  There is some index
$\bar n$ such that $w(h_n-h_{\bar n})>R$ for all~$n>\bar n$.  So
$v(d_i-c_i^{(\bar n)})>R$ for all~$i$.  Since $h_{\bar n}$ is in $W$, there is
some index $\bar\imath$, depending upon $R$, such that $v(c_i^{(\bar n)})>R$
for all $i>\bar\imath$.  It follows that $v(d_i)>R$ for all $i>\bar\imath$.
But $R$ is arbitrary, so $\lim_{i\to+\infty}d_i=0$, and hence $h\in W$, as
claimed.

Finally, $K[T]$ is dense in $W$ because every $f\in W$ is the limit of the
sequence $(f_n)_n$, where $f_n\in K[T]$ is the sum of the first $n+1$ terms
$a_0,a_1T,\ldots,a_nT^n$ of $f$.

\medskip

Let $g=b_0+b_1T+\cdots+b_NT^N$ be such that $v(b_N)=\inf_iv(b_i)=w(g)$. 
\th LEMMA 5
\enonce
For every\/ $\varphi\in W$, there is a pair\/ $q\in W$, $r\in K[T]$ consisting
of a power series and a polynomial of degree\/ $<N$ such that\/ $\varphi=qg+r$.
Moreover,   
$$
w(qg)\ge w(\varphi),\quad w(r)\ge w(\varphi).
$$
\endth
Let $\varphi_i\in K[T]$ ($i>0$) be a (fundamental) sequence of polynomials
such that $\varphi=\lim_{i\to+\infty}\varphi_i$.   Write $\varphi_i=q_ig+r_i$,
with $\deg(r_i)<N$.  We have
$$
(\varphi_j-\varphi_i)=(q_j-q_i)g+(r_j-r_i),
$$
to which we may apply Lemma~8, Lecture~9~; it shows that the sequences
$(q_i)_i$, $(r_i)_i$ are fundamental.  Taking $q=\lim_{i\to+\infty}q_i$,
$r=\lim_{i\to+\infty}r_i$ proves the lemma.

{\it Proof of th.} 3\pointir The proof is similar to that of Lemma~9,
Lecture~9.  Let $f$ and $N$ be given.  Letting $f_N=a_0+a_1T+\cdots+a_NT^N$,
we have $w(f-f_N)>w(f)$ by hypothesis, so there exists a real $\Delta>0$ such
that $w(f-f_N)=w(f)+\Delta$.

Consider a polynomial $G\in K[T]$ of degree~$N$ and a power series $H\in W$
such that  
$$
w(f-G)\ge w(f)+\Delta,\quad H(0)=1,\quad w(H-1)\ge\Delta.
$$
Adding $w(G)$ to the second inequality, we get $w(GH-G)\ge w(G)+\Delta$, and
hence $w(f-GH)\ge w(f)+\Delta$.  Defining $\delta$ by $w(f-GH)=w(f)+\delta$,
we have $\delta\ge\Delta$. 

One choice of $G,H$ would be $f_N,1$.  We shall show that if $f\neq GH$,
equivalently if $\delta<+\infty$, then we can find $G',H'$ which satisfy the
same conditions as $G,H$ but for which $\delta'\ge\delta+\Delta$.

It follows from the hypotheses on $f$, $N$ and $G$ that $w(G)=v(b_N)$.
Indeed, as $w(f)<w(f-G)$, we have $w(G)=w(f)$.  If we also had
$w(G)<v(b_N)$, then we would have $v(a_N-b_N)=w(f)$~; but $w(f-G)\le
v(a_N-b_N)=w(f)$ is a contradiction.  Using Lemma~5, write
$$
f-GH=qG+r,\quad q\in W,\ \ \deg(r)<N,
$$
so that $w(q)\ge\delta$ and $w(r)\ge w(f)+\delta$.  Now take $G'=G+r$,
$H'=H+q$.  Clearly $w(f-G')\ge w(f)+\Delta$ and $w(H'-1)\ge\Delta$.  It is
clear that $G'$ has degree $N$, so the pair $G',H'$
satisfies all the conditions that $G,H$ was required to satisfy, with the
possible exception of $H'(0)=1$.  But now,
$$\eqalign{
w(f-G'H')&=w((H-1)r+qr)\cr
&\ge\inf(w(H-1)+w(r),w(q)+w(r))\cr
&\ge w(f)+\delta+\Delta,\cr
}
$$
so $G'H'$ is a ``better approximation'' to $f$ than $GH$.  If $f\neq G'H'$,
the process can be repeated, to get an even better approximation.  In the
limit we get $g\in K[T]$ of degree $N$, with valuation equal to that of its
leading term, and $h\in W$, with $h(0)$ a unit in $K$ and $w(h-1)>0$, such
that $f=gh$.  To meet the requirement $h(0)=1$, multiply $g$ by $h(0)$ and $h$
by $h(0)^{-1}$.

\vfill\eject

\centerline{\bf Lecture 11}
\bigskip
\centerline{\it Valuations on algebraic extensions of complete fields}
\bigskip

Let $K$ be a field with a valuation $v$.  We have seen some ways of extending
$v$ to a valuation $w$ on purely transcendental extensions of $K$ (Lecture~8).

Suppose that $K$ is complete for $v$, and let $L|K$ be a finite extension.  We
shall see that there is a unique valuation $w$ on $L$ extending $v$, and that
$L$ is complete for $w$.

\th THEOREM 1
\enonce
Suppose that the field\/ $K$ is complete for the valuation\/ $v$, and let\/
$L$ be a finite extension of\/ $K$.  Then\/ $w(x)=v(N_{L|K}(x))/[L:K]$ is a
valuation on\/ $L$ extending\/ $v$.  The field\/ $L$ is complete for\/ $w$,
and\/ $w$ is the only extension of\/ $v$ to\/ $L$.   
\endth

Let $K$ be a field with an absolute value $|\ |$, and let $E$ be a vector
$K$-space.

\th DEFINITION 2
\enonce
A norm on\/ $E$ is a map\/ $|\!|\ |\!|:E\to[0,+\infty[$ satisfying
$$
|\!|x|\!|=0\Leftrightarrow x=0,\quad |\!|x+y|\!|\le|\!|x|\!|+|\!|y|\!|,\quad
|\!|\alpha x|\!|=|\alpha||\!|x|\!|,
$$
for all\/ $x,y\in E$ and every\/ $\alpha\in K$.
\endth
A norm defines a metric $(x,y)\mapsto|\!|x-y|\!|$, and hence a topology,
on~$E$.  
\th DEFINITION 3
\enonce
Two norms $|\!|\ |\!|_1,|\!|\ |\!|_2$ on\/ $E$ are said to be equivalent if
there exist constants\/ $C_1>0$, $C_2>0$ such that, for all\/ $x\in E$,
$$
|\!|x|\!|_1\le C_1|\!|x|\!|_2,\qquad |\!|x|\!|_2\le C_2|\!|x|\!|_1.
$$
\endth
Two norms are equivalent if and only if they define the same topology.  Every
norm on $K$ is equivalent to $|\ |$.

\th LEMMA 4
\enonce
Suppose that\/ $K$ is complete for\/ $|\ |$, and that\/ $E$ is
finite-dimensional.  Then there are norms on\/ $E$, any two norms are
equivalent, and\/ $E$ is complete for every norm.
\endth
Let $(e_\alpha)_{\alpha\in A}$ be a $K$-basis of $E$ and, for every
$x=\sum_{\alpha\in A}\xi_\alpha e_\alpha$,  define
$$|\!|x|\!|_1=\sup_{\alpha\in A}|\xi_\alpha|.$$ 
Clearly this is a norm on $E$ and $E$ is complete for $|\!|\ |\!|_1$.  It
suffices to show that every norm $|\!|\ |\!|$ on $E$ is equivalent to this
one.  We have
$$
|\!|x|\!|=\bigg|\!\bigg|\sum_{\alpha\in A}\xi_\alpha e_\alpha\bigg|\!\bigg|
\le\sum_{\alpha\in A}|\xi_\alpha|\,|\!|e_\alpha |\!|\le C|\!|x|\!|_1,
$$
where $C=\sum_{\alpha\in A}|\!|e_\alpha|\!|$.  It remains to find a constant
$C_1>0$ such that $|\!|x|\!|_1\le C_1|\!|x|\!|$ for every $x\in E$.

Let $\varepsilon_n>0$ ($n\in\N$) be a sequence of reals tending to~$0$.  If
$C_1$ does not exist, then, for every $n\in\N$, there would exist a $y_n\in E$
such that $|\!|y_n|\!|<\varepsilon_n|\!|y_n|\!|_1$.  For each
$y_n=\sum_{\alpha\in A}\eta_{n,\alpha} e_\alpha$, there is an $\alpha\in A$
such that $|\!|y_n|\!|_1=|\eta_{n,\alpha}|$.  As the set $A$ is finite and
$\N$ infinite, there is some $\beta\in A$ such that
$|\!|y_n|\!|_1=|\eta_{n,\beta}|$ for infinitely many $n$~; restricting
ourselves to these $n$, we may assume that this holds for {\it all\/} $n\in\N$
(and the same $\beta\in A$).  Replacing $y_n$ by $\eta_{n,\beta}^{-1}y_n$, we
have $y_n=z_n+e_\beta$, where $z_n$ is in the hyperplane $E'\subset E$
generated by the $e_\alpha$ ($\alpha\neq\beta$).

Thus, if $C_1$ does not exist, we can find a sequence $z_n$ ($n>0$) of
elements of $E'$ such that $|\!|z_n+e_\beta|\!|\rightarrow0$ as $n\to+\infty$.
This implies that $|\!|z_m-z_n|\!|$ is as small as we please if $m,n$ are
sufficiently large~: the sequence $(z_n)_{n\in\N}$ is fundamental.  By the
inductive hypothesis the space $E'$, which is of dimension $<\dim E$, is
complete for $|\!|\ |\!|$, so $t=\lim_{n\to+\infty}z_n$ exists in $E'$.

But then $|\!|t+e_\beta|\!|=\lim_{n\to+\infty}|\!|z_n+e_\beta|\!|=0$, whereas
$t+e_\beta\neq0$, in contradiction to the requirement
$|\!|x|\!|=0\Leftrightarrow x=0$.  Our assumption that $C_1$ does not exist
must therefore have been wrong.

\medskip

{\it Proof of th.} 1\pointir If two valuations $w_1,w_2$ on $L$ extend the
valuation $v$ on $K$, then the corresponding absolute values are equivalent as
norms (Lemma~4) and hence induce the same topology on $L$.  So $w_1,w_2$ are
equivalent as valuations (Lemma~2, Lecture~4)~; as they agree on $K$, they
must be equal.

It remains to show that $w(x)=v(N_{L|K}(x))/n$ is a valuation on $L$, where
$n=[L:K]$.  It is clear that $w$ extends $v$, for $N_{L|K}(x)=x^n$ for every
$x\in K$.  As $N_{L|K}:L^\times\to K^\times$ is a homomorphism and
$N_{L|K}(0)=0$, it is equally clear that $w(xy)=w(x)+w(y)$ and
$w(x)=+\infty\Leftrightarrow x=0$.

Let us show that $w(x+y)\ge\inf(w(x),w(y))$.  We have
$w(x+y)=w(y)+w(1+xy^{-1})$.  Suppose for instance that $w(x)\ge w(y)$.  By
taking $z=xy^{-1}$, it suffice to show that if $w(z)\ge0$, then $w(1+z)\ge0$,
for $z\in L$.

Let $f\in K[T]$ be the characteristic polynomial of $z$~; as
$(-1)^nf(0)=N_{L|K}(z)$, we have $f(0)\in\ogoth$ by hypothesis.  We also have
$f=g^r$, where $g\in K[T]$ is the minimal polynomial of $z$ and $r=[L:K(z)]$.
As $g$ is irreducible, it is pure (Lecture~9, Corollary~10), and hence so is
$f$ (Lecture~9, Lemma~6).  As the constant term $f(0)$ and the leading
coefficient~$1$ of $f$ are in $\ogoth$, we have $f\in\ogoth[T]$ (Lecture~9,
Corollary~11).  As the characteristic polynomial of $1+z$ is $f(T-1)$, we have
$N_{L|K}(1+z)=(-1)^nf(-1)$, and therefore $w(1+z)\ge0$, which was to be
proved.

\th COROLLARY 5
\enonce
The ring of\/ $w$-integers is the integral closure of\/ $\ogoth$ in\/ $L$.
\endth
We have already seen that if $z\in L$ is $w$-integral, then its characteristic
polynomial is in $\ogoth[T]$, so $z$ is integral over $\ogoth$.  Conversely,
if $z\in L$ satisfies $z^n+a_{n-1}z^{n-1}+\cdots+a_0=0$ with $a_i\in\ogoth$,
then 
$$w(z^n)=nw(z)\ge\inf_{i\in[0,n[}(v(a_i)+iw(z))\ge\inf(0,(n-1)w(z)),
$$
hence $w(z)\ge0$.

\th COROLLARY 6
\enonce
Let\/ $\bar K$ be an algebraic closure of the complete field\/ $K$.
There is a unique valuation\/ $\bar v$ on\/ $\bar K$ which extends\/ $v$.
\endth
For $\bar K$ is the union of the finite extensions $L\subset\bar K$ of $K$. 

\th COROLLARY 7
\enonce
If\/ $x,x'\in\bar K$ are conjugate over\/ $K$, then\/ $\bar v(x)=\bar v(x')$.
\endth
Let $\sigma$ be a $K$-automorphism of $\bar K$ such that $\sigma(x)=x'$.
Notice that $\bar v\circ\sigma$ is a valuation on $\bar K$ extending $v$, so
$\bar v\circ\sigma=\bar v$ (Corollary~6).  In particular, $\bar v(x)=\bar
v(\sigma(x))=\bar v(x')$.  Alternatively, let $L\subset\bar K$ be a finite
extension of $K$ containing $x,x'$.  As $x,x'$ have the same minimal
polynomial over $K$, they have the same characteristic polynomial (as
endomorphisms of the $K$-space $L$), and $N_{L|K}$ agrees on $x,x'$.

\th COROLLARY 8
\enonce
If\/ $x,x'\in\bar K$ are $K$-conjugate, then\/ $\bar v(x-x')\ge\bar
v(x-\alpha)$ for every $\alpha\in K$ (even if\/ $x\neq x'$).
\endth
If we had $\bar v(x-x')<\bar v(x-\alpha)$, then $\bar v(x'-\alpha)=\bar
v(x-x')<\bar v(x-\alpha)$, which is impossible by the previous corollary as
$x'-\alpha$ and $x-\alpha$ are conjugate over $K$. 

\th COROLLARY 9 (Krasner)
\enonce
Let $x\in\bar K$ be separable of degree\/ $n$ over\/ $K$, and let\/ $x=x_1, 
x_2,\ldots,x_n$ be its conjugates.  If\/ $\bar v(y-x)>\sup_{i\in[2,n]}\bar
v(x-x_i)$, then\/ $K(x)\subset K(y)$.   
\endth
Notice that $\bar v(y-x_i)=\bar v(x-x_i)$ for $i\in[2,n]$.  Let $\sigma$ be
any $K(y)$-automorphism of $\bar K$~; we have $\sigma(y)=y$ and
$\sigma(x)=x_i$ for some $i\in[1,n]$.  As $\bar v(y-x)=\bar v(y-x_i)$
(Corollary~7), the only possibility is $i=1$, in view of the hypothesis $\bar
v(y-x)>\bar v(y-x_i)$, $i\in[2,n]$.  So $x$ is fixed by every $\sigma$. Since
moreover $x$ is separable over $K(y)$, we have $x\in K(y)$, as desired.

\vfill\eject

\centerline{\bf Lecture 12}
\bigskip
\centerline{\it The algebraic closure of the completion}
\bigskip

Let $K$ be a complete unarchimedean field, $v$ its valuation, $\bar K$ an
algebraic closure of $K$, and $\bar v$ the valuation on $\bar K$ extending
$v$.  We will show that $\bar K$ need not be complete, but that its completion
$\hat{\bar K}$ is algebraically closed.

Extend the valuation $v$ to a valuation $w$ on $K[T]$ by requiring that
$w(T)=0$, and that the valuation of every polynomial is the infimum of the
valuations of its terms (Lecture~8, Lemma~1).

\th PROPOSITION 1
\enonce
Let\/ $f\in K[T]$ be a separable irreducible unitary polynomial of degree~$n$,
and let\/ $g\in K[T]$ be any unitary polynomial of the same degree.  If\/ $g$
is sufficiently close to\/ $f$, then\/ $g$ is separable, irreducible, and
defines the same extension of\/ $K$ as\/ $f$.
\endth

For the time being, let $f=a_0+\cdots+a_{n-1}T^{n-1}+T^n$ be any unitary
polynomial of degree~$n$, and let $\alpha_1,\alpha_2,\ldots,\alpha_n$ be the
roots of $f$ in $\bar K$.  Suppose that $w(f)\ge A$ for some $A\le0$.  If
$\gamma\in\bar K$ is such that $\bar v(\gamma)<A$, then, writing
$f(\gamma)=a_0+\cdots+a_{n-1}\gamma^{n-1}+\gamma^n$, we see that $\bar
v(f(\gamma))=\bar v(\gamma^n)$, for the valuation of $\gamma^n$ is strictly
less than the valuations of the other terms.  In particular, $f(\gamma)\neq0$.
It follows that $\bar v(\alpha_i)\ge A$ for every $i\in[1,n]$.

Now let $g=b_0+\cdots+b_{n-1}T^{n-1}+T^n$ be another unitary polynomial of
degree~$n$ in $K[T]$, and suppose that $w(g)\ge A$.  Let $\beta$ be a root of
$g$, so that $\bar v(\beta)\ge A$ by the foregoing.  Let $C>0$ be a constant
such that $w(f-g)>C$.  Then the equalities
$$
\prod_{i\in[1,n]}(\beta-\alpha_i)
=f(\beta)=f(\beta)-g(\beta)=\sum_{j\in[0,n]}(a_j-b_j)\beta^i
$$
taken together impliy that $\sum_{i\in[1,n]}\bar v(\beta-\alpha_i)\ge nA+C$.
Therefore there is some index $i\in[1,n]$ such that $\bar v(\beta-\alpha_i)\ge
A+C/n$.  This can be expressed by saying that by taking $g$ sufficiently close
to $f$ --- by taking $C$ to be large enough --- every root $\beta$ of $g$ can
be brought as close to some root (depending on $\beta$) $\alpha_i$ of $f$ as
we wish.

Let us take $g$ so close to $f$ that every root $\beta$ of $g$ is closer to
some root (depending on $\beta$) $\alpha_i$ of $f$ than $\alpha_i$ is to
$\alpha_j$ for $\alpha_j\neq\alpha_i$.  Suppose further that $f$ is separable
(so that the roots of $f$ are simple) and irreducible.  Then, since $\bar
v(\beta-\alpha_i)>\sup_{j\neq i}\bar v(\alpha_i-\alpha_j)$, we have
$K(\alpha_i)\subset K(\beta)$ (Lecture~11, Corollary~9).  But $f$ and $g$ have
the same degree, so $K(\alpha_i)=K(\beta)$, and the proof is over.  \medbreak

It may happen that $\bar K$ is not complete.  For example, when
$K=\Q_2$, the limit
$$
1+2\sqrt2+4\root4\of2+8\root8\of2+\cdots
$$
does not exist in $\bar\Q_2$ (E. Artin).  But

\th THEOREM 2
\enonce
The completion\/ $C$ of\/ $\bar K$ is algebraically closed.
\endth
Let $f\in C[T]$ be an irreducible unitary polynomial, and suppose first that
$f$ is separable. Let $g\in\bar K[T]$ be so close to $f$ that they define
the same extension of $C$ (Proposition~1).  But $g$ has a root in $\bar K$, and
hence in $C$, so $f$ has a root in $C$.  So $C$ is separably closed.  If the
characteristic of $K$ is $0$, the proof is over.

If the characteristic is $p\neq0$, it remains to show that $C$ is perfect. If
$\alpha\in C$ is the limit of a (fundamental) sequence $(\alpha_n)_n$ in $\bar
K$, then, because $\root p\of{\alpha_m}-\root p\of{\alpha_n}=\root
p\of{\alpha_m-\alpha_n}$, the sequence $(\root p\of{\alpha_n})_n$ is also
fundamental, and its limit $\beta\in C$ satisfies $\beta^p=\alpha$.  Therefore
the field $C$ is perfect.

\bigskip
\centerline{\hbox to5cm{\hrulefill}}
\bigbreak

Let $K$ be a field complete for a valuation $v$, $\bar K$ an algebraic closure
of $K$, and $\bar v$ the valuation on $\bar K$ extending $v$.  Recall that a
polynomial $f\in K[T]$ of degree $n>0$ is called pure of type $(n,\gamma)$ if
$f(0)\neq0$ and if the Newton polygon $\Pi_f$ of $f$ is a line segment of
slope $\gamma$.

\th PROPOSTION 3
\enonce
If\/ $f\in K[T]$ is pure of type\/ $(n,\gamma)$, then\/ $\bar v(a)=-\gamma$
for every root\/ $a\in\bar K$ of\/ $f$.
\endth
As $cf$ has the same type and the same roots as $f$ for every $c\in K^\times$,
we may assume that $f$ is unitary.  Let $f=(T-a_1)\cdots(T-a_n)$ be the
factorisation of $f$ over $\bar K$, and $(1,\gamma_i)$ the type of $T-a_i$.
If we had $\gamma_i<\gamma_j$ for some $i,j\in[1,n]$, then $f$ would not be
pure (Lemma~7, Lecture~9), contrary to the hypothesis.  So the $\gamma_i$ are
all equal to some $\delta$, and the type of $f$ is $(n,\delta)$ (Lemma~6,
Lecture~9), forcing $\delta=\gamma$. Finally, $\bar v(a_i)=-\gamma_i=-\gamma$.

\th COROLLARY 4
\enonce
If\/ $f\in K[T]$ is of type\/ $(l_1,\gamma_1;
l_2,\gamma_2;\ldots;l_r,\gamma_r)$, then\/ $f$ has\/ $l_i$ roots (counted with
multiplicity)\/ $a\in\bar K$ with\/ $\bar v(a)=-\gamma_i$ ($i\in[1,r])$.
\endth

\bigskip
\centerline{\hbox to5cm{\hrulefill}}
\bigbreak

Let $K$ be a complete unarchimedean field, $v$ its valuation, $L$ a finite
extension of $K$ of degree~$n$, and $w$ the unique extension of $v$ to a
valuation on $L$ (Lecture~11, Theorem~1).

Let $\ogoth$ be the ring of $v$-integers, $\pgoth$ its unique maximal ideal,
and similarly define $\goth O,P$ relative to $w$.  We have the inclusions
$\goth o\subset O$, $\goth p\subset P$ and
$$
{\goth o}={\goth O}\cap K\;;\qquad {\goth p}={\goth P}\cap K.
$$
If $k=\ogoth/\pgoth$ and $l={\goth O/P}$ denote the residue fields of $K$
and $L$, we get an map $k\to l$, making $l$ into an extension of $k$, called
the {\it residual extension\/} of $L|K$, and its degree $f_{L|K}=[l:k]$ is
called the {\it residual degree\/} of $L|K$.  It is clear that if $M$ is a
finite extension of $L$, then $f_{M|K}=f_{M|L}f_{L|K}$.  We say that $L|K$ is
{\it unramified\/} if $f_{L|K}=[L:K]$ and if the residual extension is {\it
  separable\/}. We say that $L|K$ is {\it totally ramified\/} if $f_{L|K}=1$.  

\th LEMMA 5
\enonce
We have\/ $f_{L|K}\le[L:K]$.  
\endth
Denote the map ${\goth O}\to l$ by $\alpha\mapsto\bar\alpha$.  Let
$\bar\alpha_1,\bar\alpha_2,\ldots,\bar\alpha_{n+1}$ be in $l$, where
$n=[L:K]$.  Among the elements $\alpha_1,\alpha_2,\ldots,\alpha_{n+1}$ of
$\goth O$, there is a linear relation
$\xi_1\alpha_1+\xi_2\alpha_2+\cdots+\xi_{n+1}\alpha_{n+1}$ ($\xi_i\in K$, not
all $0$).  Multiplying by a suitable $\eta\in K^\times$, we may assume that
$\xi_i\in\ogoth$, and that at least one of them is in $\ogoth^\times$.
Reducing modulo $\goth P$, we get the relation
$\bar\xi_1\bar\alpha_1+\bar\xi_2\bar\alpha_2+
\cdots+\bar\xi_{n+1}\bar\alpha_{n+1}$ ($\bar\xi_i\in k$, not all $0$), showing
that every family of $n+1$ elements in $l$ is linearly dependent over $k$.

\th PROPOSTION 6
\enonce
The inclusion-preserving map which sends an extension of\/ $K$ in\/ $L$ to its
residue field\/ induces a bijection between the set of unramified extensions
of\/ $K$ in\/ $L$ and the set of separable extensions of\/ $k$ in\/ $l$.
\endth
Every separable extension of $k$ in $l$ is of the form $k(a)$ for some $a\in
l$ (separable over $k$).  We show first that for any lift $\tilde a\in L$ of
$a$, there is an unramified extension $K'\subset K(\tilde a)$ whose residue
field is $k(a)$. 

Let $f\in\ogoth[T]$ be a unitary polynomial such that $\bar f\in k[T]$ is the
minimal polynomial of $a$~; by the separability hypothesis, $a$ is a {\it
  simple\/} root of $\bar f$.  Hensel's lemma (Lecture~6, Corollary~2) shows
that $f$ has a simple root $\alpha\in K(\tilde a)$ such that $\bar\alpha= a$.
Clearly, $[K(\alpha):K]=[k(a):k]$, for $f$ and $\bar f$ are irreducible (over
$K$ and $k$ respectively) of the same degree~; for $f$, use Corollary~4,
Lecture~8.  Also, the residue field of $K'=K(\alpha)$ contains $k(a)$, so must
be equal to it (Lemma~5).

Let $K''\subset L$ be another unramified extension whose residue field is
$k(a)$.  Then, in the foregoing construction, we could have taken $\tilde a$
to be in $K''$~; doing so would give the same $K'$, for $\alpha$ is the unique
lift of $a$.  We thus have $K'\subset K''$, and, as these two extensions have
the same degree, they must be equal.  There is thus a unique unramified
extensions of $K$ in $L$ with residue field $k(a)$.

\smallbreak

In particular, there is a largest unramified extension $L_0\subset L$~; its
residue field is the separable closure $l_0$ of $k$ in $l$.

\medbreak
\centerline{\hbox to5cm{\hrulefill}}
\bigbreak

We have $v(K^\times)\subset w(L^\times)$~; the index $e_{L|K}={\bf
  (}w(L^\times):v(K^\times){\bf)}$ is called the {\it ramification index\/} of
$L|K$.  It is clearly finite and divides $n$ if $v$ is discrete, for
$w(L^\times)\subset {1\over n}v(K^\times)$~; we will show that it is finite in
general.  It is clear that if $M$ is a finite extension of $L$, then
$f_{M|K}=f_{M|L}f_{L|K}$.

\th PROPOSTION 7
\enonce
We have $e_{L|K}f_{L|K}\le[L:K]$~; in particular, $e_{L|K}$ is finite.
\endth
Write more simply $e,f$ for $e_{L|K},f_{L|K}$, and let $e'\in[0,e[$ be any
integer.  Let $\omega_1,\omega_2,\ldots,\omega_f\in{\goth O}$ be lifts of a
$k$-basis of $l$, and let $\pi_0,\pi_1,\ldots,\pi_{e'}\in L^\times$ have
distinct images in $w(L^\times)/v(K^\times)$.  We shall show that the
$(e'+1)f$ elements $\pi_i\omega_j$ ($i\in[0,e'],j\in[1,f]$) of $L$ are
linearly independent over $K$.  It would follow that $e$ is finite, and,
taking $e'=e-1$, that $ef\le[L:K]$.

Suppose that we have $\sum_{i\in[0,e']}\sum_{j\in[1,f]}a_{ij}\pi_i\omega_j=0$
for some $a_{ij}\in K$, not all $0$.  Writing
$s_i=\sum_{j\in[1,f]}a_{ij}\omega_j$, we see that not all $s_i$ are $0$.
Indeed, we may assume that all $a_{ij}\in\ogoth$ and at least one
$a_{nm}\in\ogoth^\times$.  Modulo $\goth P$, we have $\bar
s_n=\sum_{j\in[1,f]}\bar a_{nj}\bar\omega_j$, and as $\bar a_{nm}\neq0$, we
have $\bar s_n\neq0$, the $\bar\omega_j$ being $k$-linearly independent.

Moreover, whenever $s_i\neq0$, we have $w(s_i)\in v(K^\times)$.  Indeed, let
$d_i\in K^\times$ be such that $d_ia_{ij}\in\ogoth$ for all $j$ and
$d_ia_{ij}\in\ogoth^\times$ for some $j$.  Reducing
$d_is_i=\sum_jd_ia_{ij}\omega_j$ modulo $\goth P$, the right-hand side is
$\neq0$, so $d_is_i\in{\goth O}^\times$, and $w(s_i)=-v(d_i)$ is in
$v(K^\times)$, as claimed.

We have seen that at least one term in the sum $\sum_{i\in[0,e']}s_i\pi_i=0$
is $\neq0$.  Therefore there must be two different terms with the same
valuation ($<+\infty$), for $w(x+y)=w(x)$ if $w(x)<w(y)$.  If $i,j\in[0,e']$
are such that $i\neq j$, $s_i\neq0$, $s_j\neq0$ and $w(s_i\pi_i)=w(s_j\pi_j)$,
then, writing
$$
w(\pi_i)-w(\pi_j)=w(s_j)-w(s_i),
$$
we conclude that $w(\pi_i)$ and $w(\pi_j)$ have the same image in
$w(L^\times)/v(K^\times)$ (because $w(s_j),w(s_i)\in v(K^\times)$, as we have
seen), which is a contradiction.  Therefore $e$ is finite and 
$ef\le[L:K]$.

\th THEOREM 8
\enonce
Suppose that\/ $v$ is discrete (and\/ $K$ complete).  Then the\/
$\ogoth$-module\/ $\goth O$ is free of rank\/ $[L:K]$.
\endth
As $v$ is discrete, so is $w$, for $v(K^\times)$ is an index-$e$ subgroup of
$w(L^\times)$ and hence $v(K^\times)=ew(L^\times)$.  Let $\Pi$ be a
uniformiser of $L$~; we will show that the $ef$ elements $\Pi^i\omega_j$
($i\in[0,e[,j\in[1,f]$) in fact constitute an $\ogoth$-basis of $\goth O$.  We
have just seen (Proposition~7) that they are $K$-linearly independent, and
hence $\ogoth$-linearly independent.  That they generate the $\ogoth$-module
$\goth O$ can be shown along the same lines as the proof that every element of
$\goth O$ has a $\Pi$-adic expansion (Theorem~6, Lecture~5).  We give a more
conceptual version of the procedure.

Let $N$ (resp.~$M$) be the sub-$\ogoth$-module of $\goth O$ generated by the
$\omega_j$ ($j\in[1,f]$) (resp. by the~$\Pi^i\omega_j$ for $i\in[0,e[$ and
$j\in[1,f]$), so that
$$
M=N+\Pi N+\Pi^2N+\cdots+\Pi^{e-1}N.
$$
By the very definition (of $N$, $\Pi$ and the $\omega_j$), we have ${\goth
  O}=N+\Pi{\goth O}$.  So
$$
{\goth O}=N+\Pi{\goth O}=N+\Pi(N+\Pi{\goth O})=N+\Pi N+\Pi^2{\goth O}
=\ldots=N+\Pi N+\cdots+\Pi^e{\goth O}.
$$
Therefore ${\goth O}=M+\Pi^e{\goth O}$.  As $v(K^\times)$ is the index-$e$
subgroup of $w(L^\times)$, we have $\Pi^e{\goth O}=\pgoth{\goth O}$.  Hence
$$
{\goth O}=M+\pgoth{\goth O}=M+\pgoth(M+\pgoth{\goth O})=M+\pgoth^2{\goth O}
=\ldots=M+\pgoth^r{\goth O}
$$
for every $r>0$.  As the $(\pgoth^r{\goth O})_{r}$ form  a fundamental
system of neighbourhoods of the origin, the above equalities shows that $M$ is
{\it dense\/} in $\goth O$. 

On the other hand, $\ogoth$ is closed in $K$, and the topology of $L$ is that
of the product of $[L:K]$ copies of $K$ (Lemma~4, Lecture~11), so $M$ is {\it
  closed\/} in the sub-$K$-space $KM\subset L$, and hence in $L$, and hence in
$\goth O$.  As $M\subset\goth O$ is both dense and closed, $M=\goth O$.  

\th COROLLARY 9
\enonce
Suppose that\/ $v$ is discrete (and\/ $K$ complete).  Then\/ $ef=[L:K]$.
\endth
We have seen that the $\Pi^i\omega_j$ are an $\ogoth$-basis of $\goth O$, so
they are a $K$-basis of $L$, and hence $ef=[L:K]$.

\medskip {\it Example}\pointir The relation $ef=[L:K]$ does not always hold
when $v$ is not discrete.  For example, let $K$ be the completion of
$\Q_2(\sqrt2,\root4\of2,\root8\of2,\ldots)$ and let $L=K(\sqrt{-1})$.  Here
$e=1$, $f=1$, but $[L:K]=2$ (E.~Artin).

\vfill\eject

\centerline{\bf Lecture 13}
\bigskip
\centerline{\it Tame ramification}
\bigskip

Let $K$ be a field complete for valuation $v$.  Recall that a finite extension
$L|K$ of degree~$n$ and residual degree~$f$ is {\it unramified\/} if $f=n$ and
the residual extension is {\it separable}.  

{\it If\/ $L|K$ and\/ $M|L$ are finite unramified extensions, then so is
  $M|K$.}  Denoting by $m,l,k$ the residue fields, we have
$$
[M:K]=[M:L][L:K]=[m:l][l:k]=[m:k].
$$
Also, $m|k$ is separable because $m|l$ and $l|k$ are separable.  Hence $M|K$
is unramified.

Similarly, {\it every subextension of an unramified extension is unramified.}

An algebraic extension (possibly of infinite degree) is called {\it
  unramified\/} if every finite subextension is unramified.

\th PROPOSITION 1
\enonce
Let\/ $L$, $K'$ be two extensions of\/ $K$ in\/ $\bar K$.  If\/ $L|K$ is
unramified, then so is\/ $LK'|K'$.
\endth
It is sufficient to treat the case $L|K$ finite.  As we have seen in the
course of the proof of Proposition~6, Lecture~12, we have $L=K(\alpha)$, where
$\alpha$ is a (simple) root of a unitary polynomial~$f$ lifting the minimal
polynomial~$\bar f$ of a primitive element $a\in l$.

Now $LK'=K'(\alpha)$~; let $g\in K'[T]$ be the minimal polynomial of $\alpha$
over $K'$~; we have $g\in\ogoth'[T]$ (Lecture~9, Corollary~11), where
$\ogoth'$ is the ring of integers of $K'$.  Notice that $\bar g|\bar f$ in
$k'[T]$, so $k'(a)$ is separable over $k'$.  Also, $\bar g$ is irreducible by
Hensel's lemma (Lecture~20, Theorem~1), bacause $g$ is irreducible.  Finally,
$l'$ being the residue field of $LK'=K'(\alpha)$ (all we know for now is that
$k'(a)\subset l'$\/),% we have
$$
[l':k']\le[K'(\alpha):K']=\deg g=\deg\bar g=[k'(a):k']\le [l':k'],
$$
so $l'=k'(a)$ is separable over $k'$ and $[K'(\alpha):K']=[l':k']$, showing
that the extension $LK'|K'$ is unramified.

\th PROPOSITION 2
\enonce
The residue field\/ $\bar k$ of\/ $\bar K$ is an algebraic closure of\/ $k$.
\endth
It is enough to show that every unitary polynomial $\varphi\in k[T]$ of degree
$n>0$ has a root in $\bar k$.  Let $f\in\ogoth[T]$ be a unitary polynomial
lifting $\varphi$, and let $\alpha\in\bar K$ be a root of $f$.  As $\alpha$ is
integral over $\ogoth$, we have $\bar v(\alpha)\ge0$ (Corollary~5,
Lecture~11).  Clearly $\varphi(\bar\alpha)=0$.

\th DEFINITION 3
\enonce
Let\/ $L$ be an algebraic extension of\/ $K$.  The compositum\/ $L_0$ of all
unramified extensions of\/ $K$ in\/ $L$ is called the maximal unramified
extension of\/ $K$ in\/ $L$.   
\endth

\th PROPOSITION 4
\enonce
The residue field of\/ $L_0$ is the separable closure (\/{\rm fermeture})\/
of\/ $k$ in\/ $l$, and\/ $w(L_0^\times)=v(K^\times)$.
\endth

The residue field $l_0$ of $L_0$, being separable over $k$, is contained in
the separable closure $\tilde k$ of $k$ in $l$.  Conversely, for every
$a\in\tilde k$, there is an unramified extension of $K$ in $L$ with residue
field $k(a)$ (Lecture~12, Proposition~6), so $a\in l_0$, and hence $l_0=\tilde
k$.  To show that $w(L_0^\times)=v(K^\times)$, we may assume that $L_0$ is
finite over $K$ ; let $e_0={\bf(}w(L_0^\times):v(K^\times){\bf)}$.  That
$e_0=1$ now follows from Theorem~7, Lecture~12, which says that
$$
[L_0:K]\ge e_0[l_0:k]
=e_0[L_0:K].
$$
\vskip-\baselineskip
\th DEFINITION 5
\enonce
A finite extension\/ $L|K$ is said to be tamely ramified if the residual
extension\/ $l|k$ is separable and if\/ $[L:L_0]$ is prime to\/~$p$, the
characteristic exponent of\/ $k$.  An algebraic extension is said to be tamely
ramified if every finite subextension is tamely ramified.
\endth
Note that the quadratic extension $K(\sqrt{-1})$ of the completion $K$ of
$\Q_2(\sqrt2,\root4\of2,\ldots)$ is not tame, although the ramification index
$e=1$ is prime to~$2$.

It is clear that every subextension of $L$ is then tamely ramified.  Moreover,
if $M|L$ is tamely ramified, then so is $M|K$.

% Let $L|K$ be a finite tamely ramified extension, and let $L_0$ be the maximal
% unramified subextension~; the extension $[L:L_0]$ is totally ramified $(f=1)$
% of degree prime to $p$.

\th PROPOSITION 6
\enonce
Let\/ $K$ be complete for a valuation.  For\/ $a\in K^\times$ and\/ $m>0$
prime to the characteristic exponent\/~$p$ of the residue field, the
extension\/ $K(\!\root m\of  a)|K$ is tamely ramified. 
\endth
Consider the binomial $\varphi=T^m-a$, which may be reducible in $K[T]$.  Let
$\alpha$ be a root of $\varphi$, put $L=K(\alpha)$, and let $w$ be the unique
extension of $v$ to $L$.  Let $d$ be the order of the image of~$\alpha$ in the
group $w(L^\times)/v(K^\times)$.  As $\alpha^m\in K^\times$, the said image is
killed by~$m$, so $m$ is a multiple of $d$, say $m=dn$~; clearly, $n$ is prime
to~$p$. 

As $dw(\alpha)\in v(K^\times)$, there is a $b\in K^\times$ be such that
$v(b)+dw(\alpha)=0$~; put $\beta=b\alpha^d$, so that $w(\beta)=0$.  Also
$\beta^{n}=b^{n}a$, so $\beta$ is a root of $\psi=T^{n}-c$, with $c=b^na$ a
unit ($c\in\ogoth^\times$).  As $\psi'(\beta)=n\beta^{n-1}$, and $w(n)=0$, we
have $w(\psi'(\beta))=0$, and hence $\bar\psi\in k[T]$ is separable.
Therefore the extension $L_0=K(\beta)$, obtained as it is by adjoining a root
of a polynomial which is separable modulo~$\pgoth$, is unramified over~$K$.
In particular, we have $w(L_0^\times)=v(K^\times)$ (Proposition~4).

Now consider the binomial $bT^d-\beta$ over $L_0$~; it has $\alpha$ as a root,
so $[L:L_0]\le d$.  Because $w(L_0^\times)=v(K^\times)$, the group
$w(L^\times)/w(L_0^\times)$ is the same as $w(L^\times)/v(K^\times)$, which
has order $\ge d$, for it has an element of order~$d$~; the ramification index
of $L|L_0$ is therefore $\ge d$.  It follows that $L|L_0$ is totally ramified
of degree~$d$ (Lecture~12, Proposition~7)~; in particular, $bT^d-\beta$ is
irreducible in $L_0[T]$.

In summary, the extension $K(\!\root m\of a)|K$ is tamely ramified for any
$a\in K^\times$ and any $m$ prime to the characteristic exponent $p$ of the
residual field~$k$.  We shall see that iterating this process is essentially
the only way to get tamely ramified extensions.

\th LEMMA 7
\enonce
Suppose that\/ $L|K$ is tamely ramified.  If\/ $e=1$ and\/ $f=1$, then\/
$L=K$.   
\endth
Because $f=1$, we have $L_0=K$, and because $L|K$ is tamely ramified,
$n=[L:K]$ is prime to~$p$.  In particular, $L|K$ is separable.  Choosing an
algebraic closure $\bar L$ of $L$, with residue field $\bar l$, we have a
reduction map $\Hom_K(L,\bar L)\to\Hom_k(l,\bar l)$ sending $\sigma$ to
$\hat\sigma$. As $l=k$, $\hat\sigma$ is the inclusion of $l$ in $\bar l$, for
every $\sigma$.

The trace map $S_{L|K}:L\to K$ is $K$-linear, and surjective because $n$ is
invertible in $K$.  Let, if possible, $\alpha\in L^\times$ have trace
$S_{L|K}(\alpha)=0$.  As $e=1$, there is an element $a\in K^\times$ such that
$v(a)+w(\alpha)=0$.  Replacing $\alpha$ by $a\alpha$, we may assume that
$\alpha\in{\goth O}^\times$.  We thefore have
$S_{L|K}(\alpha)=\sum_{\sigma}\sigma(\alpha)=0$.  Reading this in $\bar l$, we
get $\sum_\sigma\hat\sigma(\hat\alpha)=0$.  But, as we have just remarked,
$f=1$ implies $\hat\sigma(\hat\alpha)=\hat\alpha$ for every~$\sigma$, leading
to $\hat n\hat\alpha=0$.  But $\hat n\hat\alpha=0$ is impossible because $\hat
n\neq0$ and $\hat\alpha\neq0$ in $k$.  Therefore $S_{L|K}$ is injective, and
hence an isomorphism.

\th PROPOSITION 8
\enonce
Let\/ $L$ be a totally but tamely ramified finite extension of a complete
unarchimedean field\/ $K$.  Then\/ $L=K(\!\root m_1\of{a_1},\ldots,\root
m_e\of{a_e}\,)$ for some $a_i\in K^\times$ and some \/ $m_i$ dividing\/ $e$. 
\endth
Let $\gamma_i\in L^\times$ be a system of representatives of the order-$e$
group $w(L^\times)/v(K^\times)$, and $m_i$ the order of the image of
$\gamma_i$ in the said group.  As $n=[L:K]$ kills every element of
$w(L^\times)/v(K^\times)$ (for $nw(L^\times)\subset v(K^\times)$), we have
$m_i|n$, and so $\gcd(m_i,p)=1$.

Let $c_i\in K^\times$ be such that $w(\gamma_i^{m_i})=v(c_i)$, and let
$\varepsilon_i\in{\goth O}^\times$ be such that
$\gamma_i^{m_i}=\varepsilon_ic_i$.  Because the residual extension $l|k$ is
trivial, we may write $\varepsilon_i=b_iu_i$ with $b_i\in\ogoth^\times$ and
$u_i\in\Ker({\goth O}^\times\to l^\times)$.  Now, modulo~$\goth P$, the
binomial $T^{m_i}-u_i$ is separable and has the root~$1$, so it has a root
$\beta_i\in\Ker({\goth O}^\times\to l^\times)$ (Hensel's Lemma)~; we have
$\beta_i^{m_i}=u_i$.

Put $\alpha_i=\gamma_i/\beta_i$~; notice that $w(\alpha_i)=w(\gamma_i)$, so
the ramification index of $L$ over $K(\alpha_1,\ldots,\alpha_e)$ is~$1$.  We
also have $\alpha_i^{m_i}=b_ic_i$, so we may say that $\alpha_i=\root
m_i\of{a_i}$, with with $a_i=b_ic_i$ in $K^\times$.

Now the extension $L|K(\!\root m_1\of{a_1},\ldots,\root m_e\of{a_e}\,)$ is
trivial because it is tamely ramified with ramifcation index~$1$ and residual
degree~$1$ (Lemma~7).  This completes the proof.

Remark that we could cut down on the number of radicals $\root m_i\of{a_i}$
that need to be adjoined to $K$ to get $L$ by fixing a system~$S$ of
generators of the group $w(L^\times)/v(K^\times)$ and limiting the $\gamma_i$
to representatives of~$S$.  For example, when $w(L^\times)/v(K^\times)$ is
cyclic, a single $\gamma$ would do.

Such is always the case when $v$ is discrete~; we can then take $\gamma$ to be
a uniformiser of $L$. Then $\alpha$ is also a uniformiser of $L$, for
$w(\alpha)=w(\gamma)$, and $a$ is a uniformiser of $K$, for
$a=N_{L|K}(\alpha)$, and hence $v(a)=ew(\alpha)$.  We have shown that every
totally but tamely ramified extension $L$ of a complete discretely valued
field $K$ is of the form $L=K(\root e\of\pi)$ for some uniformier $\pi$ of $K$
and some $e$ prime to~$p$.  Conversely, for every uniformier $\pi$ of $K$ and
every $e$ prime to $p$, the extension $K(\root e\of\pi)$ is totally tamely
ramified of ramification index~$e$ over $K$.  Note that, as the binomial
$T^e-\pi$ is Eisenstein, it is irreducible (Lecture~8, Lemma~7), and hence
$[K(\root e\of\pi):K]=e=ef$.  We shall shortly see that $[L:K]=ef$ holds for
any tamely ramified extension $L|K$ (Proposition~10).

\th COROLLARY 9
\enonce
A finite extension\/ $L$ of a complete unarchimedean field\/ $K$ is tamely
ramified precisely when\/ $L=L_0(\!\root m_1\of{a_1},\ldots,\root
m_r\of{a_r}\,)$, where\/ $L_0$ is the maximal unramified subextension, $a_i\in
L_0^\times$, and\/ $\gcd(m_i,p)=1$.
\endth
That every such extension is tamely ramified follows from Proposition~6 (the
case $r=1$) by induction on $r$ and the fact that tameness is transitive in a
tower of extensions.
\th PROPOSITION 10
\enonce
For every tamely ramified finite extension $L$ of a complete unarchimedean
field\/ $K$, we have\/ $ef=[L:K]$.
\endth
Let us say that an extension $L|K$ of unarchimedean fields, of degree~$n$,
ramification index~$e$, and residual degree~$f$, is {\it egalitarian\/} if
$ef=n$.  If $M|L$ is also egalitarian, then so is $M|K$.  Moreover, an
unramified extension is always egalitarian.

This remark and the characterisation of tamely ramified extensions
(Corollary~9), taken together, imply that it is sufficient to show that every
tame extension $L|K$ with $f=1$ and $w(L^\times)/v(K^\times)$ cyclic of
order~$e$ prime to~$p$ is egalitarian.  We may write, as we have seen,
$L=K(\root e\of a)$ for some $a\in K^\times$.  Therefore $n\le e$.  But we
know that $e\le n$ (Lecture~12, Proposition~7), hence $e=n$, and the proof is
complete.

{\eightpoint Some authors say that $L|K$ is tame if $ef=n$, if $l|k$ is
  separable, and if $\gcd(e,p)=1$.  We have seen that this is equivalent to
  Definition~5.}

\th PROPOSITION 11
\enonce
Let\/ $L$, $M$ be two extensions of\/ $K$ in\/ $\bar K$.  If\/ $L|K$ is
tamely ramified, then so is\/ $LM|M$.
\endth
As we already know this for $L|K$ unramified (Proposition~1), we may suppose
that $L|K$ is totally ramified.  We may further suppose that $L|K$ is finite,
and indeed of the form $L=K(\root m\of a)$ for some $m$ prime to~$p$
(Proposition~8).  Then $LM=M(\root m\of a)$ is tamely ramified over $M$ by
Proposition~6, completing the proof.

\th DEFINITION 12
\enonce
Let\/ $L$ be an algebraic extension of\/ $K$.  The compositum\/ $L'$ of all
tamely ramified extensions of\/ $K$ in\/ $L$ is called the maximal tamely
ramified extension of\/ $K$ in\/ $L$.   
\endth

Let $w(L^\times)'$ be the kernel of the map
$w(L^\times)\to\Gamma\to\Gamma\otimes_{\Z}\Z_p$, where
$\Gamma=w(L^\times)/v(K^\times)$~; it is the subgroup of those $x\in
w(L^\times)$ such that $mx\in v(K^\times)$ for some $m$ prime to~$p$.  It is
also the smallest subgroup of $w(L^\times)$ containing $v(K^\times)$ and of
index a power of~$p$.

\th PROPOSITION 13
\enonce
The residue field of\/ $L'$ is the separable closure 
of\/ $k$ in\/ $l$, and\/ $w(L'^\times)=w(L^\times)'$.
\endth
We need only consider the case $L|K$ finite.  Let $L_0$ be the maximal
unramified extension of $K$ in $L$, its residue field $l_0$ is the separable
closure of $k$ in $l$ (Proposition~4).  Because $L'|L_0$ is tamely ramified,
the residue field $l'$ of $L'$ contains $l_0$ and is separable over $l_0$,
hence $l'=l_0$.

The inclusion $w(L'^\times)\subset w(L^\times)'$ follows from the fact that
the subgroup $w(L'^\times)/v(K^\times)\subset w(L^\times)/v(K^\times)$ is of
order prime to~$p$.  Conversely, as we have seen while proving Proposition~8,
every element $\omega\in w(L^\times)'$ gives rise to an $\alpha\in L^\times$
such that $w(\alpha)=\omega$ and such that $\alpha^m\in K^\times$, where $m$
is the order of $\omega$ in $w(L^\times)/v(K^\times)$.  Proposition~6 shows
that $K(\root m\of a)$ is tamely ramified, so $\alpha\in L'$, and
consequently $\omega\in w(L'^\times)$.

Thus, in every finite extension $L$ of a complete unarchimedean field $K$,
there are two canonical subextensions, $L_0\subset L'\subset L$, the maximal
unramified extension and the maximal tamely ramified extension of $K$ in $L$.
They have the same residue field, namely the separable closure $l_0$ of $k$ in
$l$.  We have $\bar v(L_0^\times)=v(K^\times)$, whereas $\bar v(L'^\times)$ is
the inverse image in $\bar v(L^\times)$ of the maximal subgroup of $\bar
v(L^\times)/v(K^\times)$ of prime-to-$p$ order.  In particular, the extension
$L|L'$ is of $p$-power degree, and its residual extension is purely
inseparable.  

Passing to residue fields, the tower $K\to L_0\to L'\to L$ gives rise to the
tower $k\to l_0= l_0\to l$~; the degree of the total extension is $f=f_0p^r$,
where $f_0=[l_0:k]=[L_0:K]$ need not be prime to~$p$, and $l|l_0$ is purely
inseparable of degree~$p^r$ for some $r\in\N$.  The value groups are
$v(K^\times)=\bar v(L_0^\times)\subset\bar v(L^\times)'\subset\bar
v(L^\times)$, where the middle index $e'$ is prime to~$p$ and the last index
is a power $p^i$ of~$p$ for some $i\in\N$~; we have $e=e'p^i$.  Denoting by
$p^n$ the degree of $L|L'$, we have $r+i\le n$ (Proposition~7, Lecture~12).

{\it Example\/}~14\pointir\ Let $k$ be a field and consider $k(T)$ with the
valuation $v_T$ of Theorem~4, Lecture~2~; the residue field is $k$ and the
completion is $K=k\series{T}$ (Lecture~5, Corollary~7).  The maximal
unramified extension of $K$ (in an algebraic closure $\bar K$ of $K$) is
$\tilde K=\tilde k\series{T}$, where $\tilde k$ is the separable closure of
$k$ in $\bar K$.

Every tamely ramified extension of $\tilde K$ can be written $L=\tilde K(\root
e\of\pi)$, for some uniformiser $\pi$ of $\tilde K$ and some $e>0$ prime to
the characteristic exponent $p$ of $k$.  Every uniformiser $\pi$ is of the
form $\pi=cuT$ where $c\in{\tilde k}^\times$ and $u\in\Ker(\tilde
k[[T]]^\times\to\tilde k^\times)$ (explicitly, $u(0)=1$).  As $\root e\of c\in
{\tilde k}^\times$ and $\root e\of u\in{\tilde K}^\times$ exist, we may write
$L=\tilde K(\root e\of T)$.  If moreover $k$ is of characteristic~$0$, these
are the only finite extensions of $\tilde K$, for every finite extension is
tame.  In this case, $\bar K=\tilde K(\root2\of T,\root3\of T,\root4\of
T,\ldots)$.

\th PROPOSITION 15 (Abhyankar)
\enonce
Let\/ $K$ be a complete discretely valued field, $L|K$ a tamely ramified
extension of ramification index\/ $e$, and\/ $M|K$ a finite extension of
ramification index $de$ for some $d>0$.  Then the extension\/ $LM|M$ is
unramified.
\endth
First notice that the residue field of $LM$ is $lm$, the compositum of the
residue fields $l,m$ of $L,M$.  As $l|k$ is separable, so is $lm|m$. 

Recall that, the valuation being discrete, $L=L_0(\root e\of\lambda)$
for some uniformiser $\lambda$ of $L_0$~; we have $LM=L_0M(\root
e\of\lambda)$. 

Next, note that the extension $L_0M|M$ is unramified (Proposition~1).  By the
multiplicativity of the ramification index in the towers $L_0M|M|K$,
$L_0M|L_0|K$, the ramification index of $L_0M|L_0$ equals $de$.

Finally, let $\pi$ be a uniformiser of $L_0M$ and write $\lambda=u\pi^{de}$
for some unit $u$ of $L_0M$ (Lecture~5, Corollary~4).  As $e$ is prime to~$p$,
the binomial $T^e-u$ is separable over the residue field, and so the extension
$L_0M(\root e\of u)$ is unramified over $L_0M$.  But $(\pi^d\root e\of
u)^e=\lambda$, which we take to mean that $\root e\of\lambda\in L_0M(\root
e\of u)$.  Hence $LM$ is unramified over~$L_0M$, and, by transitivity,
over~$M$.

% \vfill\eject
% 
% Let $L$ be a totally tamely ramified extension of a complete unarchimedean
% field $K$.  Denote the valuations by $v$, $w$.  Show that there is a natural
% bijection between subgroups of $w(L^\times)/v(K^\times)$ and subextensions of
% $L|K$.  Does it preserve inclusions~?

\vfill\eject

\centerline{\bf Lecture 14}
\bigskip
\centerline{\it The inertia and the ramification groups}

\bigskip

Let $K$ be a complete unarchimedean field, $L$ a finite extension of $K$, and
$\bar L$ an algebraic closure of $L$~; denote the residue fields by $k,l,\bar
l$~; we know that $\bar l$ is algebraically closed (Lecture~13,
Proposition~2).

For any subextension $E\subset L$, denote by $\Hom_E(L,\bar L)$ the pointed
set of all $E$-morphisms of $L$ in $\bar L$, the ``\thinspace base
point\thinspace'' being the inclusion of $L$ in $\bar L$.  The canonical tower
$$
K\subset L_0\subset L'\subset L
$$
in which $L_0$ is the maximal unramified subextension, and $L'$ the maximal
tamely ramified subextension, gives rise to the inclusions of pointed sets
$$
\Hom_L(L,\bar L)\subset\Hom_{L'}(L,\bar L)\subset
\Hom_{L_0}(L,\bar L)\subset\Hom_K(L,\bar L)
$$
in which the smallest subset is reduced to its base point.

We also have the pointed set $\Hom_k(l,\bar l)$ of $k$-morphisms of $l$ in
$\bar l$, and a map of pointed sets $\Hom_K(L,\bar L)\to \Hom_k(l,\bar l)$,
denoted $\sigma\mapsto\hat\sigma$ (the ``\thinspace reduction\thinspace'' of
$\sigma$).

\th PROPOSITION 1
\enonce
The reduction map\/ $\Hom_K(L_0,\bar L)\to\Hom_k(l_0,\bar l)$ is a bijection.
\endth
(Notice that the two sets have the same finite number $[L_0:K]=[l_0:k]$ of
elements).  Let $a\in l_0$ be primitive over $k$, so that $l_0=k(a)$, and
$\hat\varphi\in k[T]$ its minimal polynomial.  Let $\varphi\in{\goth O}_0[T]$
be a unitary polynomial lifting $\hat\varphi$, and $\alpha\in{\goth O}_0$ the
root of $\varphi$ lifting $a$~; we have $L_0=K(\alpha)$ (Proposition~6,
Lecture~12).

Let $R(\varphi)\subset\bar L$ be the set of roots of $\varphi$~; they are
integral over ${\goth o}$ because $\varphi$ is unitary.  Similarly, let
$R(\hat\varphi)\subset\bar l$ be the set of roots of $\hat\varphi$~; the
reduction map induces a bijection $R(\varphi)\to R(\hat\varphi)$.  The
polynomials $\varphi$ and $\hat\varphi$ being separable, we also have the
bijection $\Hom_K(L_0,\bar L)\to R(\varphi)$ (resp.~$\Hom_k(l_0,\bar l)\to
R(\hat\varphi)$) which sends $f$ to $f(\alpha)$ (resp.~$f(a)$).  The theorem
now follows from these observations and the commutativity of the diagram
$$
\diagram{
\Hom_K(L_0,\bar L)&\droite{}&R(\varphi)\phantom{.}\cr
\vfl{}{}{6mm}&&\vfl{}{}{6mm}\cr
\Hom_k(l_0,\bar l)&\droite{}&R(\hat\varphi).\cr
\abovedisplayskip=5.0pt plus 3.0pt minus 3.0pt
\belowdisplayskip=5.0pt plus 3.0pt minus 3.0pt
}$$
\th PROPOSITION 2
\enonce
$\Hom_{L_0}(L,\bar L)=\Ker(\Hom_K(L,\bar L)\to\Hom_k(l,\bar l))$.
\endth
The kernel of a map $\sigma:(X,a)\to(Y,b)$ of pointed sets means, of course,
the set of those $x\in X$ such that $\sigma(x)=b$.

Let $l_0$ be the separable closure of $l$ in $\bar l$.  As $l_0$ is the
residue field of $L_0$ (Lecture~12, Proposition~6), it is clear that if
$\sigma\in\Hom_{L_0}(L,\bar L)$, then $\hat\sigma|_{l_0}$ is the inclusion
(of $l_0$ in $\bar l$).  Let us show that $\hat\sigma(a)=a$ even when $a\in
l$ is inseparable over $k$.

Let $\alpha\in\goth O$ be such that $\hat\alpha=a$.  As $a^{p^s}\in l_0$ for
some $s\in\N$, we have $\alpha^{p^s}=\beta+\gamma$, where $\beta\in{\goth
  O}_0$ is such that $\hat\beta=a^{p^s}$ and $w(\gamma)>0$.  As $\sigma$ fixes
$\beta$, we have $\hat\sigma(a^{p^s})=a^{p^s}$.  But we also have
$\hat\sigma(a^{p^s})=\hat\sigma(a)^{p^s}$.  Because $p$-th roots are unique in
a field of characteristic~$p$, we must have $\hat\sigma(a)=a$.

Now for the converse.  We know that $\Hom_K(L_0,\bar L)\to\Hom_k(l_0,\bar l)$
is a bijection (Proposition~1).  This being so, let $\sigma\in\Hom_K(L,\bar
L)$ be such that $\hat\sigma$ is the inclusion (of $l$ in $\bar l$)~; we have
to show that $\sigma\in\Hom_{L_0}(L,\bar L)$, or, equivalently, that
$\sigma|_{L_0}$ is the inclusion.  But this follows from the foregoing because
$\hat\sigma|_{l_0}$ is the inclusion, completing the proof that $\sigma$ is
the inclusion on $L_0$ if and only if $\hat\sigma$ is the inclusion (of $l$).

Let us next characterise $\Hom_{L'}(L,\bar L)$ as a subset of $\Hom_K(L,\bar
L)$.  Suppose first that $L'=L$.

\th PROPOSITION 3
\enonce
Suppose that\/ $L|K$ is tamely ramified.  If a $K$-morphism\/ $\sigma:L\to\bar
L$ is such that $\bar w(\sigma(x)-x)>w(x)$ for every\/ $x\in L^\times$, then
$\sigma$ is the inclusion.
\endth
Recall that $\bar w(\sigma(x))=w(x)$ for every $\sigma\in\Hom_K(L,\bar L)$ and
every $x\in L^\times$ (Lecture~11, Corollary~7), so $\bar w(\sigma(x)-x)\ge
w(x)$.  We have to show that if $\sigma$ is not the inclusion (of $L$ in $\bar
L$), then there is an $x\in L^\times$ such that $\bar w(\sigma(x)-x)=w(x)$.

Suppose first that $\hat\sigma$ is not the inclusion (of $l$ in $\bar l$), and
let $a\in l$ be such that $\hat\sigma(a)\neq a$.  Let $x\in{\goth O}^\times$
be a lift of $a$.  Clearly, $\sigma(x)-x$ is a unit, for its reduction is
$\neq0$.  We have found an $x\in L^\times$ such that $\bar
w(\sigma(x)-x)=w(x)$.

Suppose next that $\hat\sigma$ is the inclusion~; we know that
$\sigma\in\Hom_{L_0}(L,\bar L)$ (Proposition~1).  Also, $L=L_0(\!\root
m_1\of{a_1},\ldots,\root m_r\of{a_r}\,)$, where\/ $a_i\in L_0^\times$, and\/
$\gcd(m_i,p)=1$ (Corollary~9, Lecture~13).  As $\sigma$ is not the inclusion,
there is an $i\in[1,r]$ such that $\sigma(x)\neq x$ for $x=\root m_i\of{a_i}$.
We must have $\sigma(x)=\zeta x$ for some root $\zeta\neq1$ of $T^{m_i}-1$.
As $m_i$ is prime to~$p$, we have $\hat\zeta\neq1$, and
$$
\bar w(\sigma(x)-x)=w(x)+\bar w(\zeta-1)=w(x).
$$
We are now in a position to characterise the subset $\Hom_{L'}(L,\bar L)$ in
general, when $[L:L']$ is allowed to be $>1$.
\th PROPOSITION 4
\enonce
Let\/ $L$ be a finite extension of a complete unarchimedean field\/ $K$.  A\/
$K$-morphism\/ $\sigma\in\Hom_K(L,\bar L)$ is in\/ $\Hom_{L'}(L,\bar L)$
precisely when\/ $\bar w(\sigma(x)-x)>w(x)$ for every\/ $x\in L^\times$. 
\endth
If $\bar w(\sigma(x)-x)>w(x)$ for every\/ $x\in L^\times$, then it holds in
particular for every $x\in L^{\prime\times}$, and hence $\sigma|_{L'}$ is the
inclusion (Proposition~3).

Conversely, supposing that $\bar w(\sigma(x)-x)=w(x)$ for some\/ $x\in
L^\times$, we have to show that $\sigma|_{L'}$ is not the inclusion.  There is
nothing to prove if the characteristic~$p$ of $k$ is $0$, so assume that
$p\neq0$.  In the binomial expansion
$$
(\sigma(x)-x)^p=\sigma(x^p)+(-x)^p
  +\sum_{i\in[1,p[}{p\choose   i}\sigma(x)^i(-x)^{p-i},
$$
we have $w({p\choose i})>0$ (cf.~Lecture~8) and $\bar w(\sigma(x))=w(x)$,
so the valuation of ${p\choose i}\sigma(x)^i(-x)^{p-i}$ is $>w(x^p)$ for
every $i\in[1,p[$.  Combined with the hypothesis $\bar w(\sigma(x)-x)
=w(x)$, we get
$$
w(x^p)=\bar w((\sigma(x)-x)^p) =\bar w(\sigma(x)^p+(-x)^p)
$$ 
which implies $w(x^p)=\bar w(\sigma(x)^p-x^p)$ if $p\neq2$.  But even if
$p=2$, we can derive the same conclusion from 
$$
w(x^2)=\bar w(\sigma(x)^2+x^2)=\bar w(\sigma(x)^2-x^2+2x^2)
=\bar w(\sigma(x)^2-x^2),
$$
as $w(2x^2)>w(x^2)$.  Hence $w(x^p)=\bar w(\sigma(x)^p-x^p)$, whether $p$ is
even or uneven.  The process can be repeated to get $w(x^{p^i})=\bar
w(\sigma(x)^{p^i}-x^{p^i})$ for every $i\in\N$.

The image of $x^{p^i}$ in $w(L^\times)/v(K^\times)$ has order prime
to~$p$ if $i$ is large enough ($i\ge r$ would do, where $e=e'p^r$ is the
ramification index of $L|K$ and $e'$ is prime to $p$).  We may then write
$w(x^{p^i})=w(y)$ for some $y\in L'^{\times}$ (Proposition~13, Lecture~13), so
that $x^{p^i}/y=u$ is in ${\goth O}^\times$.

The image $\hat u\in l^\times$ of $u$ may or may not be in $l'^\times$, but
${\hat u}^{p^j}\in l'^\times$ for some $j\in\N$, for the extension $l|l'$ is
purely inseparable.  Thus, by taking $i$ to be large enough, we may assume
that $\hat u\in l'^\times$~; there is then a $z\in{\goth O}'^\times$ such that
$\hat z=\hat u$.  In other words, $w(x^{p^i}\!/y-z)>0$.  Putting things
together, we get
$$
w(x^{p^i}-yz)>w(y)=w(x^{p^i}),
$$
so that, posing $t=x^{p^i}-yz$, we have $w(t)>w(x^{p^i})$.  Applying $\sigma$,
we get $\sigma(t)=\sigma(x^{p^i})-\sigma(yz)$, with $\bar
w(\sigma(t))=w(t)>w(x^{p^i})$.  As a result, we have 
$$
\eqalign{
w(yz)=w(x^{p^i})=\bar w(\sigma(x^{p^i})-x^{p^i})
&=\bar w(\sigma(yz)-yz+\sigma(t)-t)\cr
&=\bar w(\sigma(yz)-yz),\cr
}$$
because $\bar w(\sigma(t)-t)>w(x^{p^i})$.  Having found an element $yz\in
L'^\times$ such that $w(yz)=\bar w(\sigma(yz)-yz)$, we may conclude that
$\sigma|_{L'}$ is not the inclusion (Proposition~3), completing the proof.

\medbreak
\centerline{\hbox to5cm{\hrulefill}}
\bigbreak

Now suppose that the extension $L|K$ is galoisian, and let $G=\Gal(L|K)$ be
the group of $K$-automorphisms of $L$.  To the subextensions $L_0\subset L'$
of $L|K$ correspond the subgroups $G'\subset G_0$ of $G$.  As we have seen
(Proposition~3), $G_0$ consists of those $K$-automorphisms of $L$ which induce
$\Id_l$ on the residue field~; it is called the {\it inertia subgroup\/} of
$G$. We have also seen that $G'$, the {\it ramification subgroup\/} of $G$,
consists of those $K$-automorphisms $\sigma$ of $L$ such that
$w(\sigma(x)-x)>w(x)$ for every $x\in L^\times$ (Proposition~4).

For every $\sigma\in G$, the subextension $M_0$ (resp.~$M'$) corresponding to
$\sigma G_0\sigma^{-1}$ (resp.~ $\sigma G'\sigma^{-1}$) is an unramified
(resp.~tamely ramified) extension of $K$ in $L$, hence $M_0=L_0$
(resp.~$M'=L'$).  In other words, the subgroups $G_0$ and $G'$ are invariant
under conjugation in $G$.  We wish to determine the quotient groups
$G/G_0=\Gal(L_0|K)$, $G/G'=\Gal(L'|K)$ and $G_0/G'=\Gal(L'|L_0)$.

Notice first that the residual extension $l_0|k$ is also galoisian.  Indeed,
it is separable because $L_0|K$ is unramified.  Let $a\in l_0$ be a primitive
element, $\hat g\in k[T]$ its minimal polynomial, $\alpha\in{\goth O}_0$ a
lift of $a$, and $f\in K[T]$ the minimal polynomial of $\alpha$~; we have
$f\in\ogoth[T]$ because $\alpha$ is integral over $\ogoth$ (Lecture~9,
Corollaries~10,~11).  As $L_0|K$ is normal, $f$ factors completely over
$L_0$, and therefore $\hat f$ factors completely over $l_0$.  As $a$ is a root
of $\hat f$, we have $\hat g|\hat f$, and hence $\hat g$ factors completely
over $l_0$, showing that $l_0|k$ is normal.

Proposition~1 implies that the reduction map $\Gal(L_0|K)\to\Gal(l_0|k)$ is an
isomorphism of groups.  Indeed, we have a commutative square 
$$
\diagram{
\Gal(L_0|K)&\droite{}&\Hom_K(L_0,\bar L)\cr
\vfl{}{}{6mm}&&\vfl{}{}{6mm}\cr
\Gal(l_0|k)&\droite{}&\Hom_k(l_0,\bar l)\cr
\abovedisplayskip=5.0pt plus 3.0pt minus 3.0pt
\belowdisplayskip=5.0pt plus 3.0pt minus 3.0pt
}$$
in which the horizontal arrows ($\sigma\mapsto i\circ\sigma$, $i$ being the
inclusion) are bijective, the vertical arrows are the reduction maps
($\sigma\mapsto\hat\sigma$), of which the one on the right is bijective
(Proposition~1), and hence so is the one on the left.  Thus our interpretation
of $G/G_0$ is that the map which sends $\sigma\in G$ to the reduction of 
$\sigma|_{L_0}$ induces an isomorphism $G/G_0\to\Gal(l_0|k)$.  We have proved
the following proposition.

\th PROPOSITION 5
\enonce
The reduction map\/ $\Gal(L|K)\to\Gal(l_0|k)$ is surjective, with kernel \/
$G_0=\Gal(L|L_0)$~; the sequence $1\to G_0\to G\to\Gal(l_0|k)\to1$ is exact.  
\endth

Next, we will define a bimultiplicative pairing $\Phi:w(L^\times)\times G_0\to
l^\times$.  For $x\in L^\times$ and $\sigma\in G_0$, define $\Phi(x,\sigma)$
to be the image in $l^\times$ of the unit $\sigma(x)/x\in{\goth O}^\times$.
For $u\in{\goth O}^\times$ and $\sigma\in G_0$, we have $\Phi(u,\sigma)=1$,
because $\hat\sigma$ is an $l_0$-automorphism of the purely inseparable
extension $l$, and hence $\hat\sigma=\Id_l$.  The identity
$$
{\sigma(ux)\over ux}={\sigma(u)\over u}{\sigma(x)\over x}.
$$
shows that $\Phi(x,\sigma)$ depends only on the image of $x$ in
$w(L^\times)=L^\times\!/{\goth O}^\times$~; we thus get a map
$\Phi:w(L^\times)\times G_0\to l^\times$.

It is easy to see that $\Phi$ is bimultiplicative.  Indeed, for $\bar x,\bar
y\in w(L^\times)$ and $\sigma,\tau\in G_0$, we have
$$
{\sigma(xy)\over xy}
={\sigma(x)\over x}{\sigma(y)\over y}~;\quad
{\sigma\tau(x)\over x}={\sigma(\tau(x))\over\tau(x)}{\tau(x)\over x}
={\sigma(u)\over u}{\sigma(x)\over x}{\tau(x)\over x}
$$
where $u=\tau(x)/x$ is in ${\goth O}^\times$.  These imply that
$\Phi(xy,\sigma)=\Phi(x,\sigma)\Phi(y,\sigma)$ and
$\Phi(x,\sigma\tau)=\Phi(x,\sigma)\Phi(x,\tau)$, because $\Phi(u,\sigma)=1$. 

Let us determine the left- and
the right-kernel of $\Phi$.  The right kernel consists of those $\sigma\in
G_0$ such that $\Phi(\bar x,\sigma)=1$ for every $\bar x\in w(L^\times)$.
This is the same as saying that $\sigma(x)/x\equiv1\pmod{\goth P}$, which is
equivalent to
$$
w\left({\sigma(x)\over x}-1\right)>0\ \ \Leftrightarrow\ \ 
w(\sigma(x)-x)>w(x)\ \ \Leftrightarrow\ \ 
\sigma\in G'
$$
by Proposition~4.  As for the left kernel, the set of $\bar x\in w(L^\times)$
such that $\Phi(\bar x,\sigma)=1$ for every $\sigma\in G_0$, it clearly
contains $v(K^\times)=w(L_0^\times)$, as $\sigma(x)=x$ for every $x\in
L_0^\times$. 
We have thus proved the following result.
\th PROPOSITION 6
\enonce
The ramification group\/ $G'$ is the kernel of the canonical homomorphism
$$
G_0\to\Hom(w(L^\times)/v(K^\times),\,l^\times),\quad 
\sigma\mapsto\left(\bar x\mapsto{\sigma(x)\over x}\pmod{{\goth P}}\right). 
$$
\endth
Notice that, as the only element of $p$-power order in $l^\times$ is~$1$, we
have 
$$
\Hom(w(L^\times)/v(K^\times),\,l^\times)=
\Hom(w(L'^\times)/v(K^\times),\,l^\times),
$$
by Lecture~13, Proposition~13.  Also, the order $e'$ of
$w(L'^\times)/v(K^\times)$ is prime to~$p$, hence so is that of
$\Hom(w(L'^\times)/v(K^\times),\,l^\times)$.  It follows that every
sub-$p$-group of $G_0$ is contained in the $p$-group $G'$.

\th PROPOSITION 7
\enonce 
The ramification subgroup\/ $G'\subset G_0$ is the only maximal
sub-$p$-group of\/ $G_0$.  
\endth 

The group $G_0/G'=\Gal(L'|L_0)$ has order $e'=(w(L'^\times):v(K^\times))$.  If
it has a cyclic quotient of order~$m>0$, then $L_0$ has a degree-$m$ extension
which, being cyclic and totally tamely ramified, is of the form $L_0(\!\root
m\of a)$ for some irreducible $T^m-a\in L_0[T]$.  It follows that
$L_0^\times$, and hence $l^\times$, has an element of order~$m$.  As a result,
$\Hom(w(L'^\times)/v(K^\times),\,l^\times)$ is the dual
$(w(L'^\times)/v(K^\times))^\vee$ of $w(L'^\times)/v(K^\times)$ and has order
$e'$.  This can be expressed as follows.

\th PROPOSITION 8
\enonce
The sequence\/ $1\to G'\to G_0\to
(w(L'^\times)/v(K^\times))^\vee\to1$ is exact.
\endth

Let us finally compute $\Phi(x,\sigma\tau\sigma^{-1})$ for $\sigma\in G$,
$\tau\in G_0$.  Writing $\sigma^{-1}(x)=ux$ ($u\in{\goth O}^\times$), we have 
$$
{\sigma\tau\sigma^{-1}(x)\over x}=
\sigma\left({\tau(\sigma^{-1}(x)\over\sigma^{-1}(x)}\right)=
\sigma\left({\tau(ux)\over ux}\right)=
\sigma\left({\tau(u)\over u}{\tau(x)\over x}\right)
$$
Reducing modulo~$\goth P$ gives $\Phi(x,\sigma\tau\sigma^{-1})
=\hat\sigma(\Phi(u,\tau)\Phi(x,\tau))=\hat\sigma(\Phi(x,\tau))$, for
$\Phi(u,\tau)=1$.  In particular, if $G$ is commutative, so that
$\sigma\tau\sigma^{-1}=\tau$, the element $\Phi(x,\tau)\in l^\times$ is
invariant under every $k$-automorphism $\hat\sigma$ of $l$.  At the same time,
as $\Phi(x,\tau)\in l^\times$ has (finite) order prime to~$p$, it is
separable over~$k$, and hence belongs to $k^\times$. 

% We also have the exact sequence $1\to G_0/G'\to G/G'\to G/G_0\to1$.  

\vfill\eject

\centerline{\bf Lecture 15}
\bigskip
\centerline{\it Higher ramification groups}

\bigskip

Let $K$ be a complete unarchimedean field and $L|K$ a finite galoisian
extension, of group $G=\Gal(L|K)$.  Let $w$ be the valuation of $L$, and
suppose that $w(L^\times)=\Z$ if $w$ is discrete.  For each real $r\in\R$, we
shall define a subgroup $G_{[r}\subset G$ ({\it vee-sub-ar-inclusively}).

Define $G_{[r}$ to be the set of those $\sigma\in G$ such that
$w(\sigma(x)-x)-w(x)\ge r$ for every $x\in L^\times$~; clearly $G_{[0}=G$.
For $r\le s$, we have $G_{[s}\subset G_{[r}$.

\th LEMMA 1
\enonce
Let\/ $S$ be a system of generators of the group\/ $L^\times$, and let\/
$\sigma\in G$.  If\/ $w(\sigma(x)-x)-w(x)\ge r$ for every\/ $x\in S$, then
$\sigma\in G_{[r}$.
\endth
We have $(\sigma-1)(xy)=\sigma(x)(\sigma-1)(y)+y(\sigma-1)(x)$.  Dividing
throughout by $xy$ and noting that $w(\sigma(x))=w(x)$, we get
$$
w\left({\sigma(xy)\over xy}-1\right)
\ge\inf\left(w\left({\sigma(x)\over x}-1\right),
 w\left({\sigma(y)\over y}-1\right)\right).
$$
Also, $\displaystyle w\left({\sigma(x^{-1})\over x^{-1}}-1\right)
=w\left({x\over\sigma(x)}-1\right)=w(x-\sigma(x))-w(\sigma(x))$.  From these
two facts it follows that if $w(\sigma(x)-x)-w(x)\ge r$ holds for every\/
$x\in S$, then it holds for every $x$ in the subgroup generated by $S$.

\th LEMMA 2
\enonce For every\/ $r$, $G_{[r}$ is a subgroup of\/ $G$~; it is invariant
under conjugation in\/ $G$.
\endth

If $\sigma,\tau\in G_{[r}$, then $w(\sigma(x)-x)-w(x)\ge r$ and
$w(\tau(x)-x)-w(x)\ge r$ for every $x\in L^\times$.  Notice that
$(\sigma\tau-1)(x)=(\sigma-1)(\tau(x))+(\tau-1)(x)$, and $w(\tau(x))=w(x)$,
hence $\sigma\tau\in G_{[r}$.  Also, every $y\in L^\times$ can be written
$y=\sigma(x)$ for some $x\in L^\times$, which implies that
$$
w(\sigma^{-1}(y)-y)=w(x-\sigma(x))\ge w(x)+r=w(y)+r
$$
and $\sigma^{-1}\in G_{[r}$.  Finally, because
$(\tau\sigma\tau^{-1}-1)(x)=\tau(\sigma(\tau^{-1}(x))-\tau^{-1}(x))$, these
subgroups are invariant under conjugation by $\tau\in G$.

\th LEMMA 3
\enonce
For sufficiently large\/ $r$, we have $G_{[r}=\{\Id_L\}$.
\endth
Write $L=K(z)$ for some element $z$ primitive over $K$, and, for $\sigma\in
G$, define $i_{L|K}(\sigma)=w(\sigma(z)-z)-w(z)$.  If
$i_{L|K}(\sigma)=+\infty$, then $\sigma=\Id_L$, and conversely.  Thus if $r$
is very large (say, strictly larger than $i_{L|K}(\sigma)$ for every
$\sigma\neq\Id_L$ in $G$), then the subgroup $G_{[r}$ is reduced to the
neutral element $\Id_L$.

\medskip

It follows that the filtration $(G_{[r})_r$ of $G$ has ``\thinspace
jumps\thinspace'' at some $r_1,r_2,\ldots,r_n\in\R$.  This means that
$G_{[r_i+\varepsilon}\neq G_{[r_i}$ for every $\varepsilon>0$, and the $r_i$
are the only real numbers with this property.  In other words, the real line
$\R$, where the indices of the filtration lie, gets partitioned into a certain
number of intervals $]-\infty,r_1]$, $]r_1,r_2]$, $\ldots$, $]r_n,+\infty[$
(with $r_1\ge0$) in each of which the group $G_{[r}$ is constant, but gets
strictly smaller as we move from one interval to the next, as indicated below
$$
\vbox{\halign{&\hfil$#$\hfil\quad\cr
r\in&]-\infty,r_1]&]r_1,r_2]&\cdots&]r_{n-1},r_n]&]r_n,+\infty[\cr
\noalign{\vskip-5pt}
\multispan6\hrulefill.\cr
G_r=&G&G_{[r_2}&\cdots&G_{[r_n}&\{\Id_L\}\cr
}}
$$

\th LEMMA 4
\enonce
If\/ $\sigma\in G_{[r}$ and\/ $\tau\in G_{[s}$, 
then\/ $\sigma\tau\sigma^{-1}\tau^{-1}\in G_{[r+s}$.  
\endth
For every $x\in L^\times$, we have
$$\eqalign{
w((\sigma-1)(\tau-1)(x))&\ge w((\tau-1)(x))+r\cr
 &\ge w(x)+s+r.
}$$
Similarly, $w((\tau-1)(\sigma-1)(x))\ge w(x)+s+r$, and, in view of the
identity $\sigma\tau-\tau\sigma=(\sigma-1)(\tau-1)-(\tau-1)(\sigma-1)$, we
conclude that $w((\sigma\tau-\tau\sigma)(x))\ge w(x)+s+r$.

Taking $x=\sigma^{-1}\tau^{-1}(y)$ in the above estimate, so that $w(x)=w(y)$,
we get $w((\sigma\tau\sigma^{-1}\tau^{-1}-1)(y))\ge w(y)+r+s$ for every $y\in
L^\times$.  Hence $\sigma\tau\sigma^{-1}\tau^{-1}\in G_{[r+s}$.

\th COROLLARY 5
\enonce
The derived subgroup of\/  $G_{[r}$ is contained in\/ $G_{[2r}$.
\endth
This is the case $r=s$ of Lemma~4.

\th COROLLARY 6
\enonce
If a jump occur at\/ $r_i>0$, then the quotient group\/
$G_{[r_i}/G_{[r_{i+1}}$ is commutative. 
\endth
Indeed, we have $G_{[2r_i}\subset G_{[r_i+\varepsilon}=G_{[r_{i+1}}\subset
G_{[r_i}$ for sufficiently small $\varepsilon>0$.  But $G_{[r_i}/G_{[2r_i}$ is
commutative (Lemma~4), therefore so is $G_{[r_i}/G_{[r_{i+1}}$.

Let us show that these quotients are killed by $p$, the characteristic
exponent of the residue field.

\th PROPOSITION 7
\enonce
If a jump occur at\/ $r_i>0$, then the group\/ $G_{[r_i}/G_{[r_{i+1}}$ is an 
elementary abelian $p$-group.  
\endth
Notice that if $\sigma\in G_{[r}$, then $\sigma^p\in G_{[\lambda(r)}$, where
$\lambda(r)=\inf(pr,w(p)+r)$.  Indeed, we may write
$$\eqalign{
\sigma^p-1&=((\sigma-1)+1)^p-1\cr
 &=(\sigma-1)^p+p(\sigma-1)^{p-1}+\cdots+p(\sigma-1)
}$$
so that
$(\sigma^p-1)(x)=(\sigma-1)^p(x)+p(\sigma-1)^{p-1}(x)+\cdots+p(\sigma-1)(x)$,
from which it follow that $w((\sigma^p-1)(x))\ge w(x)+\lambda(r)$, as was to be
shown. 

Next, as $w(p)>0$, we have $\lambda(r)>r$ for every $r\in\,]0,+\infty[$.  In
particular, we have the inclusions $G_{[\lambda(r_i)}\subset
G_{[r_i+\varepsilon}=G_{[r_{i+1}}\subset G_{[r_i}$ for sufficiently small
$\varepsilon>0$.  We have seen that the group $G_{[r_i}/G_{[\lambda(r_i)}$ is
killed by $p$, and hence so is the quotient $G_{[r_i}/G_{[r_{i+1}}$ (which is
commutative by Corollary~6).

\medbreak
\centerline{\hbox to5cm{\hrulefill}}
\bigbreak

Let us now assume, in addition to the field $K$ being complete for $v$, that
$v$ is discrete, and allow $L|K$ to be any finite extension, galoisian or not.

\th THEOREM 8
\enonce
If the residual extension $l|k$ is separable, then ${\goth O}=\ogoth[\omega]$
for some\/ $\omega\in{\goth O}$. 
\endth
Consider first the case when $L=L_0$ is unramified over $K$.  Let $a\in l$ be
an element primitive over $k$~; then $1,a,\ldots,a^{f-1}$ is a $k$-basis of
$l$.  If $\omega\in{\goth O}^\times$ is a lift of $a$, then $1,\omega,\ldots,
\omega^{f-1}$ is an $\ogoth$-basis of $\goth O$ (Lecture~12, Theorem~8). We
have ${\goth O}=\ogoth[\omega]$.

Next consider the totally ramified case $L_0=K$.  If $\Pi$ is a uniformiser of
$L$, then we know that $1,\Pi,\ldots,\Pi^{e-1}$ is an $\ogoth$-basis of $\goth
O$ (Lecture~12, Theorem~8). We may take $\omega=\Pi$.

In the general case, write $l=k(a)$, and let $\varphi\in\ogoth[T]$ be any lift
of the minimal polynomial of $a$ over $k$~; we have $\hat\varphi(a)=0$ but
$\hat\varphi'(a)\neq0$.  For any lift $\omega\in{\goth O}^\times$ of $a$, we
have $w(\varphi(\omega))>0$, hence $w(\varphi(\omega))\ge1$.  If
$w(\varphi(\omega))>1$, notice that for the lift $\omega'=\omega+\Pi$ of $a$,
we have $w(\varphi(\omega'))=1$. Indeed,
$$
\varphi(\omega')
=\varphi(\omega+\Pi)=\varphi(\omega)+\varphi'(\omega)\Pi+\alpha\Pi^2
$$
for some $\alpha\in{\goth O}$.  As $w(\varphi(\omega))>1$ and
$w(\alpha\Pi^2)>1$, but $w(\varphi'(\omega)\Pi)=1$ since
$\hat\varphi(a)\neq0$, we have $w(\varphi(\omega'))=1$.  Replacing $\omega$ by 
$\omega'$, we may assume that $w(\varphi(\omega))=1$, in addition to $\hat
\omega=a$.

Consider $\omega^{j-1}$ ($j\in[1,f]$)~; their reduction is the $k$-basis
$1,a,\ldots,a^{f-1}$ of $l$.  Also, $\varphi(\omega)$ is a uniformiser of $L$,
by cour choice of $\omega$.  We have seen that $\varphi(\omega)^i
\omega^{j-1}$ ($i\in[0,e[$, $j\in[1,f]$) is an $\ogoth$-basis of $\goth O$
(Lecture~12, Theorem~8).  But each of these elements is in $\ogoth[\omega]$,
hence ${\goth O}=\ogoth[\omega]$.

\th COROLLARY 9
\enonce
We have $\bar w(\sigma(\omega)-\omega)=\inf_{x\in{\goth O}}\bar
w(\sigma(x)-x)$ for every\/ $K$-morphism\/ $\sigma:L\rightarrow\bar L$. 
\endth
Write $x=a_0+a_1\omega+\cdots+a_{n-1}\omega^{n-1}$ ($a_i\in\ogoth$).  We have
$\sigma(a_i)=a_i$, so
$$
\eqalign{
\sigma(x)-x
&=a_1(\sigma(\omega)-\omega)+\cdots
  +a_{n-1}(\sigma(\omega)^{n-1}-\omega^{n-1})\cr 
&=(\sigma(\omega)-\omega)(a_1+a_2(\sigma(\omega)+\omega)+\cdots),\cr
}$$
where the second factor is in $\sigma(\goth O)$.  Hence $\bar
w(\sigma(x)-x)\ge\bar w(\sigma(\omega)-\omega)$, with equality for
$x=\omega$. 

\th PROPOSITION 10
\enonce
Let\/ $L|K$ be a totally ramified extension, and\/ $\Pi$ a uniformiser of\/
$L$.  Then ${\goth O}=\ogoth[\Pi]$, and the minimal polynomial of\/ $\Pi$ is
Eisenstein.  Conversely, for every Eisenstein polynomial\/
$\varphi\in\ogoth[T]$, the extension\/ $K[T]/\varphi K[T]$ is totally
ramified, and the image of\/ $T$ is a uniformiser.
\endth

Let $\varphi=T^e+a_{e-1}T^{e-1}+\cdots+a_0$ be the minimal polynomial of a
uniformiser $\Pi$ of $L$~; it has degree $e=[L:K]=(w(L^\times):v(K^\times))$
because ${\goth O}=\ogoth[\Pi]$ (Theorem~8).  As $N_{L|K}(\Pi)=(-1)^e a_0$, we
have $v(a_0)=ew(\Pi)$, and hence $a_0$ is a uniformiser of $K$.  Moreover,
$v(a_i)>0$ for $i\in[1,e[$, because the Newton polygon of $\varphi$ is the
line segment joining $(0,v(a_0))$ and $(n,0)$ (cf.\ Corollary~11, Lecture~9).
Hence $\varphi$ is Eisenstein.

Conversely, let $\varphi\in\ogoth[T]$ be Eisenstein of degree~$n$~; in
particular $\varphi(0)$ is a uniformiser of $K$.  We know that $\varphi$ is
irreducible (Lecture~8, Lemma~7), so $L=K(t)$ is of degree~$n$ over $K$, where
$t$ is a root of $\varphi$.  Extending $v$ to $w$ on $L$, we have
$nw(t)=v(\varphi(0))$, so the ramification index $e$ of $L|K$ is $\ge n$, and
hence $e=n$ (Corollary~9, Lecture~12).  Also, $t$ is a uniformiser of $L$
because $w(t)=v(\varphi(0))/e$.

\th COROLLARY 11
\enonce
Suppose that\/ $L|K$ is totally ramified, and let\/ $\Pi_1,\Pi_2$ be 
uniformisers of\/ $L$.  Then\/ $\bar w(\sigma(\Pi_1)-\Pi_1)=\bar
w(\sigma(\Pi_2)-\Pi_2)$  for every\/ $K$-morphism\/ $\sigma:L\to\bar L$. 
\endth
Indeed, we have ${\goth O}=\ogoth[\Pi_i]$ (Proposition~10) on the one hand,
and $w(\sigma(\Pi_i)-\Pi_i)=\inf_{x\in{\goth O}}w(\sigma(x)-x)$ (Corollary~9)
on the other.

\medbreak
\centerline{\hbox to5cm{\hrulefill}}
\bigbreak

Assume now that the (finite) extension $L|K$ is galoisian with separable
residual extension $l|k$, where $K$ is complete for the discrete valuation
$v$~; put $G=\Gal(L|K)$.  The valuation $w$ on $L$ is discrete~; normalising
it so that $w(L^\times)=\Z$, we have $G_{[r}=G_{[n}$, where $n=\lceil r\rceil$
is the smallest integer $\ge r$.  When the residual extension is separable, we
also have ${\goth O}=\ogoth[\omega]$ for some $\omega\in{\goth O}$
(Theorem~8).

\th LEMMA 12
\enonce
Suppose that\/ $l|k$ is separable.  Then\/ $\sigma\in G$ is in $G_0$ if and
only if\/ $w(\sigma(\omega)-\omega)>0$, where\/ $G_0$ is the inertia subgroup
of\/ $G$. 
\endth
Indeed, $w(\sigma(\omega)-\omega)>0$ is equivalent to
$\sigma(\omega)\equiv\omega\pmod{\goth P}$, which is equivalent to
$\hat\sigma(a)=a$ for every $a\in l$, because $l=k(\hat\omega)$, and thus to
$\sigma\in G_0$ (Lecture~14, Proposition~5).

\th LEMMA 13
\enonce
The ramification subgroup\/ $G'\subset G$ is the same as\/ $G_{[1}$.
\endth
Indeed, in view of Proposition~4, Lecture~14, we have
$$
\eqalign{
\sigma\in G'
& \Leftrightarrow\  w(\sigma(x)-x)-w(x)>0\ (\hbox{for all }x\in L^\times)\cr
& \Leftrightarrow\  w(\sigma(x)-x)-w(x)\ge1\cr
& \Leftrightarrow\ \sigma\in G_{[1}.\cr
}$$
It follows that $G_{[n}\subset G_0$ for every $n>0$.  Indeed, we have the
inclusions $G_{[n}\subset G_{[1}=G'\subset G_0$.

\medbreak
\centerline{\hbox to5cm{\hrulefill}}
\bigbreak

Classically, when $v$ is discrete as here, one defines the decreasing sequence
of subgroups $G_n\subset G$ ($n=-1,0,1,\ldots$) as
$$
G_n=\{\sigma\in G\;|\;w(\sigma(x)-x)\ge n+1\hbox{ for every }x\in{\goth O}\}.
$$
It is clear that $G_{-1}=G$, and that the new $G_0$ as defined here is the
inertia group $G_0$ as defined earlier as the kernel of the map
$\Gal(L|K)\to\Aut_k(l)$ (Lecture~14, Proposition~5).

\th PROPOSITION 14
\enonce
Suppose that the residual extension is separable (and hence galoisian).  We
have\/ $G_{[n}=G_n$ for every\/ $n>0$, and also for\/ $n=0$ if\/ $L|K$ is 
totally ramified (but $G_{[0}\neq G_0$ if\/ $G_0\neq G$).
\endth
We have $G_{[0}=G$ in general, and $G_0=G$ if and only if $L|K$ is totally
ramified, so the case $n=0$ is clear.  For $n>0$, both $G_{[n}$ and $G_n$ are
subgroups of $G_0$, so their elements are $L_0$-automorphisms (of $L$).  

Let $\sigma\in G_{[n}$.  Let $\Pi$ be a uniformiser of $L$, so that ${\goth
  O}={\goth O}_0[\Pi]$ (Theorem~8).  Then, in particular for $\Pi$, we have
$w(\sigma(\Pi)-\Pi)-w(\Pi)\ge n$, which implies, in view of $w(\Pi)=1$, that
$w(\sigma(x)-w(x)\ge n+1$ for every $x\in{\goth O}$ (Corollary~9).  Therefore
$\sigma\in G_n$.

Conversely, take a $\sigma\in G_n$~; we have, in particular,
$w(\sigma(\Pi)-\Pi)\ge1+n$ for every uniformiser $\Pi$ of $L$.  This means
that $w(\sigma(\Pi)-\Pi)-w(\Pi)\ge n$ for such $\Pi$.  But the set $w^{-1}(1)$
of uniformisers of $L$ generates the group $L^\times$, so we have $\sigma\in
G_{[n}$ (Lemma~1).

\medbreak
\centerline{\hbox to5cm{\hrulefill}}
\bigbreak

Let us return to the general situation of a field $K$ complete for a valuation
$v$ (discrete or not), and a finite galoisian extension $L|K$, of group
$G=\Gal(L|K)$.

Define $G_{]r}$ ({\it vee-sub-ar-exclusively}) to be the set of those
$\sigma\in G$ such that $w(\sigma(x)-x)-w(x)>r$ for every $x\in L^\times$.
Being the union of the increasing sequence of subgroups $G_{[r+\varepsilon}$
(Lemma~2) as $\varepsilon\to0+$, each $G_{]r}$ is a subgroup of $G$ invariant
under conjugation in $G$.  For $r<0$, the group $G_{]r}$ contains $G_{[0}$,
therefore $G_{]r}=G$, because $G_{[0}=G$.

Many of the properties of the subgroups $G_{[r}$ extend to the newly-defined
$G_{]r}$, by a sort of limiting process.  We skip the details.

\th PROPOSITION 15
\enonce
We have $G_{]0}=G'$, the ramification group.
\endth
This holds because a $\sigma\in G$ is in $G'$ if and only if
$w(\sigma(x)-x)>w(x)$ for every $x\in L^\times$ (Proposition~4, Lecture~14). 

\medbreak
\centerline{\hbox to5cm{\hrulefill}}
\bigbreak

In the special case $e=1$, $f=1$ (in which we have the equality
$w(L^\times)=v(K^\times)$ of the value groups and $l=k$ of residue fields),
let us show that $G_{]r}=G_{[r}$ for every $r\in\R$.  This may fail in
general. 

For every $x\in L^\times$, there exists a $b\in K^\times$ such that
$w(x)=v(b)$.  Further, there is a unit $c\in\ogoth^\times$ such that
$x/b\equiv c$ in $l=k$.  Therefore $w(x-bc)>v(b)=w(x)$.  Writing $a=bc$, we
have $w(x-a)>w(x)$, and, writing $y=x-a$, we have $w(y)>w(x)$.  Now,
$$
(\sigma-1)(x)=(\sigma-1)(a)+(\sigma-1)(y)=(\sigma-1)(y),
$$
since $a\in K^\times$.  Therefore $w(\sigma(x)-x)=w(\sigma(y)-y)$, which
implies that 
$$
w(\sigma(y)-y)-w(y)<w(\sigma(x)-x)-w(x)\leqno{(2)}
$$
since $w(y)>w(x)$. In other words, for every $x\in L^\times$, we have found
a $y\in L^\times$ for which $(2)$ holds.  This implies that if $\sigma\in
G_{[r}$, then there is no $x\in L^\times$ such that $w(\sigma(x)-x)-w(x)=r$.
In other words, we have proved the following proposition.

\th PROPOSITION 16
\enonce
If\/ $e=1$ and\/ $f=1$, then\/ $G_{[r}=G_{]r}$ for every\/ $r\in\R$.
\endth

\vfill\eject

\centerline{\bf Lecture 16}
\bigskip
\centerline{\it The upper numbering}

\bigskip

Let $K$ be a complete discretely valued field and $L|K$ a finite galoisian
extension, of group $G=\Gal(L|K)$.  Let $w$ be the valuation of $L$ such that
$w(L^\times)=\Z$.  For every integer $n\in\{-1\}\cup\N$, we have the subgroup
$G_{n}\subset G$, invariant under conjugation in $G$.  We have $G_{-1}=G$,
whereas $G_0$ is the inertia group (Lecture~14, Proposition~5) and $G_1$ the
ramification group (Lecture~15, Proposition~14 and Lemma~13).

Recall that $G_n$ is the set of those $\sigma\in G$ which operate trivially on
${\goth O}/{\goth P}^{n+1}$, or equivalently, for which $w(\sigma(x)-x)\ge
n+1$ for every $x\in{\goth O}$.  Suppose that the residual extension $l|k$ is
separable, and write ${\goth O}=\ogoth[\omega]$ for some $\omega\in{\goth O}$
(Lecture~15, Theorem~8), then $G_n$ is the set of those $\sigma\in G$ for which
$w(\sigma(\omega)-\omega)\ge n+1$ (Lecture~15, Corollary~9).  

We have the short exact sequences $1\to G_0\to G\to\Gal(l|k)\to1$ and
$$
%1\to G_1\to G_0\to\Hom(w(L^\times)/v(K^\times),l^\times)\to1,
1\to G_1\to G_0\to\Hom(\Z/e'\Z,l^\times)\to1,
$$ 
where $e'$ is the prime-to-$p$ part of the ramification index $e=e'p^r$ of
$L|K$ (Lecture~14, Proposition~8), and $p$ is the characteristic exponent of
$k$~; the order of $G_0$ (resp.~$G_1$) is $e$ (resp.~$p^r$).

\th THEOREM 1
\enonce
The inertia group\/ $G_0$ is solvable.  If the residue field\/ $k$ is finite,
then so is\/ $G$.
\endth
As we have just recalled, the group $G_0/G_1$ is cyclic, and hence solvable.
We have also shown that $G_1=G'$ is a $p$-group (Proposition~7, Lecture~14)
and hence solvable. Being an extension of solvable groups, $G_0$ is solvable.
If moreover $k$ is finite, the group $G/G_0=\Gal(l|k)$ is cyclic, and hence
solvable.  It follows as before that $G$ is solvable, in the case at hand.

\medbreak

For $\sigma\in G$, define $i_G(\sigma)=w(\sigma(\omega)-\omega)$~; this
integer does not depend on the choice of $\omega$ (Lecture~15, Corollary~9).
The function $i_G:G\to\N$ determines the filtration on $G$~; indeed, we have
$\sigma\in G_n\Leftrightarrow i_G(\sigma)\ge n+1$, Also, because $G_n\subset G$
is a normal subgroup,
$$
i_G(\sigma^{-1})=i_G(\sigma),\quad
i_G(\sigma\tau)\ge\inf(i_G(\sigma),i_G(\tau)),\quad 
i_G(\tau\sigma\tau^{-1})=i_G(\sigma),
$$
for every $\sigma,\tau\in G$.

\medbreak

Let $H\subset G$ be a subgroup of $G$, so that $H=\Gal(L|L^H)$, and we have
the subgroups $H_n\subset H$ for $n\in\{-1\}\cup\N$.  Notice that
$i_H=i_G|_H$. 

\th PROPOSITION 2
\enonce
For every\/ $n$, we have\/ $H_n=G_n\cap H$.
\endth
For every $\sigma\in H$, we have $\sigma\in H_n\Leftrightarrow i_H(\sigma)\ge
n+1\Leftrightarrow i_G(\sigma)\ge n+1\Leftrightarrow \sigma\in G_n\cap H$.

It follows, in particular, that if $H=G_0$, then $H_n=G_n$ for every
$n\in\N$. 

Now suppose that the subgroup $H\subset G$ is normal, so that the extension
$L^H$ of $K$ galoisian, of group $G/H=\Gal(L^H|K)$.  The filtration $(G/H)_n$
is determined by the filtration $(G)_n$, as the following proposition shows.

\th PROPOSITION 3
\enonce
Let\/ $e_{L|L^H}$ be the ramification index of\/ $L|L^H$.  Then,
for  every\/ $s\in G/H$, we have\/  
$$
i_{G/H}(s)={1\over e_{L|L^H}}
 \sum_{\sigma\mapsto s}i_G(\sigma)\qquad(\sigma\in G).
$$
\endth
This is clear if $s=\Id_{L^H}$, in which case each side equals $+\infty$.
Suppose that $s\neq\Id_{L^H}$.  Denote by $u$ the valuation of
$L^H$ such that $u(L^{H\times})=\Z$~; we have $w=e_{L|L^H}.u$. 

Let $\eta$ be a $u$-integer of $L^H$ such that $\ogoth[\eta]$ is the
ring of $u$-integers of $L^H$ (Lecture~15, Theorem~8).  We have
$i_{G/H}(s)=u(s(\eta)-\eta)$, so that $e_{L|L^H}.i_{G/H}=w(s(\eta)-\eta)$.

Fix a $\sigma\in G$ such that $\sigma\mapsto s$, which means that
$\sigma|_{L^H}=s$.  Then $\sigma H$ is the set of $\tilde s\in G$ with image
$s$ in $G/H$.  So everything comes down to showing that the elements
$$
a=\sigma(\eta)-\eta,\qquad b=\prod_{\tau\in H}(\sigma\tau(\omega)-\omega)
$$
generate the same ideal of the ring $\goth O=\ogoth[\omega]$ of $w$-integers
of $L$.  This is what we will show --- that $a|b$ and $b|a$.

Let $f$ be the minimal polynomial of $\omega$ over $L^H$~; it has coefficients
in $\ogoth[\eta]$, and in fact $f=\prod_{\tau\in H}(T-\tau(\omega))$.
Applying $\sigma$ to the coefficients of $f$, we get $f^\sigma=\prod_{\tau\in
  H}(T-\sigma\tau(\omega))$.  Writing $f=\sum_ic_i(\eta)T^i$ for some
$c_i\in\ogoth[T]$, we see that the coefficients of $f^\sigma-f$ are all
divisible by $a=\sigma(\eta)-\eta$, because
$f^\sigma=\sum_ic_i(\sigma(\eta))T^i$.  Hence $a$ divides
$f^\sigma(\omega)-f(\omega)=\pm b$.

Conversely, to show that $b$ divides $a$, write $\eta=g(\omega)$ for some
$g\in\ogoth[T]$.  As $\omega$ is a root of the polynomial
$g-\eta\in\ogoth[\eta][T]$, we may write $g-\eta=fh$ for some
$h\in\ogoth[\eta][T]$ (Lemma~3, Lecture~8), because $f$ is the minimal
polynomial of $\omega$ over $L^H$.  Applying $\sigma$ and making the
substitution $T\mapsto\omega$, we see that
$$
a=\sigma(\eta)-\eta
=f^\sigma(\omega).h^\sigma(\omega)-f(\omega).h(\omega)
=\pm b.h^\sigma(\omega).
$$
As $h^\sigma(\omega)$ is in $\goth O$, we conclude that $b$ divides $a$, 
completing the proof. 

\th COROLLARY 4
\enonce
If\/ $H=G_m$ for some\/ $m\in\N$, then\/ $(G/H)_n=G_n/H$ for\/ $n\le m$ and\/
$(G/H)_n=\{1\}$ for $n\ge m$.
\endth
The subgroups $G_n/H$ ($n\le m$) form a derceasing filtration of $G/H$.  For
every $s\neq1$ in $G/H$, there is a unique index $n<m$ such that $s\in G_n/H$
but $s\notin G_{n+1}/H$.  If $\sigma\in G$ lifts $s$, then it is clear that
$\sigma\in G_n$ but $\sigma\notin G_{n+1}$, therefore $i_G(\sigma)=n+1$.
Using Proposition~3 we get $i_{G/H}(s)=n+1$ because, as $H\subset G_0$, the
extension $L|L^H$ is totally ramified of degree $e_{L|L^H}=\Card H$.  Hence
$G_n/H$ coincides with $(G/H)_n$ for $n\le m$.  In particular, we have
$(G/H)_m=\{1\}$, and hence $(G/H)_n=\{1\}$ for every $n\ge m$.

\medbreak
\centerline{\hbox to5cm{\hrulefill}}
\bigbreak

For every real $u\in[-1,+\infty[$, define the ramification group $G_u$ to be
$G_i$, where $i=\lceil u\rceil$ is the smallest integer $\ge u$, so that
$\sigma\in G_u\Leftrightarrow i_G(\sigma)\ge u+1$ for every $u$.  Define
$$
\varphi_{L|K}(u)=\int_0^u{dt\over(G_0:G_t)},
$$
where, for $t\in[-1,0]$, we define $(G_0:G_t)$ to be $1/(G_t:G_0)$, so that it
is equal to $1/(G_{-1}:G_0)$ if $t=-1$, and equal to $1=1/(G_t:G_0)$ if
$t\in\;]-1,0]$.

More explicitly, $\varphi_{L|K}(u)=u$ for $u\in[-1,0]$, and, for every
$m\in\N$ and every $u\in[m,m+1]$,
$$
\varphi_{L|K}(u)={1\over g_0}(g_1+g_2+\cdots+(u-m)g_{m+1}),\quad
g_i=\Card G_i.
$$
In particular, $\varphi_{L|K}(m)+1=(1/g_0)\sum_{i=0}^mg_i$ for evey $m\in\N$.
\th PROPOSITION 5
\enonce
The function\/ $\varphi_{L|K}$ is continuous, piecewise linear, increasing and
concave.  We have\/ $\varphi_{L|K}(0)=0$.  It is differentiable at every\/
$u\notin\{-1\}\cup\N$, with derivative\/ $1/(G_0:G_u)$~; at\/ $u=-1$, the right
derivative is\/ $1$~; at\/ $u\in\N$, the left and right derivatives are
respectively 
$$
{1\over(G_0:G_u)},\quad {1\over(G_0:G_{u+1})}.
$$
\endth
This follows directly from the definition of $\varphi_{L|K}$.  Note that
$\varphi_{L|K}$ is a homeomorphism of $[-1,+\infty[$ onto itself.  Denote by
$\psi_{L|K}=\varphi_{L|K}^{-1}$ the inverse.
\th PROPOSITION 6
\enonce
The function\/ $\psi_{L|K}$ is continuous, piecewise linear, increasing and
convex.  We have\/ $\psi_{L|K}(0)=0$.  If\/ $v=\varphi_{L|K}(u)$ for some
$u\notin\{-1\}\cup\N$, then $\psi_{L|K}$ is differentiable with derivative
$(G_0:G_u)$.  If\/ $u=-1$, the right derivative at $v=-1$ is\/ $1$~; if\/
$u\in\N$, the left and right derivatives at $v$ are respectively $(G_0:G_u)$
and $(G_0:G_{u+1})$.  Finally, if $v\in\N$ is an integer, then so is
$u=\psi_{L|K}(v)$.
\endth
Only the last statement needs a proof.  Let $m$ be an integer such that
$u\in[m,m+1]$.  Then, by definition,
$$
g_0v=g_1+g_2+\cdots+g_m+(u-m)g_{m+1}.
$$
As $G_{m+1}$ is a subgroup of the groups $G_0$, $G_1$, $\ldots$, $G_m$, its
order $g_{m+1}$ divides $g_0$, $g_1$, $\ldots$, $g_m$.  Therefore $u-m$ is an 
integer, and hence so is $u$.

\th PROPOSITION 7
\enonce
We have\/ $\displaystyle\varphi_{L|K}(u)+1={1\over g_0}\sum_{\sigma\in
  G}\inf(i_G(\sigma),u+1)$. 
\endth
Denote by $\theta(u)$ the function on the right~; it is a continuous,
piecewise linear function such that $\theta(0)=1$.  For $m\in\{-1\}\cup\N$ and
$u\in\,]m,m+1[$, the derivative $\theta'(u)$ is equal to the number of
$\sigma\in G$ with $i_G(\sigma)\ge m+2$, multiplied by $1/g_0$, so
$\theta'(u)=1/(G_0:G_{m+1})$.  As the function on the left has the same value
at $u=0$, and the same derivative at every $u\notin\{-1\}\cup\N$
(Proposition~5), the two functions must be equal.
\th PROPOSITION 8
\enonce
Let\/ $s\in G/H$, and let\/ $j(s)=\sup_{\sigma\mapsto s}i_G(\sigma)$
($\sigma\in G$).  Then\/ $i_{G/H}(s)-1=\varphi_{L|L^H}(j(s)-1)$.
\endth
Suppose that the supremum is attained at $\sigma\in G$, so that $\sigma\mapsto
s$ and $j(s)=i_G(\sigma)$~; put $m=i_G(\sigma)$.  Let $\tau\in H$.  If
$\tau\in H_{m-1}$, then $i_G(\tau)\ge m$, hence $i_G(\sigma\tau)\ge m$, so
$i_G(\sigma\tau)=m$, because $m$ is the supremum.  If, on the other hand,
$\tau\notin H_{m-1}$, then $i_G(\tau)<m$, and $i_G(\sigma\tau)=i_G(\tau)$.  In
either case, $i_G(\sigma\tau)=\inf(i_G(\tau),m)$.  Applying Proposition~3, we
get
$$
i_{G/H}(s)={1\over e_{L|L^H}}
 \sum_{\tau\in H}\inf(i_G(\tau),m).
$$
But observe that $i_G(\tau)=i_H(\tau)$ and that $e_{L|L^H}=\Card H_0$.
Applying Proposition~7 to the group $H$, we get
$i_{G/H}(s)=\varphi_{L|L^H}(m-1)+1$, which was to be proved.
\th PROPOSITION 9 (Herbrand)
\enonce
If\/ $v=\varphi_{L|L^H}(u)$, then\/ $G_uH/H=(G/H)_v$.
\endth
In the notation of Proposition~8, we have, for $s\in G/H$,
$$\eqalign{
s\in G_uH/H
&\Leftrightarrow j(s)-1\ge u\cr
&\Leftrightarrow \varphi_{L|L^H}(j(s)-1)\ge\varphi_{L|L^H}(u)
            \quad\hbox{(Proposition~5)}\cr
&\Leftrightarrow
i_{G/H}(s)-1\ge\varphi_{L|L^H}(u)\quad\hbox{(Proposition~8)}\cr 
&\Leftrightarrow s\in (G/H)_v.\cr
}$$
\vskip-\baselineskip
\th PROPOSITION 10
\enonce
We have the transitivity formulae\/ 
$$
\varphi_{L|K}=\varphi_{L^H|K}\circ\varphi_{L|L^H},\quad
\psi_{L|K}=\psi_{L|L^H}\circ\psi_{L^H|K}.
$$
\endth
The derivative of the function on the right in the first purported equality,
at any $u>-1$ not an integer, is $\varphi_{L^H|K}'(v)\varphi_{L|L^H}'(u)$, with
$v=\varphi_{L|L^H}(u)$.  By Propositions~5 and~9, it can therefore be written
as 
$$
{\Card(G/H)_v\over e_{L^H|K}}{\Card H_u\over e_{L|L^H}}
={\Card G_u\over e_{L|K}},
$$
which is also the derivative of $\varphi_{L|K}$ at $u$.  Hence the first
equality, of which the second is a consequence.

\medbreak

One defines the ramification filtration in the {\it upper numbering\/} by
$G^v=G_{\psi(v)}$ for $v\in[-1+\infty[$, or equivalently, by
$G^{\varphi(u)}=G_u$.  We have $G^{-1}=G$, $G^0=G_0$, and $G^v=\{1\}$ for $v$
sufficiently large.  The filtration in the {\it lower numbering\/} $(G_u)_u$
can be recovered from $(G^v)_v$ because
$$
\psi(v)=\int_0^v(G^0:G^w)\,dw.
$$
\vskip-\baselineskip
\th THEOREM 11
\enonce
We have\/ $G^vH/H=(G/H)^v$ for every\/ $v\in[-1,+\infty[$.
\endth
We have $(G/H)^v=(G/H)_x$, with $x=\psi_{L^H|K}(v)$.  By Proposition~9, we
have $(G/H)_x=G_wH/H$, with $w=\psi_{L|L^H}(x)=\psi_{L|K}(v)$
(Proposition~10).  Therefore $G_w=G^v$, and the proof is over.

\vfill\eject

\centerline{\bf Lecture 17}
\bigskip
\centerline{\it The different and the discriminant}

\bigskip

Let $K$ be a field complete for a discrete valuation $v$, and let $L|K$ be a
finite {\it separable\/} extension~; the $K$-linear trace map $S_{L|K}:L\to K$
is surjective, and we have a $K$-linear isomorphism
$$
x\mapsto(y\mapsto S_{L|K}(xy)):L\to\Hom_K(L,K).\leqno{(1)}
$$

Recall that a sub-$\goth O$-module ${\goth A}\subset L$ is a {\it fractional
  ideal\/} if ${\goth A}\neq0,L$.  Every fractional ideal is of the form
${\goth P}^n$ for some $n\in\Z$~; we say that the fractional ideal ${\goth
  P}^{-n}$ is the {\it inverse\/} of ${\goth P}^n$.  For every fractional
ideal ${\goth A}$, we have $K{\goth A}=L$.
  
An $\ogoth$-{\it lattice\/} in $L$ is a finitely generated sub-$\ogoth$-module
$\goth A$ such that $K{\goth A}=L$~; every fractional ideal is an
$\ogoth$-lattice.  For such a lattice, consider the subset ${}^*{\goth A}$ of
all $x\in L$ such that $S_{L|K}(x{\goth A})\subset\ogoth$~; it is called the
{\it complementary lattice\/} of $\goth A$ --- it is clearly an
$\ogoth$-lattice, and indeed a fractional ideal if $\goth A$ is a fractional
ideal.  If ${\goth A}\subset{\goth B}$, then ${}^*{\goth B}\subset{}^*{\goth
  A}$.  If $e_1,\ldots,e_n$ is a $K$-basis of $L$, $e_1^*,\ldots,e_n^*$ the
dual basis with respect to the bilinear form $(x,y)\mapsto S_{L|K}(xy)$, and
$\goth A$ the sub-$\ogoth$-module generated by the $e_i$, then $\goth A$ is an
$\ogoth$-lattice and the complementary lattice is the sub-$\ogoth$-module
generated by the $e_i^*$.

As an $\ogoth$-module, ${}^*{\goth A}$ is canonically isomorphic to
$\Hom_{\ogoth}({\goth A},\ogoth)$ by the map $x\mapsto(y\mapsto S_{L|K}(xy))$.
Indeed, because ${\goth A}\otimes_\ogoth K=L$, the map $f\mapsto
f\otimes\Id_K$ identifies $\Hom_{\ogoth}({\goth A},\ogoth)$ with a
sub-$\ogoth$-module of $\Hom_K(L,K)$ --- the image of $\goth A$ under the
$K$-isomorphism $L\to\Hom_K(L,K)$ (1).

\th DEFINITION 1
\enonce
The different\/ ${\goth D}_{L|K}=({\goth C}_{L|K})^{-1}$ of\/ $L|K$ is the
inverse of the complementary lattice\/ ${\goth C}_{L|K}={}^*{\goth O}$ of the
$\ogoth$-lattice\/ $\goth O$.   
\endth
Note that ${}^*{\goth O}$ is the largest sub-$\goth O$-module $E$ of $L$ such
that $S_{L|K}(E)\subset\ogoth$.  Also, as ${\goth O}\subset{}^*{\goth O}$, we
have ${\goth D}_{L|K}\subset{\goth O}$, so the different is an ideal of $\goth
O$. 
\th LEMMA 2
\enonce
If\/ $M|L$ is also finite separable, then\/ 
${\goth D}_{M|K}={\goth D}_{M|L}{\goth D}_{L|K}$.
\endth
This is equivalent to ${\goth C}_{M|K}={\goth C}_{M|L}{\goth C}_{L|K}$.
As $S_{M|K}=S_{L|K}\circ S_{M|L}$, we have
$$\eqalign{
S_{M|K}({\goth C}_{M|L}{\goth C}_{L|K})
&=S_{L|K}(S_{M|L}({\goth C}_{M|L}{\goth C}_{L|K}))\cr
&=S_{L|K}({\goth C}_{L|K}\, S_{M|L}({\goth C}_{M|L}))\cr
&\subset S_{L|K}({\goth C}_{L|K}\,{\goth O}_L)\subset\ogoth,\cr
}$$
and hence the inclusion ${\goth C}_{M|L}{\goth C}_{L|K}\subset{\goth
  C}_{M|K}$.  Conversely, as $\ogoth$ contains $S_{M|K}({\goth
  C}_{M|K})=S_{L|K}(S_{M|L}({\goth C}_{M|K}))$, we have $S_{M|L}({\goth
  C}_{M|K})\subset{\goth C}_{L|K}$.  In other words, $S_{M|L}({\goth
  C}_{L|K}^{-1}{\goth C}_{M|K})\subset{\goth O}_L$.  This means that ${\goth
  C}_{L|K}^{-1}{\goth C}_{M|K}\subset{\goth C}_{M|L}$, which gives the
other inclusion ${\goth C}_{M|K}\subset{\goth C}_{M|L}{\goth C}_{L|K}$, and
hence their equality.

\th PROPOSITION 3
\enonce
Suppose that\/ ${\goth O}=\ogoth[\alpha]$ for some\/ $\alpha\in{\goth O}$.
Then\/ ${\goth D}_{L|K}=f'(\alpha){\goth O}$, where\/ $f\in\ogoth[T]$ is the
minimal polynomial of\/ $\alpha$.
\endth
Write $f=a_0+a_1T+\cdots+a_nT^n$, and define the $b_i\in{\goth O}$
($i\in[0,n[$) by 
$$
{f(T)\over T-\alpha}=b_0+b_1T+\cdots+b_{n-1}T^{n-1}.
$$
We claim that, with respect to the bilinear form $(x,y)\mapsto S_{L|K}(xy)$,
the $K$-basis of $L$ dual to the basis $1,\alpha,\ldots,\alpha^{n-1}$ is
$b_0/f'(\alpha),\ldots,b_{n-1}/f'(\alpha)$. 

Indeed, denoting by $\alpha_1,\ldots,\alpha_n$ the $n$ roots of $f$, we have,
for $r\in[0,n[$,  
$$
\sum_{i=1}^n{f(T)\over T-\alpha_i}{\alpha_i^r\over f'(\alpha_i)}=T^r,
$$
because the difference of the two sides is a polynomial of degree $<n$ having
$\alpha_1,\ldots,\alpha_n$ as roots, and hence must be identically zero.  The
above equation can be rewritten as 
$$
S_{L|K}\!\left({f(T)\over T-\alpha}{\alpha^r\over f'(\alpha)}\right)
=T^r,\quad(r\in[0,n[).
$$
As $f(T)/(T-\alpha)=b_0+b_1T+\cdots+b_{n-1}T^{n-1}$, we conclude with Euler
that  
$$
S_{L|K}\!\left(\alpha^i{b_j\over f'(\alpha)}\right)
=\cases{1&if $i=j$\cr
 0&if $i\neq j$,\cr}
$$
thereby establishing the claim that
$b_0/f'(\alpha),\ldots,b_{n-1}/f'(\alpha)$ is the $K$-basis of~$L$ dual to
$1,\alpha,\ldots,\alpha^{n-1}$, with respect to the trace form.

From the fact that ${\goth O}=\ogoth+\ogoth\alpha+\cdots+\ogoth\alpha^{n-1}$,
it now follows that ${\goth
  C}_{L|K}=f'(\alpha)^{-1}{\goth B}$, where ${\goth
  B}=\ogoth b_0+\ogoth b_1+\cdots+\ogoth b_{n-1}$ is the sub-$\ogoth$-module
generated by the $b_i$.  To complete the proof, it remains to show that $\goth
B=O$. 

This follows from the sequence of relations $b_{n-1}=1$,
$b_{n-2}=\alpha+a_{n-1}$, and generally, 
$$
b_{n-i}=\alpha^{i-1}+a_{n-1}\alpha^{i-2}+\cdots+a_{n-i+1}\quad(i\in[1,n]).
$$
Therefore the $\alpha^i$ and the $b_j$ generate the same
sub-$\ogoth$-module of $L$, and hence $\goth B=O$ and ${\goth
  D}_{L|K}={\goth C}_{L|K}^{-1}=f'(\alpha){\goth O}$.

\th THEOREM 4
\enonce
Let\/ ${\goth D}_{L|K}={\goth P}^s$ be the different of the extension\/
$L|K$, and suppose that\/ $l|k$ is separable.  Then\/ $s=0$ if and only if\/
$L|K$ is unramified.  More generally, $s=e-1$ if and only if\/ $L|K$ is tamely
ramified.  Otherwise, $s\in[e,e+w(e)[$, where $w$ is the valuation on\/ $L$
such that $w(L^\times)=\Z$.
\endth

Write ${\goth O}=\ogoth[\alpha]$ for a suitable $\alpha\in{\goth O}$
(Lecture~15, Theorem~8)~; denoting by $f$ $(\in\ogoth[T])$ the minimal
polynomial of $\alpha$, we know that $s=w(f'(\alpha))$ (Proposition~3).

If $L|K$ is unramified, then $\hat\alpha\in l$ is a simple root of $\hat f\in
k[T]$, therefore $f'(\alpha)\in{\goth O}^\times$, hence $s=0$ and conversely.

This allows us to replace $L|K$ by $L|L_0$ in general, by the transitivity of
the different (Lemma~2).  In other words, we may suppose that $L|K$ is totally
ramified of degree~$e$, and that $\alpha$ is a uniformiser of $L$ (Lecture~15,
Proposition~10).  Its minimal polynomial $f=a_eT^e+\cdots+a_1T+a_0$ ($a_e=1$)
is then Eisenstein, so that $v(a_i)>0$ for $i\in[1,e[$, and $a_0$ is a
uniformiser of~$K$.  We have
$$
f'(\alpha)=ea_e\alpha^{e-1}+(e-1)a_{e-1}\alpha^{e-2}+\cdots+a_1.
$$
For $i\in[0,e[$, we have
$w((e-i)a_{e-i}\alpha^{e-i-1})=w(e-i)+w(a_{e-i})+(e-i-1)$, which is $\equiv
-i-1\pmod e$ because $e-i,a_{e-i}$ are in $K^\times$ and hence
$w(e-i),w(a_{e-i})\in e\Z$, unless $(e-i)a_{e-i}=0$.  This means that any two
different terms in the above expression for $f'(\alpha)$ have different
valuations, unless both terms vanish.  Therefore
$$
s=w(f'(\alpha))=\inf_{i\in[0,e[}(w((e-i)a_{e-i}\alpha^{e-i-1})).
$$
Now it is clear that the infimum is $e-1$ if and only if $e$ is prime to~$p$,
which is the same as saying that $L|K$ is tamely ramified. Otherwise, when
$v(e)>0$, if the infimum $s$ occurs at the first term $a_1$, then
$s=w(a_1)=ev(a_1)$ is $>e-1$~; if the infimum occurs at the last term
$ea_e\alpha^{e-1}$, then $s=e-1+w(e)$~; otherwise $s$ lies between these
two values.  So we get $s\in[e,e+w(e)[$, as claimed, if $L|K$ is wildly
ramified.  

\medbreak
\centerline{\hbox to5cm{\hrulefill}}
\bigbreak

Recall that given a commutative ring $A$ and an $A$-algebra $B$, an $A$-{\it
  derivation\/} on $B$ is a $B$-module $M$ together with an $A$-linear map
$d:B\to M$ such that $d(xy)=xd(y)+yd(x)$ for every $x,y\in B$.  It follows
that $d(a)=0$ for every $a\in A$.

An $A$-derivation $d_{B|A}:B\to\Omega^1_{B|A}$ is said to be {\it universal\/}
if for every $A$-derivation $d:B\to M$, there exists a unique $B$-linear map
$f:\Omega^1_{B|A}\to M$ such that $d=f\circ d_{B|A}$.  If a universal
$A$-derivation on $B$ exists, it is unique, up to a unique $B$-isomorphism.

As for the existence of the universal $A$-derivation on $B$, it can be
verified that if $E$ is the free $B$-module on the symbols $d(x)$ ($x\in B$)
and $R$ the sub-$B$-module generated by
$$
d(a),\quad d(x+y)-d(x)-d(y),\quad d(xy)-xd(y)-yd(x),\quad (a\in A,\;x,y\in B),
$$
then the $B$-module $\Omega^1_{B|A}=E/R$, together with the 
map $d_{B|A}:B\to\Omega^1_{B|A}$ which sends $x\in B$ to the
image of $d(x)\in E$ in $\Omega^1_{B|A}$, is universal.

\smallskip

Let us return to our extension $L|K$ and suppose that the residual extension
$l|k$ is {\it separable}.
\th PROPOSITION 5
\enonce
The different\/ ${\goth D}_{L|K}$ is the annihilator of the\/ $\goth
O$-module\/ $\Omega^1_{{\goth O}|\ogoth}$~: we have\/ ${\goth
  D}_{L|K}=\Ann_{\goth O}\Omega^1_{{\goth O}|\ogoth}$. 
\endth
Indeed, write ${\goth O}=\ogoth[x]$ for a suitable $x\in{\goth O}$
(Lecture~15, Theorem~8), and let $f$ ($\in\ogoth[T]$) be the minimal polynomial
of $x$.  Then the $\goth O$-module $\Omega^1_{{\goth O}|\ogoth}$ is generated
by $d_{{\goth O}|\ogoth}(x)$, whose annihilator is therefore the ideal in
$\goth O$  generated by $f'(x)$, namely the different ${\goth D}_{L|K}$
(Proposition~3). 

\medbreak
\centerline{\hbox to5cm{\hrulefill}}
\bigbreak

Now assume, in addition, that the extension $L|K$ is galoisian of group
$G=\Gal(L|K)$.  We wish to express the exponent $s=w({\goth D}_{L|K})$ in the
different ${\goth D}_{L|K}={\goth P}^s$ of $L|K$ in terms of the function
$i_G:G\to\N$ and thereby in terms of the ramification filtration $(G_n)_n$ on
the group $G$.

Recall that $G_n=\{1\}$ for large $n$ (Lemma~3 and
Proposition~14, Lecture~15), so the sum $\sum_{n\in\N}(\Card G_n-1)$ makes
sense. 
\th PROPOSITION 6
\enonce
We have\/ $\displaystyle w({\goth D}_{L|K})
 =\sum_{\sigma\neq1}i_G(\sigma)=\sum_{n\in\N}(\Card G_n-1)$.
\endth
Write ${\goth O}=\ogoth[x]$ for a suitable $x\in{\goth O}$ (Lecture~15,
Theorem~8), and let $f$ ($\in\ogoth[T]$) be the minimal polynomial of $x$~; we
have $f=\prod_{\sigma\in G}(T-\sigma(x))$ and hence
$$
f'(x)=\prod_{\sigma\neq1}(x-\sigma(x)).
$$
As we have $w({\goth D}_{L|K})=w(f'(x))$ (Proposition~3), the above expression
gives $w({\goth D}_{L|K})=\sum_{\sigma\neq1}w(x-\sigma(x))$, which is the
first equality, in view of the definition $i_G(\sigma)=w(x-\sigma(x))$.

Put $r_n=\Card G_n-1$.  For every $n\in\N$, the function $i_G:G\to\N$ takes
the value $n$ on the set $G_{n-1}-G_n$, which has $r_{n-1}-r_n$ elements.
Therefore
% $$\eqalign{
% \sum_{\sigma\neq1}i_G(\sigma)=\sum_{n\in\N}n(r_{n-1}-r_n)
% &=(r_0-r_1)+2(r_1-r_2)+3(r_2-r_3)+\cdots\cr
% &=r_0+r_1+r_2+\cdots,\cr
% }$$
$$
\sum_{\sigma\neq1}i_G(\sigma)=\sum_{n\in\N}n(r_{n-1}-r_n)
=r_0+r_1+r_2+\cdots=\sum_{n\in\N}(\Card G_n-1),
$$
thereby proving the second equality and completing the proof of the
proposition. 

\th COROLLARY 7
\enonce
Let\/ $H\subset G$ be any subgroup, $K'=L^H$ the fixed field,  and\/ $v'$ the
normalised valuation of\/ $K'$.  Then\/ 
$\displaystyle e_{L|K'}.v'({\goth D}_{K'|K})=\sum_{\sigma\notin H}i_G(\sigma)$.
\endth
Indeed, Proposition~6, applied to the galoisian extensions $L|K$ and $L|K'$,
of groups $G$ and $H$ respectively, gives, in view of the fact that
$i_H=i_G|_H$,  
$$
w({\goth D}_{L|K})=\sum_{\sigma\neq1,\sigma\in G}i_G(\sigma),\quad
w({\goth D}_{L|K'})=\sum_{\sigma\neq1,\sigma\in H}i_G(\sigma).
$$
These, togther with the relation $w=e_{L|K'}.v'$ and the transitivity formula  
${\goth D}_{L|K}={\goth D}_{L|K'}{\goth D}_{K'|K}$ (Lemma~2), give the desired
result. 

\medbreak
\centerline{\hbox to5cm{\hrulefill}}
\bigbreak

Let $K$ be any field and $L$ a finite separable extension of $K$.  For a
$K$-basis $b_1,\ldots,b_n$ of $L$, consider the determinant
$$
d_{L|K}(b_1,\ldots,b_n)=\det(S_{L|K}(b_ib_j)_{i,j\in[1,n]})~;
$$
it is $\neq0$ because $L|K$ is separable.  If we change the basis of
$L$ to another $K$-basis $c_1,\ldots,c_n$, and if $\hbox{\it
  M\/}=(m_{ij})_{i,j\in[1,n]}$ is the transition matrix, so that
$c_i=\sum_jm_{ij}b_j$, then
$$
d_{L|K}(c_1,\ldots,c_n)=\det(\hbox{\it M\/})^2d_{L|K}(b_1,\ldots,b_n),
$$
because $S_{L|K}$ is $K$-linear.  Let $\bar L$ be a separable algebraic
closure of $L$, and $\sigma_1,\ldots,\sigma_n$ the $K$-morphisms $L\to\bar
L$.  Then it can checked that 
$$
d_{L|K}(b_1,\ldots,b_n)=(\det(\sigma_i(b_j))^2.
$$
In particular, if $b_1^*,\ldots,b_n^*$ is the $K$-basis of $L$ dual to
$b_1,\ldots, b_n$, then
$$
d_{L|K}(b_1,\ldots,b_n).d_{L|K}(b_1^*,\ldots,b_n^*)=1.
$$
Now suppose that $K$ is complete for a discrete valuation, and $\Lambda$ is an 
$\ogoth$-lattice in $L$.  If we take $\ogoth$-bases $(b_i)_i$, $(c_i)_i$ of
$\Lambda$, then we have $\det\hbox{\it M\/}\in\ogoth^\times$, as both
$\hbox{\it M\/}$ and $\hbox{\it M\/}^{-1}$ have entries in $\ogoth$.
Consequently, the sub-$\ogoth$-mdoule of $K$ generated by
$d_{L|K}(b_1,\ldots,b_n)$ depends only $\Lambda$, not on the choice of the
$\ogoth$-basis for $\Lambda$~; denote it ${\goth d}_\Lambda$.  The above
equality implies that ${\goth d}_{\Lambda}{\goth d}_{{}^*\Lambda}=\ogoth$.

\th DEFINITION 8
\enonce
The discriminant\/ ${\goth d}_{L|K}$ of\/ $L|K$ is the
sub-$\ogoth$-mdoule of\/ $K$ generated by\/ $d_{L|K}(b_1,\ldots,b_n)$, where\/
$b_1,\ldots,b_n$ is any $\ogoth$-basis of\/ $\goth O$.   
\endth
Clearly, ${\goth d}_{L|K}$ is an ideal of $\ogoth$, so that ${\goth
  d}_{L|K}=\pgoth^d$ for some $d\in\N$, called the {\it exponent\/} of the
discriminant.

\th PROPOSITION 9
\enonce
If\/ ${\goth D}_{L|K}={\goth P}^s$ is the different of\/ $L|K$, and\/
$f=[l:k]$ the residual degree, then\/ ${\goth d}_{L|K}={\goth p}^{sf}$.
\endth
Indeed, if $b_1,\ldots,b_n$ is an $\ogoth$-basis of $\goth O$, and $\Pi$ a
uniformiser of $L$, then $\Pi^{-s} b_1,\ldots,\Pi^{-s} b_n$ is an
$\ogoth$-basis of ${\goth C}={}^*{\goth O}$.  Denoting by $\hbox{\it
  M\/}$ the transition matrix, we have 
$$
d_{L|K}(\Pi^{-s}b_1,\ldots,\Pi^{-s}b_n)
=(\det\hbox{\it M\/})^2d_{L|K}(b_1,\ldots,b_n),
$$
so ${\goth d}_{\goth C}={\goth t}^2{\goth d}_{L|K}$, where $\goth t$ is the
ideal generated by $\det\hbox{\it M\/}=N_{L|K}(\Pi)^{-s}$.  But observe that
${\goth d}_{\goth C}={\goth d}_{L|K}^{-1}$, because $\goth C$ admits the
$\ogoth$-basis $b_1^*,\ldots,b_n^*$, and that ${\goth t}=\pgoth^{-sf}$, because
$N_{L|K}(\Pi)\ogoth=\pgoth^f$, as follows from the relations
$$
nw(\Pi)=v(N_{L|K}(\Pi)),\quad n=ef.
$$
Substituting these values in ${\goth d}_{\goth C}={\goth t}^2{\goth d}_{L|K}$,
we get ${\goth d}_{L|K}^{-1}=\pgoth^{-2sf}{\goth d}_{L|K}$, which was to be
proved. 

The above proposition is usually expressed as by saying that the norm of the
different is the discriminant, and written
$$
N_{L|K}({\goth D}_{L|K})={\goth d}_{L|K},
$$ 
where the definition $N_{L|K}({\goth P})=\pgoth^f$ is extended by
multiplicativity to all fractional ideals of $\goth O$, so that
$N_{L|K}(x{\goth O})=N_{L|K}(x)\ogoth$ for every $x\in L^\times$.

\th COROLLARY 10
\enonce
Suppose that\/ $l|k$ is separable.  The extension\/ $L|K$ is unramified if and
only if\/ ${\goth d}_{L|K}=\ogoth$.
\endth
Indeed, the exponent $s$ of the different is $0$ if and only if $L|K$ is
unramified (Theorem~4), and the exponent of the discriminant is $sf$
(Proposition~9). 

\th COROLLARY 11
\enonce
For a tower of finite separable extensions $M|L|K$, we have\/ ${\goth
  d}_{M|K}={\goth d}_{L|K}^{[M:L]}N_{L|K}({\goth d}_{M|L})$.
\endth
This follows from applying $N_{M|K}=N_{L|K}\circ N_{M|L}$ to the transitivity
formula ${\goth D}_{M|K}={\goth D}_{M|L}{\goth D}_{L|K}$ (Lemma~2) for the
different.  Indeed,
$$
{\goth d}_{M|K}=N_{L|K}({\goth d}_{M|L})N_{L|K}({\goth D}_{L|K}^{[M:L]})
 =N_{L|K}({\goth d}_{M|L}){\goth d}_{L|K}^{[M:L]}.
$$

\medbreak
\centerline{\hbox to5cm{\hrulefill}}
\bigbreak

Let us illustrate these ideas by computing the discriminant of the extension
of $\Q_p$ obtained by adjoining a primitive $p^n$-th root of~$1$.

\th PROPOSITION 12
\enonce
Let\/ $\zeta\in\bar\Q_p^\times$ be an element of order\/ $p^n$ ($n>0$).
Then\/ $\Q_p(\zeta)$ is a totally ramified galoisian extension of\/ $\Q_p$,
with group of automorphisms\/ $(\Z/p^n\Z)^\times$.  Moreover, $1-\zeta$ is a
uniformiser with norm\/ $p$, and the ring of integers is\/ $\Z_p[\zeta]$.
\endth
We already know that $\Q_p(\zeta)$ is totally ramified of degree
$p^n-p^{n-1}$ over $\Q_p$ (Example~10, Lecture~8, and Proposition~10,
Lecture~15).  Let's give a new proof, which proceeds by induction on $n$.  Put
$\xi_n=\zeta$, $\xi_{n-1}=\xi_n^p$, $\ldots$, $\xi_1=\xi_2^p$.  Also put
$K_i=\Q_p(\xi_i)$ and $\pi_i=1-\xi_i$.

Let us show that $\pi_1$ is a uniformiser of $K_1$, which is totally ramified
of degree $p-1$ over $\Q_p$.  As $\xi_1\in K_1^\times$ is an element of order
$p$, we have
$$
{1-\xi_1^p\over1-\xi_1}=1+\xi_1+\xi_1^2+\cdots+\xi_1^{p-1}=0,
$$
which, in terms of $\pi_1=1-\xi_1$, means that
$$
{1-(1-\pi_1)^p\over\pi_1}=p-{p\choose 2}\pi_1+\cdots+(-1)^{p-1}\pi_1^{p-1}=0. 
$$
Thus $\pi_1$ is a root of an Eisenstein polynomial, and hence
$K_1\,|\,\Q_p$ is totally ramified of degree~$p-1$, and $\pi_1$ is a
uniformiser of $K_1$ of norm $p$, so the ring of integers is
$\Z_p[\pi_1]=\Z_p[\xi_1]$.  Also, the extension is galoisian, and the 
natural map 
$$
\chi:\Gal(K_1|\Q_p)\to\Aut{}_pK_1^\times=(\Z/p\Z)^\times,\quad
\sigma(\xi_1)=\xi_1^{\chi(\sigma)}
$$ 
is injective, and hence an isomorphism because the two groups have the same
order $p-1$.  This establishes the case $n=1$.

Suppose that we have proved the proposition for $n-1$.  This means that
$K_{n-1}\,|\,\Q_p$ is totally ramified galoisian of
group~$(\Z/p^{n-1}\Z)^\times$, with $\pi_{n-1}=1-\xi_{n-1}$ as a uniformiser
of norm $p$ in $\Q_p$.  From $\xi_n^p-\xi_{n-1}=0$ it follows that
$$
(1-\pi_n)^p-(1-\pi_{n-1})=\pi_{n-1}-
{p\choose 1}\pi_n+
{p\choose 2}\pi_n^2-+\cdots 
+(-1)^p\pi_n^p=0,
$$
which means that $\pi_n$ is a root of a degree-$p$ Eisenstein polynomial
with coefficients in $K_{n-1}$.  Therefore $\pi_n$ is a uniformiser of $K_n$,
of norm $\pi_{n-1}$ in $K_{n-1}$, and hence of norm $p$ in $\Q_p$.  Also,
$K_n$ is totally ramified of degree~$p$ over $K_{n-1}$, hence totally ramified
of degree~$(p-1)p^{n-1}=p.(p-1)p^{n-2}$ over $\Q_p$, with ring of integers
$\Z_p[\pi_n]=\Z_p[\zeta]$.  As before, $K_n|\Q_p$ is galoisian, and the map
$\Gal(K_n|\Q_p)\to(\Z/p^n\Z)^\times$ which sends $\sigma$ to $\chi(\sigma)$ if
$\sigma(\zeta)=\zeta^{\chi(\sigma)}$ is an isomorphism because it is an
injective homomorphism and the groups have the same order~; the iniverse is
$a\mapsto(\zeta\mapsto\zeta^a)$.  This completes the proof by induction.

\smallskip

Now fix a prime number $p$, an integer $n>0$ and let $K=\Q_p$,
$L=\Q_p(\zeta)$, where $\zeta$ is a primitive $p^n$-th root of~$1$,
and $G=\Gal(L|K)=(\Z/p^n\Z)^\times$.  We wish to determine the ramification
filtration on $G$.

For every $s\in[0,n]$, denote by $G(s)\subset G$ the kernel of the canonical
projection $G\to(\Z/p^s\Z)^\times$, so that $G(0)=G$, $G(n)=\{1\}$, and, for
every~$s$, $G/G(s)=(\Z/p^s\Z)^\times$.  Let $w$ be the valuation on $L$ such
that $w(L^\times)=\Z$.

\th PROPOSITION 13
\enonce
The ramification filtration in the lower numbering\/
$(G_u)_{u\in\{-1\}\cup\N}$ on\/ $G$ is given by 
$$
\vbox{\halign{&\hfil$#$\hfil\quad\cr
u\in&\{-1,0\}&[1,p[&[p,p^2[&\cdots&[p^{n-1},p^n[\cr
\noalign{\vskip-5pt}
\multispan6\hrulefill.\cr
G_u=&G(0)&G(1)&G(2)&\cdots&G(n)\cr
}}
$$
\endth
For $a\in G$, denote the corresponding $K$-automorphism of $L$ by $\sigma_a$.
Suppose that $a\neq1$, and let $s$ be the largest integer such that
$a\in G(s)$~; we have $s\in[0,n[$.  Then
$$
i_G(\sigma_a)=w(\sigma_a(\zeta)-\zeta)=w(\zeta^a-\zeta)=w(\zeta^{a-1}-1).
$$

Notice that $\zeta^{a-1}$ is a primitive $p^{n-s}$-th root of\/~$1$, so
$\zeta^{a-1}-1$ is a uniformiser of $K(\zeta^{p^s})=L^{G(n-s)}$, of which $L$
is a totally ramified extension of degree $[L:L^{G(n-s)}]=\Card G(n-s)=p^s$.
Hence $i_G(\sigma_a)=p^s$.  

In other words, an $a\in G$ is in $G(m+1)$ for some $m\in[0,n[$ if and only if
$i_G(\sigma_a)>p^m$.  Now, if the integer $u$ is in $[p^m,p^{m+1}[$, then
$$
\sigma_a\in G_u\ \Leftrightarrow\ 
i_G(\sigma_a)>u\ \Leftrightarrow\ 
i_G(\sigma_a)>p^m\ \Leftrightarrow\
\sigma_a\in G(m+1)
$$
which shows that $G_u=G(m+1)$, as displayed in the table above. Cqfd.
\th COROLLARY 14
\enonce
The exponent of the different of $L|K$, as well as the exponent of the
discriminant, equals
$$\eqalign{
&(p^n-p^{n-1}-1)+(p-1)(p^{n-1}-1)+\cdots+(p^{n-1}-p^{n-2})(p-1)\cr
&\ \ \ =np^n-(n+1)p^{n-1}.\cr
}$$
\endth
For the different, this follows upon applying Proposition~6 to the table in
Proposition~13~; as the extension $L|K$ is totally ramified, the exponent of
the discriminant is the same.  Notice that it vanishes precisely when $n=1$,
$p=2$. 

\th COROLLARY 15
\enonce
Suppose that $n>1$, so that $L|K(\zeta^p)$ is cyclic of degree~$p$.  The
unique ramification break then occurs at $p^{n-1}-1$.
\endth

Indeed, we have $\Gal(L|K(\zeta^p))=G(n-1)$, whose ramification filtration (in
the lower numbering) is the restriction of the filtration $(G_u)_u$, by
Proposition~2, Lecture~16. 

\smallskip

Alternatively, $\sigma$ being a generator of $H=\Gal(L|K(\zeta^p))$, we have
$\sigma(\zeta)=\alpha\zeta$ for some order-$p$ element $\alpha\in L^\times$,
because $\zeta$ is a $p$-th root of $\zeta^p$.  The valuation of $1-\alpha$ in
$L$ is $p^{n-1}$.  Therefore
$$
w(\sigma(\zeta)-\zeta)=w(\alpha\zeta-\zeta)=w(\alpha-1)=p^{n-1},\quad
i_H(\sigma)=p^{n-1}, 
$$
so $\sigma\in H_{p^{n-1}-1}$ but $\sigma\notin H_{p^{n-1}}$, and the
ramification break occurs at $p^{n-1}-1$.  This allows one to compute the
exponent of the different (and of the discriminant) of $L|K(\zeta^p)$ as
$(p-1).p^{n-1}$.  Applying the transitivity formula for the discriminant to
the tower $L|K(\zeta^p)|K$, we get
$$
np^n-(n+1)p^{n-1}=(p-1).p^{n-1}+p.((n-1)p^{n-1}-np^{n-2})
$$
which is indeed true, and provides an alternative computation of the
discriminant of $L|K$ by induction on~$n$.

\smallskip

% \th PROPOSITION 16
% \enonce
% The ramification filtration in the upper numbering on\/ $G$ is given by
% $$
% G^m=G(m),\quad (m\in[0,n]).
% $$
% \endth
% Let $\psi=\psi_{L|K}$~; we have $G^v=G_{\psi(v)}$ for every
% $v\in[-1,+\infty[$.  The integral formula for $\psi$ (Chapter~16) gives,
% successively, for every integer $r\in[0,n[$,
% $$
% \psi(1)=p-1,\quad
% \psi(r+1)=\psi(r)+p^{r}(p-1),\ 
% $$
% and hence $\psi(m)=p^m-1$ for $m\in[0,n]$.  Therefore, by
% Proposition~13, 
% $$
% G^m=G_{\psi(m)}=G_{p^m-1}=G(m).
% $$ 

\th PROPOSITION 16
\enonce
The breaks in the ramification filtration in the upper numbering on\/ $G$
occur at the integers\/ $1,\ldots,n-1$, and also at\/ $0$ if\/ $p\neq2$.
Moreover, for every integer\/ $m\in[0,n]$, we have\/ $G^m=G(m)$.
\endth
Let $\varphi=\varphi_{L|K}$~; we have $G^{\varphi(u)}=G_u$ for every
$u\in[-1,+\infty[$, so the upper ramification breaks occur at $\varphi(p^m-1)$
for $m\in[0,n[$ (Proposition~13), but not at $0$ if $p=2$, for then
$G(0)=G(1)$.   Denoting by $g_i$ the order of $G_i$, the definition of
$\varphi$ (Chapter~16) gives, for $m\in[0,n]$, 
$$
g_0.\varphi(p^m-1)=
g_1+g_2+\cdots+g_{p^m-1}
=m.g_0
$$
and hence $\varphi(p^m-1)=m$, which also gives $G^m=G_{p^m-1}=G(m)$, as
desired. 

\medbreak
\centerline{\hbox to5cm{\hrulefill}}
\bigbreak

\vfill\eject

\centerline{\bf Lecture 18}
\bigskip
\centerline{\it The case of local nuber fields}

\bigskip

In this lecture we want to determine the maximal unramified extension $\tilde
K$, and the maximal tamely ramified extension $K'$, of a finite extension $K$
of $\Qp$ in a given algebraic closure $\bar K$ of $K$.
% Denote by $d=[K:\Q_p]$ the degree of $K$, by $e=(v(K^\times):v(\Qp^\times))$
% the ramification index, and by $f=[k:\F_p]$ the residual degree~; we have
% $d=ef$ and $q=p^f$, where $q=\Card k$.

Recall that the residue field $\bar k$ of $\bar K$ is an algebraic closure of
$k$ (Lecture~13, Proposition~2).  Recall also that, for each $m\in\N^*$, there
is a unique (separable) extension $k_m$ of $k$ in $\bar k$ of degree~$m$, that
$k_m^\times$ is cyclic of order $q^m-1$, that $k_m|k$ is cyclic, and in fact
$1\mapsto\varphi_m$, where $\varphi_m$ is the $k$-automorphism $x\mapsto x^q$
of $k_m$, is an isomorphism $\Z/m\Z\to\Gal(k_m|k)$.  

For $m'\in\N^*$, the relation ``\thinspace$m'|m$\thinspace'' is equivalent to
``\thinspace$k_{m'}\subset k_m$\thinspace''~; when this is so, the projection
$\Gal(k_m|k)\to\Gal(k_{m'}|k)$ sends $\varphi_m$ to $\varphi_{m'}$ and gets
identified with the natural projection $\pi_{m,m'}:\Z/m\Z\to\Z/m'\Z$.
Finally, $\bar k$ is the union, for $m\in\N^*$, of the $k_m$, and $\Gal(\bar
k|k)$ gets identified with $\hat\Z$, the inverse limit of the system
$(\Z/m\Z,\pi_{m,m'})$~; under this identification, the generator $1$ of the
profinite group $\hat\Z$ goes to the $k$-automorphism $\varphi:x\mapsto x^q$
of $\bar k$.

If follows (Lecture~12, Proposition~6) that, for every $m\in\N^*$, there is a
unique unramified extension $K_m$ of $K$ in $\bar K$~; it is cyclic, and in
fact the reduction map $\Gal(K_m|K)\to\Gal(k_m|k)$ is an isomorphism
(Lecture~14, Proposition~5).  For $m'\in\N^*$, the relation
``\thinspace$m'|m$\thinspace'' is equivalent to ``\thinspace$K_{m'}\subset
K_m$\thinspace''.  The maximal unramified extension $\tilde K$ of $K$ in $\bar
K$ is the union, for $m\in\N^*$, of the $K_m$, and $\Gal(\tilde
K|K)\to\Gal(\bar k|k)$ is an isomorphism.

For every divisor $n$ of $q^m-1$, the polynomial $T^n-1$ has $n$ simple roots
in $k_m$, so it has $n$ simple roots in $K_m$ (Lecture~6, Corollary~2).

\th PROPOSITION 1
\enonce
Let\/ $\zeta\in\bar K^\times$ be an element of order\/ $n$ prime to\/ $p$.
Then\/ $K(\zeta)=K_g$, where $g$ is the order of\/ $q$ in the group\/
$(\Z/n\Z)^\times$.  
\endth
As $n|q^g-1$, we have $K(\zeta)\subset K_g$ and hence $K(\zeta)=K_{g'}$ for
some $g'|g$. (In particular, $K(\zeta)$ is unramified over $K$). As
$K_{g'}^\times$ contains an element --- $\zeta$, for example --- of order $n$
prime to $p$, we have $n|q^{g'}-1$, which means that $q^{g'}\equiv1\pmod n$
and $g|g'$, because $g$ is the order of $q$ in $(\Z/n\Z)^\times$.  Thus
$g'=g$.

\th COROLLARY 2
\enonce
The maximal unramified extension\/ $\tilde K$ of\/ $K$ in\/ $\bar K$ is
obtained by adjoining to\/ $K$ the\/ $n$-th roots of\/ $1$ for every\/ $n$
prime to\/ $p$.
\endth

\th COROLLARY 3
\enonce
If\/ $K$ contains the roots of\/ $1$ of order\/ $n$ prime to\/ $p$, then
$n\,|\,q-1$. 
\endth
Indeed, the order of $q$ in $(\Z/n\Z)^\times$ is $1$ (Proposition~1), meaning
$q\equiv1\pmod n$.  Applied to $K_m$, we see that if $K_m^\times$ has an
element of order $n$ prime to $p$, then $n|q^m-1$.

\medbreak

Let us now determine the maximal tamely ramified extension $K'$ of $K$ in
$\bar K$.  Fix a uniformiser $\pi$ of $K$.  For every $n\in\N^*$ prime to $p$,
choose an $n$-th root $\pi_n=\root n\of\pi$ of $\pi$ in $\bar K$ such that
$\pi_{nd}^d=\pi_n$ for every $d\in\N^*$ (prime to $p$).

For every $n\in\N^*$ prime to $p$, the polynomial $T^n-\pi$ is irreducible
over $K$ (Lecture~8, Lemma~7)~; it has the root $\pi_n$ in $\bar K$, and the
extension $Z_n=K(\pi_n)$ is (totally and) tamely ramified of degree~$n$
(Lecture~13, Proposition~6)~; for $n'\in\N^*$, the relation $n'|n$ is
equivalent to $Z_{n'}\subset Z_n$.

Let $Z_\infty$ be the union, in $\bar K$, of the $Z_n$, for $n\in\N^*$ prime
to $p$.  The extension $Z_\infty|K$ is totally but tamely ramified.

\th PROPOSITION 4
\enonce
The compositum\/ $K'=\tilde KZ_\infty=\tilde K(\pi_n)_{\gcd(p,n)=1}$ is the
maximal tamely ramified extension of\/ $K$ in\/ $\bar K$.
\endth
We have to show that every finite tamely ramified extension $L$ of $K$ in
$\bar K$ is contained in $K'$.  Let $e$ be the ramification index of $L|K$~;
it is prime to $p$ (Lecture~13, Definition~5).  We may write $L=L_0(\root
e\of\varpi)$, where $L_0$ is the maximal unramified subextension of $L|K$, 
and $\varpi$ is a suitable uniformiser of $L_0$ (Lecture~13, Corollary~9) 

As we have seen above, we have $L_0=K_f$, where $f$ is the residual degree of
$L|K$. Also, $\pi$ is a uniformiser of $L_0$ (Lecture~13, Proposition~4), so
write $\varpi=u\pi$, for some unit $u$ of $L_0$.  As $e$ is prime to $p$, the
extension $L_0(\root e\of u)$ is unramified over $L_0$, and hence over $K$.
So is the extension $L(\root e\of1)$ (Proposition~1).  Therefore
$L_0(\zeta,\root e\of u)\subset\tilde K$.

As $\root e\of\varpi=\zeta\root e\of u\root e\of\pi$ for some $e$-th root
$\zeta$ of~$1$, it is clear that $Z_e(\zeta,\root e\of u)$ contains $L$ and is
contained in $\tilde K(\pi_e)=\tilde KZ_e$, hence $L\subset K'=\tilde
KZ_\infty$, as claimed.

\th COROLLARY 5
\enonce
The maximal tamely ramified extension\/ $K'$ of\/ $K$ in\/ $\bar K$ is
obtained by adjoining\/ $\root n\of1$ and\/ $\root n\of\pi$ for every\/
$n\in\N^*$ prime to\/ $p$. 
\endth

\medskip

In the exact sequence $1\to\Gal(K'|\tilde K)\to\Gal(K'|K)\to\Gal(\tilde
K|K)\to1$, we have seen that the quotient $\Gal(\tilde K|K)$ is canonically
isomorphic to the procyclic group $\hat\Z$.  As for the kernel
$G_0=\Gal(K'|\tilde K)$, the inertia subgroup for the extension $K'|K$, one
can see that it is canonically isomorphic to $\Hom(\Z[1/l]_{l\neq p}/\Z,\bar
k^\times)$ (Lecture~14, Proposition~8), where $l$ runs through primes $\neq
p$~; it is isomorphic to the procyclic group $\prod_{l\neq p}\Z_l$.

\medskip

\th PROPOSITION 6
\enonce
Let\/ $u\in\ogoth^\times$ be a unit of\/ $K$, and let\/ $e\in\N^*$ be an
integer prime to~$p$.  The extensions\/ $L=K(\root e\of\pi)$ and\/
$M=K(\root e\of{u\pi})$ are\/ $K$-isomorphic if and only if\/
$u\in\ogoth^{\times e}$. 
\endth
It is clear that $L$ and $M$ are $K$-isomorphic if $u\in\ogoth^{\times e}$.
Conversely, suppose that $L$ and $M$ are $K$-isomorphic, and let $\eta$ be a
unit of $L=M$ such that $\eta\root e\of\pi=\root e\of{u\pi}$.  Raising to the
$e$-th power, we get $\eta^e=u$.  As the extension $K(\eta)$ is unramified
(Lecture~13, Proposition~6), and $L=M$ totally ramified, we have $K(\eta)=K$,
and hence $u\in\ogoth^{\times e}$.

\th COROLLARY 7
\enonce
There are exactly\/ $\gcd(e,q-1)$ totally ramified extensions of\/ $K$ of
degree\/ $e$ (prime to\/ $p$).
\endth
In view of Proposition~5, it suffices to prove that the group
$\ogoth^\times\!/\ogoth^{\times e}$ is (cyclic) of order $\gcd(e,q-1)$. As
$k^\times$ is cyclic of order $q$, this follows from the the
following proposition.

\th PROPOSITION 8 
\enonce
The exact sequence\/ $1\to U_1\to\ogoth^\times\to k^\times\to1$ has a
canonical  splitting.  Moreover, the map\/ $(\ )^e:U_1\to U_1$ is an
isomorphism. 
\endth 
For $x\in k^\times$, let $y\in\ogoth^\times$ be any lift, so that $\bar y=x$.
The limit $\omega(x)=\lim_{n\to+\infty}y^{q^n}$ depends only on $x$, not on
the choice of the lift $y$, and $\omega$ is a section (``Teichm{\"u}ller'') of
the projection $\ogoth^\times\to k^\times$~; we have
$\ogoth^\times=\omega(k^\times)\times U_1$.

Denote by $U_n$ the kernel of $\ogoth^\times\to(\ogoth/\pgoth^n)^\times$.  It
is easily seen that the group $U_n/U_{n+1}$ is isomorphic to the additive
group $k$, a finite commutative $p$-group.  Also, $U_1$ is the projective
limit of the $\Z_p$-modules $U_1/U_n$ under the natural maps $U_1/U_{n+1}\to
U_1/U_n$ ($n\in\N$), so it is a $\Z_p$-module.  As $e$ is invertible in $\Z_p$
by hypothesis, the map $(\ )^e:U_1\to U_1$ is bijective.

Identify $k^\times$ with $\omega(k^\times)\subset\ogoth^\times$, and recall
that $\pi$ is a uniformiser of $K$.  If $u\in k^\times\subset\ogoth^\times$
generates this group, so that its image generates $k^\times\!/k^{\times e}$,
which is cyclic of order $g=\gcd(e,q-1)$, we have seen that the $g$ totally
tamely ramified extensions of $K$ of degree $e$ (prime to $p$) are $K(\root
e\of{u^r\pi})$, with $r\in[0,g[$.

Applying this to $K_f$, and fixing a generator $u\in k_f^\times\subset
K_f^\times$, it follows that every tamely ramified extension $L|K$ of residual
degree\/ $f$ and ramification index\/ $e$ (prime to\/ $p$) can be written\/
$L=K_f(\root e\of{u^r\pi})$, where $K_f=L_0$ is the maximal unramified
subextension and $r\in[0,g[$, where\/ $g=\gcd(e,q^f-1)$.

\th PROPOSITION 9
\enonce
For\/ $L=K_f(\root e\of{u^r\pi})$ to be galoisian over $K$, it is necessary
and sufficient that\/ $e\,|\,q^f-1$ (in which case\/ $g=e$) and\/
$e\,|\,r(q-1)$.  
\endth

If $L|K_f$ is galoisian, then $K_f$ contains a primitive $e$-th root $\zeta$
of $1$, and hence $e\,|\,q^f-1$ (Corollary~3).  Indeed, the conjugates of
$\alpha=\root e\of{u^r\pi}$ are $\zeta^i\alpha$ ($i\in\Z/e\Z$), and the group
$\Gal(L|K_f)$ is generated by $\tau:\alpha\mapsto\zeta\alpha$, so $\zeta$ is
fixed by $\tau$ and hence $\zeta$ is in $K_f$.

The extension $K_f|K$ is galoisian and the order-$\!f$ cyclic group
$\Gal(K_f|K)$ is generated by the automorphism $\varphi$ which induces the map
$x\mapsto x^q$ on the residual extension $k_f|k$.

If $L|K$ is galoisian, then so is $L|K_f$, and $e\,|\,q^f-1$.  For every
$K$-automorphism $\sigma$ of $L$ which restricts to $\varphi$ on $K_f$
--- they are $e$ in number --- we have $\sigma(\alpha)=\root
e\of{u^{qr}\pi}$, where $\alpha=\root e\of{u^r\pi}$, because $\varphi(u)=u^q$.
It follows that $u^{(q-1)r}\in l_0^{\times e}$ (Proposition~6), which is
equivalent to $e\,|\,(q-1)r$.

Notice that, when $e\,|\,q^f-1$, the requirement $u^{(q-1)r}\in l_0^{\times
  e}$ is also equivalent to $(u^r)^{{q^f-1\over e}}$ being a $(q-1)$-th root
of $1$, so being in $k^\times$.

The argument can be reversed.  Namely, if\/ $e\,|\,q^f-1$, then
$k_f^\times\subset K_f^\times$ has an element of order $e$ and the extension
$L|K_f$ is galoisian, indeed cyclic of degree~$e$.  If moreover\/
$e\,|\,r(q-1)$, then $L$ contains all the $K$-conjugates of $\root
e\of{u^r\pi}$, and hence $L|K$ is galoisian.

\th COROLLARY 10
\enonce
If\/ $L=K_f(\root e\of{u^r\pi})$ if galoisian over\/ $K$, then $\Gal(L|K)$
admits the presentation\/ 
$$
\langle
\sigma,\tau\,|\,
\tau^e=1,\sigma^f=\tau^r,\sigma\tau\sigma^{-1}=\tau^q
\rangle. 
$$
\endth
Write $(q-1)r=ne$ for some $n\in\N$ (Proposition~9), so that
$(u^r)^{q-1}=(u^n)^e$, and let $\sigma$ be the extension of the canonical
generator (``Frobenius'') $\varphi$ of $\Gal(K_f|K)$ to $L$ such that
$\sigma(\root e\of{u^r\pi})=u^n\root e\of{u^r\pi}$.  Also, let $\tau$ be the
$K_f$-automorphism of $L$ such that $\tau(\root e\of{u^r\pi})=u^{{q^f-1\over
    e}}\root e\of{u^r\pi}$.  It can be verified that $\sigma,\tau$ generate
$\Gal(L|K)$ subject only to the stated relations
$\tau^e=1$, $\sigma^f=\tau^r$, $\sigma\tau\sigma^{-1}=\tau^q$. 

\th COROLLARY 11
\enonce
For\/ $L=K_f(\root e\of{u^r\pi})$ to be abelian, it is necessary and
sufficient that $e\,|\,q-1$. 
\endth
Indeed, the group defined by the above presentation is commutative if and only
if $\tau^q=\tau$, which is equivalent to $e\,|\,q-1$.  

Notice that the subgroup generated by $\tau$ is the inertia subgroup
$\Gal(L|K_f)$ of $\Gal(L|K)$.

\th COROLLARY 12
\enonce
The only totally tamely ramified abelian extensions of\/ $K$ are $K(\root
q-1\of{\rho\pi})$, as $\rho$ runs through\/ $k^\times\subset K^\times$, and
their subextensions. 
\endth

{\it Exercise}\pointir Show that $\Q_p(\root p-1\of{-p})$ contains a primitive
$p$-th root of $1$.

\th COROLLARY 13
\enonce
The maximal tamely ramified abelian extension of\/ $K$ is\/ $\tilde
K(\root{q-1}\of\pi)$, where\/ $\tilde K$ is the maximal unramified extension
of\/ $K$.
\endth

\medskip

Let us compute the number of extensions of $K$ with given residual degree $f$
and given ramification index $e$ (prime to $p$), in terms of
$g=\gcd(e,q^f-1)$.  For\/ a divisor $t$ of $g$, let $\chi_q(t)$ denote the
order of $q$ in $(\Z/t\Z)^\times$, a group of order $\phi(t)$~; we have
$\chi_q(t)\,|\,\phi(t)$.

\th PROPOSITION 14
\enonce
The number of extensions\/ $L$ of\/ $K$ of residual degree\/ $f$ and
ramification index\/ $e$ (prime to\/ $p$) is
$$
\sum_{t|g}{\phi(t)\over \chi_q(t)},
$$
where $g=\gcd(e,q^f-1)$.  If\/ $g=e$, then precisely\/ $g_1=\gcd(e,q-1)$ of
these are galoisian over\/ $K$.  If\/ $g_1=e$ (in which case\/ $g=e$), then
all of them are abelian.  These are the only galoisian or abelian cases.  
\endth

We have seen that, once we fix a uniformiser $\pi$ of $K$, the totally
ramified degree-$e$ extensions of $K_f$ are indexed by the group
$k_f^\times/k_f^{\times e}$, which is cyclic of order $g$.  To $a\in
k_f^\times/k_f^{\times e}$ corresponds the extension $L_a=K_f(\root
e\of{\alpha\pi})$, where $\alpha$ is any unit of $K_f$ with image $a$~; we
agree to write $L_a=K_f(\root e\of{a\pi})$.

When are two such extensions $K$-isomorphic~? Let $\theta:L_a\to L_b$ be a
$K$-isomorphism, and let $\sigma\in\Gal(K_f|K)$ be its restriction to $K_f$.
As $\theta(\root e\of{a\pi})$ is an $e$-th root of $\sigma(a)\pi$, we have
$b=\sigma(a)$.  Conversely, if $b=\sigma(a)$ for some $\sigma\in\Gal(K_f|K)$,
then there is a $K$-isomorphism $L_a\to L_b$ whose restriction to $K_f$ is
$\sigma$ and which sends $\root e\of{a\pi}$ to $\root e\of{b\pi}$.  So the
  $K$-isomorphism classes of the $L_a$ are indexed by the orbits of
  $k_f^\times/k_f^{\times e}$ under the action of $\Gal(K_f|K)$.
  
  The result will now follow from the following lemma, Proposition~9 and
  Corollary~11.

\th LEMMA 15
\enonce
Let\/ $G$ be a cyclic group of order\/ $g$, written multiplicatively.  Let\/
$q>0$ be an integer prime to\/ $g$, and make\/ $\Z$ act on\/ $G$ by\/
$1\mapsto(x\mapsto x^q)$.  Then the number of orbits is\/
$\sum_{t|g}\phi(t)/\chi_q(t)$. 
\endth
Notice that if $x$, $y$ are in the same orbit, then they have the same order
in $G$.  The possible orders are the divisors of $g$~; for each divisor $t$ of
$g$, there are $\phi(t)$ elements of order $t$.  It is therefore sufficient to
show that the orbit of an order-$t$ element has $\chi_q(t)$ elements.

Indeed, if $x$ has order $t$ in $G$, and if its orbit consists of the $n$
elements \hbox{$x$, $x^q$, $\ldots$, $x^{q^{n-1}}$}, then $n$ is the smallest
integer $>0$ such that $x^{q^n}=x$, which is the same as saying that $n$ is
the smallest integer $>0$ such that $q^n-1\equiv0\pmod t$, so $n=\chi_q(t)$
equals the order of $q$ in $(\Z/t\Z)^\times$.

We have seen (Proposition~9) that for $L_a|K$ to be galoisian, it is necessary
and sufficient that $g=e$ and $a\in k^\times\!/k^{\times e}$, a group of order
$g_1=\gcd(e,q-1)$.  As the action of $\Gal(K_f|K)$ on $k^\times\!/k^{\times
  e}$ is trivial --- this also follows from the fact that $q\equiv1\pmod{g_1}$
---, the number of orbits (for the trivial action of $\Gal(K_f|K)$ on
$k^\times\!/k^{\times e}$) is $g_1=\sum_{t|g_1}\phi(t)$, which is therefore
the number of galoisian extensions among the $L_a$ (when $g=e$).

Finally, if $g_1=e$, then every $L_a|K$ is abelian, and, if $L_a|K$ is abelian
for some $a\in k_f^\times\!/k_f^{\times e}$, then $e\,|\,q-1$ (Corollary~11).

\vfill\eject

\centerline{\bf Lecture 19}
\bigskip
\centerline{\it The maximal abelian extension of\/ $\Q_p$}

\bigskip

In this final lecture, we want to prove that the maximal abelian extension of
$\Q_p$ is obtained by adjoining $\root n\of1$ for every $n\in\N^*$ (local
``Kronecker-Weber'').  We begin with a series of lemmas.  

Let $\zeta_n$ be a primitive $p$-th root of $1$, for every $n\in\N^*$, and put
$\xi_n=\zeta_{p^n}$.  The word {\it extension\/} means an extension of $\Q_p$,
barring explicit mention of some other base field.  We will say that an
extension is {\it cyclotomic\/} if every finite subextension is contained in
$L(\xi_n)$ for suitable~$n$, where $L$ is an unramified extension.  We have to
show that every abelian extension is cyclotomic.

We denote by $L_n=\Q_p(\zeta_{p^{p^n}-1})$ the unramified degree-$p^n$
extension.

\th LEMMA 1 
\enonce 
Suppose that\/ $p=2$ and\/ $n>1$.  The extension\/ $\Q_2(\xi_n)$ is the
compositum of\/ $\Q_2(\xi_2)$ and a totally ramified cyclic extension\/
$M_{n-2}$ of degree\/ $2^{n-2}$.  
\endth 
In view of Lecture~17, Proposition~12, it is sufficient to prove that the
group $(\Z/2^n\Z)^\times$ is the direct product of its quotient
$(\Z/2^2\Z)^\times$ and a cyclic group of order $2^{n-2}$.  Consider the exact
sequence
$$
1\to U_2/U_n\to(\Z_2/2^n\Z_2)^\times\to(\Z_2/2^2\Z_2)^\times\to1~;
$$ 
it clearly has the section $-1\mapsto-1$.  The $\Z_2$-module $U_2=1+2^2\Z_2$
is free ($1+2^2$ is a basis, for example) and the map $(\ )^{2^{n-2}}:U_2\to
U_n$ is an isomorphism (Lecture~8, Corollary~9), so the the group $U_2/U_n$,
being isomorphic to $\Z_2/2^{n-2}\Z_2=\Z/2^{n-2}\Z$, is cyclic of order
$2^{n-2}$.  

{\it Exercise}\pointir Identify $M_1$ among the seven quadratic extensions of
$\Q_2$.  Solution~: clearly, $M_1$ is the subextension of $\Q_2(\xi_3)$ fixed
by the involution $\sigma_{-1}:\xi_3\mapsto\xi_3^{-1}$ (just as $\Q_2(\xi_2)$
is fixed by $\sigma_{5}:\xi_3\mapsto\xi_3^{5}$), so $\xi_3+\xi_3^{-1}\in M_1$.
As $(\xi_3+\xi_3^{-1})^2=2$, we have $M_1=\Q_2(\sqrt2)$.  Notice that
$\Q_2(\xi_3)$ also contains $\Q_2(\sqrt{-2})$.

\th LEMMA 2
\enonce
Suppose that\/ $p$ is odd and\/ $n\in\N^*$.  The extension\/ $\Q_p(\xi_n)$ is
the compositum of\/ $\Q_p(\xi_1)$ and a totally ramified cyclic extension\/
$M_{n-1}$ of degree\/ $p^{n-1}$.
\endth
As in Lemma~1, it is sufficient to prove that $(\Z/p^n\Z)^\times$ is the
direct product of its quotient $(\Z/p\Z)^\times$ and a cyclic group of order
$p^{n-1}$.  We have seen that the exact sequence
$$
1\to U_1/U_n\to(\Z_p/p^n\Z_p)^\times\to(\Z_p/p\Z_p)^\times\to1,
$$ 
has the section $\omega$ (``Teichm{\"u}ller'').  The $\Z_p$-module
$U_1=1+p\Z_p$ is free ($1+p$ is a basis, for example) and the map $(\ 
)^{p^{n-1}}:U_1\to U_n$ is an isomorphism (Lecture~8, Corollary~9), so the
result follows from the fact that the group $\Z_p/p^{n-1}\Z_p=\Z/p^{n-1}\Z$ is
cyclic of order $p^{n-1}$.

\th LEMMA 3
\enonce
Every tamely ramified abelian extension\/ $K$ (of\/ $\Q_p$) is cyclotomic.
\endth
Indeed, $K|\Q_p$, being tame and abelian, is contained in the maximal tamely
ramified abelian extension $\tilde\Q_p(\root{p-1}\of p)$ (Lecture~18,
Corollary~13).  The latter, being the same as $\tilde\Q_p(\xi_1)$ (exercise),
is cyclotomic.  (Recall (Lecture~18, Corollary~2) that the maximal unramified
extension $\tilde\Q_p$ is obtained upon adjoining $\root n\of1$ for every
$n\in\N^*$ prime to $p$.)   

\medbreak

We need a lemma of a general nature.  Let $F$ be a field in which $p$ is
invertible, $\zeta$ a primitive $p$-th root of $1$, and let $E|F$ be a cyclic
extension of degree $p$.  The extension $E(\zeta)|F(\zeta)$ is also cyclic of
degree~$p$~; it corresponds therefore (``Kummer theory'') to an $\F_p$-line
$D\subset F(\zeta)^\times\!/F(\zeta)^{\times p}$.  The group
$G=\Gal(F(\zeta)|F)$ acts on the latter space.  Let $\chi:G\to\F_p^\times$ be
the cyclotomic character giving the action of $G$ on the $p$-th roots of $1$,
so that $\sigma(\zeta)=\zeta^{\chi(\sigma)}$ for every $\sigma\in G$.

\th LEMMA 4
\enonce
The\/ $\F_p$-line\/ $D$ is\/ $G$-stable, and\/ $G$ acts on\/ $D$ via\/
$\chi$. 
\endth
Let $ a\in F(\zeta)^\times$ be such that its image $\bar a$ modulo
$F(\zeta)^{\times p}$ generates $D$, so that $E(\zeta)=F(\zeta)(\root
p\of a)$.  We have $\sigma(\bar a)=\overline{\sigma( a)}$ for
every $\sigma\in G$.  We have to first show that $\sigma(\bar a)\in D$.
Identify $G$ with $\Gal(E(\zeta)|E)$.  For every $\sigma\in G$, we have 
$$
F(\zeta)(\root p\of a)=
F(\zeta,\root p\of a)=
F(\sigma(\zeta),\sigma(\root p\of a))=
F(\zeta)(\sigma(\root p\of a))
$$
and $(\sigma(\root p\of a))^p=\sigma( a)$ is in $F(\zeta)^\times$,
so $\bar a$ and $\overline{\sigma( a)}$ belong to the same
$\F_p$-line, namely $D$.  Hence $D$ is\/ $G$-stable.

Let $\eta:G\to\F_p^\times$ be the character through with $G$ acts on $D$, so
that, for a generator $\sigma$ of $G$, we have $\sigma(\root p\of a)=b(\root
p\of a)^{\eta(\sigma)}$ for some $b\in F(\zeta)^\times$.  Let $\tau$ be the
generator $\root p\of a\mapsto\zeta\root p\of a$ of the group
$\Gal(E(\zeta)|F(\zeta))$, so that $\tau(\zeta)=\zeta$.  We have
$\sigma(\tau(\root p\of a))=\sigma(\zeta\root p\of
a)=\zeta^{\chi(\sigma)}b(\root p\of a)^{\eta(\sigma)}$.

Also, $\tau(\sigma(\root p\of a)) =\tau(b(\root p\of a)^{\eta(\sigma)})
=b\zeta^{\eta(\sigma)}(\root p\of a)^{\eta(\sigma)}$.  But
$\sigma\tau=\tau\sigma$, hence $\eta=\chi$.

\th COROLLARY 5
\enonce
Let\/ $H$ be the\/ $\chi$-eigenspace for the action of\/ $G$ on\/
$F(\zeta)^\times\!/F(\zeta)^{\times p}$.  The map\/ $E\mapsto D$ is an
injection of the set of degree-$p$ cyclic extensions\/ $E|F$ into the set of\/
$\F_p$-lines in\/ $H$.
\endth

{\it Exercise}\pointir Conversely, show that the above map is surjective.  In
other words, every $G$-stable $\F_p$-line $D\subset
F(\zeta)^\times\!/F(\zeta)^{\times p}$ on which $G$ acts via $\chi$ comes from
a degree-$p$ cyclic extension $E|F$.

\medbreak

\th LEMMA 6
\enonce
Suppose that\/ $p$ is odd.  Every degree-$p$ cyclic extension of\/ $\Q_p$ is
cyclotomic, and indeed contained in $L_1M_1$. 
\endth
The idea is to show, upon taking $F=\Q_p$, that the sub-$\F_p$-space $H$
(Corollary~5) of $F(\xi_1)^\times\!/F(\xi_1)^{\times p}$ on which $G$ acts via
$\chi$ is $2$-dimensional, so that $F$ has at most $p+1$ degree-$p$ cyclic
extensions.  Recall that $L_1$ is the unramified degree-$p$ extension, and
$M_1$ is the totally ramified degree-$p$ (cyclic) extension contained in
$F(\xi_2)$ (Lemma~2).  The compositum $L_1M_1$ contains $p+1$ degree-$p$
cyclic extensions, and all of them are cyclotomic, for they are all contained
in $L_1(\xi_2)$.  Hence every degree-$p$ cyclic extension is cyclotomic.

It remains to show that the $\F_p$-space $H$ is $2$-dimensional.  Let $D$ be
the line generated by $\bar\xi_1$ and $\bar U_p$ the line generated by
$\overline{1+\pi_1^p}$, where $\pi_1=1-\xi_1$.  It can be shown that they are
distinct and that $H$ is the plane $D\bar U_p$, that $M_1$ corresponds to $D$
and $L_1$ to $\bar U_p$.  For some hints, see {\tt arXiv\string:0711.3878},
Part~III.

Suppose first that\/ $p\neq2$.  
\th LEMMA 7
\enonce
Two elements suffice to generate the group\/
$G=\Gal(K|\Q_p)$ of automorphisms of an abelian\/ $p$-extension\/ $K|\Q_p$.
\endth
Recall that $\dim_{\F_p}A/pA$ elements suffice to generate an abelian
$p$-group $A$.  Writing $G$ additively, the extension $K^{pG}$ is abelian of
exponent~$p$, and hence $K^{pG}\subset L_1M_1$, the maximal abelian extension
of exponent~$p$ (Lemma~6).  Consequently, $\dim_{\F_p}G/pG$ is at most~$2$.

\th LEMMA 8
\enonce
Every quadratic extension of\/ $\Q_2$ is cyclotomic, and indeed contained in\/
$L_1M_1(\xi_2)$. 
\endth
The group of automorphisms $(\Z/2^3\Z)^\times$ of the totally ramified
extension $\Q_2(\xi_3)$ has exponent $2$, hence the same is true of
$L_1(\xi_3)|\Q_2$, where $L_1$ is the unramified quadratic extension.
Therefore $L_1(\xi_3)$ is contained in $\Q_2(\sqrt5,\sqrt3,\sqrt2)$, 
the maximal abelian extension of exponent~$2$ (Lecture~6, Corollary~6).  As
the two extensions have degree~$2^3$, we get
$L_1(\xi_3)=\Q_2(\sqrt5,\sqrt3,\sqrt2)$, and hence every quadratic extension
is cyclotomic.  Notice that $M_1(\xi_2)=\Q_2(\xi_3)$ (Lemma~1), and
$L_1=\Q_2(\zeta_3)$, so the maximal abelian extension of exponent~$2$ is
$\Q_2(\zeta_{24})$.

\th LEMMA 9
\enonce
Three elements suffice to generate the group\/ $G=\Gal(K|\Q_2)$ of
automorphisms of an abelian\/ $2$-extension\/ $K|\Q_2$.
\endth
As in Lemma~7, $G$ can be generated by $\dim_{\F_2}G/2G$ elements.  The
extension $K^{2G}$ is abelian of exponent~$2$, and hence
$K^{2G}\subset\Q_2(\sqrt5,\sqrt3,\sqrt2)$.  Consequently, $\dim_{\F_2}G/2G$ is
at most~$3$.

\medbreak

\th LEMMA 10
\enonce
Let\/ $F$ be a field in which\/ $2$ is invertible, and let\/ $E|F$ be a cyclic
extension of degree\/~$4$.  Every\/ $d\in F^\times$ such that\/ $d\notin
F^{\times2}$ but\/ $d\in E^{\times2}$ can be written\/ $d=a^2+b^2$ for some\/
  $a,b\in F$.
\endth
Write $F'=F(\delta)$, where $\delta^2=d$, and write $E=F'(\theta)$, with
$\theta^2\in F'$.  Let $\sigma\in\Gal(E|F)$ be a generator, so that
$\sigma(\delta)=-\delta$, $\sigma^2(\theta)=-\theta$.

Put $\theta^2=u+v\delta$ (for some $u,v\in F$), so that
$\sigma(\theta)^2=u-v\delta$.  As
$\sigma(\theta.\sigma(\theta))=-\theta.\sigma(\theta)$, we have
$\theta\sigma(\theta)=s\delta$ for some $s\in F^\times$.  Squaring the last
relation, we get
$$
u^2-dv^2=ds^2,
$$
which shows that $d^{-1}$, and hence $d$, has the required form.

\th LEMMA 11
\enonce
There is no degree-$4$ cyclic extension\/ $K$ of\/ $\Q_2$ containing\/
$\xi_2=\sqrt{-1}$. 
\endth
If there were such an extension, we would be able to write $-1=a^2+b^2$ for
some $a,b\in\Q_2$ (Lemma~10), but this is impossible.

\medbreak

\th THEOREM 12
\enonce
Every abelian extension\/ $K$ of\/ $\Q_p$ is cyclotomic.
\endth
We may suppose that $G=\Gal(K|\Q_p)$ is finite~; write $G=PT$, where $P$ is a
$p$-group and $T$ is a $p'$-group --- a group of order prime to $p$.  As the
$T$-extension $K^P$ is tame, it is cyclotomic (Lemma~3)~; we may therefore
assume that $G$ is a $p$-group.

Assume first that $p\neq2$, and let $p^t$ be the exponent of $G$.  Replacing
$K$ by $KL_tM_t$, where $L_t$ is the degree-$p^t$ unramified extension of
$\Q_p$ and $M_t$ the degree-$p^t$ extension of $\Q_p$ contained in
$\Q_p(\xi_{t+1})$ (Lemma~2), we may assume that $L_tM_t\subset K$.  

As $G$ can be generated by two elements (Lemma~7), has exponent $p^t$, and has
the quotient $\Gal(L_tM_t|\Q_p)$ of type $(p^t,p^t)$, the group $G$ is
necessarily of type $(p^t,p^t)$ and we have $K=L_tM_t$, which is cyclotomic,
for $L_tM_t\subset L_t(\xi_{t+1})$.

Now consider $p=2$, and let $2^t$ be the exponent of $G$.  As in the case of
odd primes, we may suppose that $L_tM_t\subset K$, where $L_t$ is the
degree-$2^t$ unramified extension and $M_t$ the degree-$2^t$ extension
contained in $\Q_2(\xi_{t+2})$ (Lemma~1).  We may also suppose that $\xi_2\in
K$.

As the group $G$ can be generated by three elements (Lemma~9), has exponent
$2^t$, has the quotient $\Gal(L_tM_t(\xi_2)|\Q_2)$ of type $(2^t,2^t,2)$, but
does not have a cyclic quotient $G/H$ of order $2^2$ such that $K^H$ contains
$\xi_2$ (Lemma~11), $G$ is necessarily of type $(2^t,2^t,2)$ and we have
$K=L_tM_t(\xi_2)$, which is cyclotomic, for $L_tM_t(\xi_2)\subset
L_t(\xi_{t+2})$.

\th COROLLARY 13
\enonce
The group of automorphisms of the maximal abelian extension of\/ $\Q_p$ is\/
$G=\Z_p^\times\times\hat\Z$.  The breaks in the ramification filtration\/
$(G^u)_{u\in[-1,+\infty[}$ occur at\/ $u\in\{-1\}\cup\N$, but not at\/ $u=0$
when\/ $p=2$.  We have\/ $G^{-1}=G$, $G^0=\Z_p^\times$, and\/ $G^n=U_n$ for
every integer\/ $n>0$.
\endth
Theorem~12 says that the maximal abelian extension is the compositum
$K_\infty\tilde\Q_p$ of the maximal unramified extension $\tilde\Q_p$ and the
inductive limit $K_\infty$ of the totally ramified extensions
$K_n=\Q_p(\xi_n)$, for $n\in\N$.  The group $\Gal(K_\infty|\Q_p)$ is therefore
the projective limit, for $n\in\N$, of the groups
$\Gal(K_n|\Q_p)=(\Z/p^n\Z)^\times$, namely $\Z_p^\times$.  Also,
$\Gal(\tilde\Q_p|\Q_p)=\hat\Z$ (Lecture~18).  Therefore
$G=\Z_p^\times\times\hat\Z$.

We know that the filtration in the upper numbering is compatible with the
passage to the quotient (Lecture~16, Theorem~11).  We also know that the
breaks in the ramification filtration in the upper numbering on\/
$\Gamma_{(n)}=(\Z/p^n\Z)^\times$ occur at the integers\/ $1,\ldots,n-1$, and
also at\/ $0$ if\/ $p\neq2$~; moreover, for $m\in[0,n]$, we have\/
$\Gamma_{(n)}^m=\Ker(\Gamma_{(n)}\to\Gamma_{(m)})$ (Lecture~17,
Proposition~16).  The result follows from these two facts.

\th COROLLARY 14
\enonce
Let\/ $K|\Q_p$ be an abelian extension, and\/ $H=\Gal(K|\Q_p)$.  If a break in
the ramification filtration\/ $(H^u)_{u\in[-1,+\infty[}$ occurs at\/ $u$, then
$u\in\{-1\}\cup\N$.  
\endth
Indeed, $H$ is a quotient of $G$, and, if the filtration $(H^u)_u$ has a break
at $u$, then a break occurs at $u$ for the upper-numbering filtration on
$G=\Z_p^\times\times\hat\Z$ (Lecture~16, Theorem~11).  Therefore
$u\in\{-1\}\cup\N$ (Corollary~13).

{\it Exercise}\pointir For every quadratic extension $K$ of $\Q_p$, determine
the subgroup $\Gamma\subset\Z_p^\times\times\hat\Z$ such that $M^\Gamma=K$,
where $M$ is the maximal abelian extension.

{\it Exercise}\pointir For every $n>0$, determine the subgroup
$\Gamma\subset\Z_p^\times\times\hat\Z$ such that $M^\Gamma=\Q_p(\xi_n)$, where
$M$ is still the maximal abelian extension.  Answer~: $\Gamma=U_n\times\hat\Z$
(cf.~Lecture~17, Proposition~12).  Also determine the subgroup
$\Delta\subset\Z_p^\times\times\hat\Z$ such that $M^\Delta=M_n$, where $M_n$
is the degree-$p^n$ extension introduced in Lemmas~1 and~2.

\bye